\RequirePackage{snapshot}
\documentclass[12pt]{elsarticle}
\usepackage{graphicx}
\usepackage{amssymb}    
\usepackage{epstopdf}
\usepackage{color}
\usepackage{hyperref}

\usepackage{lscape}
\usepackage{epsfig}
\usepackage{fullpage}
\usepackage[titletoc]{appendix}
\usepackage{fancyhdr}
\usepackage{amsmath,amssymb,xspace}
\usepackage{tabularx}

\newfont{\bbb}{msbm10 scaled\magstep1}

\newtheorem{prop}{Proposition}[section]

\let \leq \leqslant

\let \epsilon \varepsilon

{ \par \medskip \par
  \noindent \textit{\textbf{Demonstration\/}} : }{\null \hfill $\Box$ \par }
 
\newcommand{\R} {\ensuremath{\mathbb{R}}}

\newcommand{\N} {\ensuremath{\mathbb{N}}}
\newcommand{\C} {\ensuremath{\mathbb{C}}}

\makeatletter
\newcommand{\doublewidetilde}[1]{{%
  \mathpalette\double@widetilde{#1}%
}}
\newcommand{\double@widetilde}[2]{%
  \sbox\z@{$\m@th#1\widetilde{#2}$}%
  \ht\z@=.9\ht\z@
  \widetilde{\box\z@}%
}

\begin{document}

\begin{frontmatter}

\title{Pseudospectral computational methods for the time-dependent Dirac equation in static curved spaces}




\author[iecl]{Xavier Antoine}
\ead{xavier.antoine@univ-lorraine.fr}

\author[iqc]{Fran\c{c}ois Fillion-Gourdeau}
\ead{francois.fillion@emt.inrs.ca}

\author[carl,crm]{Emmanuel Lorin}
\ead{elorin@math.carleton.ca}

\author[iqc]{Steve McLean}

\address[iecl]{Institut Elie Cartan de Lorraine, Universit\'e de Lorraine, UMR 7502, Inria Nancy-Grand Est, F-54506 Vandoeuvre-l\`es-Nancy Cedex, France}
\address[iqc]{Institute for Quantum Computing, University of Waterloo, Waterloo, Ontario, Canada, N2L~3G1}
\address[carl]{School of Mathematics and Statistics, Carleton University, Ottawa, Canada, K1S~5B6}
\address[crm]{Centre de Recherches Math\'{e}matiques, Universit\'{e} de Montr\'{e}al, Montr\'{e}al, Canada, H3T~1J4}

\begin{abstract}
Pseudospectral numerical schemes for solving the Dirac equation in general static curved space are derived using a pseudodifferential representation of the Dirac equation along with a simple Fourier-basis technique. Owing to the presence of non-constant coefficients in the curved space Dirac equation, convolution products usually appear when the Fourier transform is performed. The strategy based on pseudodifferential operators allows for efficient computations of these convolution products by using an ordinary fast Fourier transform algorithm. The resulting numerical methods are efficient and have spectral convergence. Simultaneously, wave absorption at the boundary can be easily derived using absorbing layers to cope with some potential negative effects of periodic conditions inherent to spectral methods. The numerical schemes are first tested on simple systems to verify the convergence and are then applied to the dynamics of charge carriers in strained graphene. 
\end{abstract}

\end{frontmatter}

\section{Introduction}

The Dirac equation is one of the most important equations of Physics, giving a quantum description of relativistic spin-1/2 particles such as electrons and quarks \cite{Thaller}. For this reason, it can be found in the theoretical description of many physical systems in nuclear physics, condensed matter physics, laser-matter interaction, cosmology, and many others. However, this partial differential equation is notoriously hard to solve, and one often has to resort to numerical methods for accurate and non-perturbative solutions. Motivated by quantum electrodynamics (heavy ion collision, pair production)
\cite{PhysRevA.78.062711,PhysRevA.24.103,PhysRevC.71.024904,p13,Fillion-Gourdeau2013}, by strong field physics (Schwinger's effect, intense laser-molecule interaction) \cite{Fillion-Gourdeau2013,Fillion-Gourdeau2013b,graph2,p6,keitel4}, graphene modeling \cite{graph1, Katsnelson2006}, tremendous efforts have been put these past two decades on the development of numerical methods for the computation of the time-dependent Dirac equation in flat Minkowski space. Real space methods such as Quantum Lattice Boltzman techniques \cite{PhysRevLett.111.160602,succi,jcp2014,cpc2012,Lorin_Bandrauk}, Galerkin methods \cite{0022-3700-19-20-003,Fillion-Gourdeau2016122,8}, and pseudospectral or spectral methods \cite{16,fft1,fft2,fft3,Beerwerth2015189,Bauke20112454,keitel,grobe} were developed to efficiently and accurately solve the Dirac equation.   Semi-classical regimes were considered in \cite{shijin} using Gaussian beams or Frozen Gaussian Approximations \cite{FGA,rousse}, while the non-relativistic limit has been studied in several recent papers \cite{bao3,bao2}. The computational difficulties for solving the time-dependent Dirac equation include the {\it fermion doubling problem}, related to numerical dispersion \cite{cpc2012}, and the Zitterbewegung \cite{Thaller}, resulting in highly oscillating solutions whose origin can be traced back to the presence of the mass term $\beta mc^2$. Finally, drastic stability conditions can lead to numerical diffusion related to the finite wave propagation speed, the speed of light $c$. 

Another numerical challenge for the Dirac equation shared by any other wave equation in real space solved on a truncated domain is the need of imposing special boundary conditions in order to avoid spurious wave reflections at the computational domain boundary. Therefore, the computational methods on truncated domains require non-reflective boundary conditions \cite{jcp2014,AABES,MOLPHYS,Hammer2014728,Antoine2014268}, absorbing or perfectly matched layers (PML) \cite{MOLPHYS,pinaud,TurkelYefet,zeng,Tsynkov}, or the introduction of an artificial potential \cite{pinaud}. On the other hand, Fourier-based methods applied on bounded domains naturally induce periodic boundary conditions, which can be problematic when dealing with delocalized wave functions. In this respect, the technique developed in \cite{jcp2019,cpc2017} for the Dirac equation in flat space, where a spectral method is combined with PML, is an interesting alternative. One of the goals of this article is to extend this numerical scheme to the Dirac equation in curved space.

Recently, the Dirac equation in curved space-time has gained important interest in some applications such as condensed matter physics for describing the dynamics of charge carriers in deformed 2D Dirac materials \cite{CORTIJO2007293,Cortijo_2007}, as well as astrophysics for fermion tunneling in black holes \cite{Kerner_2008,Di_Criscienzo_2008,LI2008370,CHEN2008106}. In its discrete version, it has been considered from the lattice Boltzmann technique point of view \cite{Succi2015,PhysRevB.98.155419,Debus_2018} and as continuous limits of quantum walks \cite{PhysRevA.88.042301,Arrighi2016,Mallick_2019}. However, the literature on numerical methods for this equation is scarse. Therefore, this article is an attempt to fill this gap. In particular, we present a numerical method based on the pseudodifferential representation \cite{TaylorBook} of the Dirac equation in curved space in combination with perfectly matched layers (PML). In the pseudodifferential representation, it is possible to efficiently use Fourier-based methods, even though the Dirac equation under consideration has non-constant spatial coefficients. A similar methodology was successfully developed in \cite{jcp2019,cpc2017} for the Dirac equation in flat space and in \cite{AntoineGeuzaineTang} for the Gross-Pitaevskii equation. In these two cases, the pseudodifferential representation was used within some absorbing layers at the truncated domain boundary to implement wave absorption at the boundary using PML. Within the domain however, the scheme was an usual spectral numerical method. In curved space, as we will show below, it is beneficial to employ the pseudodifferential representation in the whole domain because convolution products do not appear explicitly. When combined with an implicit scheme for the time evolution, this allows for a Fourier-based method benefiting from spectral convergence and unconditional $\ell^2-$stability. In addition, the PML can be straightforwardly included, reducing the effects of the inherent periodic boundary conditions.

This paper is organized as follows. In Section \ref{sec:dirac}, we describe the Dirac equations under consideration. In Section \ref{sec:PML}, we present the type of absorbing (perfectly matched) layers for the Dirac equation. In Sections \ref{sec:discrete} and \ref{sec:space}, we propose and analyze two numerical methods for approximating the Dirac equation in static curved space-times.  Section \ref{sec:numerics} is dedicated to numerical experiments. We conclude in Section \ref{sec:conclusion}.

 \section{Dirac equation}\label{sec:dirac}
This section is devoted to the presentation of the Dirac equations studied in this paper. We first recall the basics of the usual Dirac equation in flat space and then, we present its extension to curved space. Finally, we reformulate the latter in ``Hamiltonian form'', similar to the one in flat space but with space-dependent coefficients.

\subsection{Dirac equation in flat space}
The time-dependent Dirac equation in Cartesian coordinates reads \cite{Itzykson:1980rh}
\begin{eqnarray}
{\tt i}\partial_t \psi(t,{\boldsymbol x}) = H_{\mathrm{flat}}(t,{\boldsymbol x}) \psi(t,{\boldsymbol x}),
\label{eq:dirac_eq}
\end{eqnarray}
where $\psi(t,{\boldsymbol x})$ is the time and coordinate dependent four-spinor, and $H_{\mathrm{flat}}$ is the Hamiltonian operator. The latter is given by
\begin{eqnarray}
H_{\mathrm{flat}}(t,{\boldsymbol x})  =  
\boldsymbol{\alpha} \cdot \left[  -{\tt i}\nabla - e\boldsymbol{A}(t,{\boldsymbol x}) \right] + \beta m + \mathbb{I}_{4}V(t,{\boldsymbol x}),
\label{eq:hamiltonian_flat}
\end{eqnarray}
where $\psi(t,{\boldsymbol x}) \in L^{2}(\mathbb{R}^{3}) \otimes \mathbb{C}^{4}$ is the time $t=x^{0}$ and coordinate (${\boldsymbol x} = (x^{1},x^{2},x^{3})$) dependent four-spinor, $\boldsymbol{A}(t,{\boldsymbol x})$ represents the three space components of the electromagnetic vector potential, $V(t,{\boldsymbol x}) = eA_{0}(t,{\boldsymbol x})+V_{\textrm{nuc.}}({\boldsymbol x})$ is the sum of the scalar and interaction potentials, $e$ is the electric charge (with $e=-|e|$ for an electron),  $\mathbb{I}_{4}$ is the $4 \times 4$ unit matrix and $\beta,\boldsymbol{\alpha}=(\alpha^{i})_{i=1,2,3}$ are the Dirac matrices. 
In this work, the Dirac representation is used, where
\begin{eqnarray}
\beta = 
\begin{bmatrix}
\mathbb{I}_{2} & 0 \\
0 & -\mathbb{I}_{2} 
\end{bmatrix}
\; \; , \; \;
 \alpha^{i} = 
 \begin{bmatrix}
 0 & \sigma^{i} \\
 \sigma^{i} & 0 
 \end{bmatrix}.
\label{eq:dirac_mat}
\end{eqnarray}
The $\sigma^{i}$ are the usual $2 \times 2$ Pauli matrices defined as
\begin{eqnarray}
\sigma^{1} = 
\begin{bmatrix}
0 & 1 \\ 1 & 0  
\end{bmatrix}
\;\; \mbox{,} \;\;
\sigma^{2} = 
\begin{bmatrix}
0 & -{\tt i} \\ {\tt i} & 0 
\end{bmatrix}
\;\; \mbox{and} \;\;
\sigma^{3} = 
\begin{bmatrix}
1 & 0 \\ 0 & -1 
\end{bmatrix},
\end{eqnarray}
while $\mathbb{I}_{2}$ is the $2 \times 2$ unit matrix. Note that natural units are used where $\hbar = c = 1$.

To simplify the notation further and to parallel the one used in the next section for the Dirac equation in curved space time, it is convenient to write the Hamiltonian as
\begin{eqnarray}
H_{\mathrm{flat}}(t,{\boldsymbol x})  =  
   -{\tt i}\boldsymbol{\alpha} \cdot \nabla + \beta m + F_{\mathrm{flat}}(t, \boldsymbol{x}),
\label{eq:hamiltonian_flat_simp}
\end{eqnarray}
where the function $F_{\mathrm{flat}}$ contains the contribution from the electromagnetic potential:
\begin{eqnarray}
F_{\mathrm{flat}}(t, \boldsymbol{x}):=  - e \boldsymbol{\alpha} \cdot \boldsymbol{A}(t,{\boldsymbol x})   + \mathbb{I}_{4}V(t,{\boldsymbol x}) .
\end{eqnarray} 


\subsection{Dirac equation in curved space}

In this section, the generalization of the Dirac equation to curved space is presented. Throughout, Einstein's notation is assumed with the following conventions: greek indices relate to general  curved space (characterized by the general metric $g^{\mu \nu}(x)$, where $x:=(t,\boldsymbol{x})$ denotes a space-time point), uppercase latin indices relate to flat space (characterized by the Minkowski metric $\eta^{AB} = \mathrm{diag}(1,-1,-1,-1)$)  while lowercase latin indices are summed over spatial coordinates only ($g^{ij}(x)$ for $i,j = 1,2,3$). 

The extension of the Dirac equation to curved space follows by imposing general covariance under arbitrary coordinate transformations. In covariant notation and general background space, the Dirac equation takes the form \cite{pollock2010dirac}
\begin{eqnarray}
\label{eq:dirac_covariant}
\biggl\{ {\tt i} \gamma^{\mu}(x) \left[ \partial_{\mu} + \Omega_{\mu}(x) -{\tt i}eA_{\mu}(x) \right] - m \biggr\} \psi(x) = 0,
\end{eqnarray}
where $A_{\mu}$ is the four vector electromagnetic potential. The generalized gamma matrices used in the Dirac equation define a Clifford algebra:
\begin{eqnarray}
\label{eq:anticomm}
\{\gamma^{\mu}(x), \gamma^{\nu}(x) \} = 2 g^{\mu \nu}(x),
\end{eqnarray}
where the notation $\{\cdot,\cdot \}$ stands for the anticommutator and $g^{\mu \nu}$ is the metric characterizing the curved space. These matrices are a generalization of the Dirac matrices in flat space, which are not space-dependent and are related to the Minkowski metric as
\begin{eqnarray}
\{\gamma^{A}, \gamma^{B} \} = 2 \eta^{AB}.
\end{eqnarray}
The two sets of matrices are related \textit{via} the tetrad formalism as
\begin{eqnarray}
\gamma^{\mu}(x) = \gamma^{A}e_{ A}^{\ \mu}(x),
\end{eqnarray}  
where $e_{ A}^{\ \mu}(x)$ is the tetrad. The tetrads are used to link the metric in curved and flat spaces and thus, obey the property:
\begin{eqnarray}
g^{\mu \nu}(x) = e_{A}^{\  \mu}(x) e_{B}^{\ \nu}(x) \eta^{AB}.
\end{eqnarray}
%

The spinorial affine connection $\Omega_{\mu}(x)$ was introduced in the Dirac equation to preserve the covariance. It is given by
\begin{eqnarray}
\Omega_{\mu}(x) = -\frac{{\tt i}}{4} \omega_{\mu}^{\ AB}(x) \sigma_{AB},
\end{eqnarray}
where $\sigma_{AB} = {\tt i}[\gamma_{A},\gamma_{B}]/2$ is the commutator of the ``flat space'' Dirac matrices while the spin connection is
\begin{eqnarray}
\omega_{\mu}^{\ AB}(x) =e_{\nu}^{\ A}(x) \left[   
\partial_{\mu}e^{\nu B}(x)  
+ \Gamma^{\nu}_{\ \mu \sigma}(x) e^{\sigma B}(x)
\right],
\end{eqnarray}
where the  Christoffel symbols $\Gamma^{\nu}_{\ \mu \sigma}(x)$ were introduced.  It is also important to notice that in curved space, the usual $\ell^2-$norm is not preserved. Instead, denoting
\begin{eqnarray}\label{l2gamma}
\left.
\begin{array}{lclcl}
\langle \psi, \psi\rangle_{\gamma} & = & \|\psi\|_{\gamma}^2 & = &\displaystyle  \int \sqrt{|g(x)|}\psi^{\dagger}[\gamma^0\gamma^0(x)]\psi d^3x \, ,
\end{array}
\right.
\end{eqnarray}
where $g(x)$ is the determinant of the metric, the norm $\|\psi\|_{\gamma}$ is preserved in time, see \cite{Leclerc_2007,PhysRevD.22.1922,PhysRevLett.44.1559}. To develop a numerical scheme, it is convenient to rewrite \eqref{eq:dirac_covariant} in a form similar to \eqref{eq:dirac_eq} in flat space, i.e. in ``Hamiltonian form''. This is performed straightforwardly multiplying \eqref{eq:dirac_covariant} by $\gamma^{0}(x)$ and using the anticommutation relation \eqref{eq:anticomm}. The Dirac equation in curved space can then be written as \cite{PhysRevD.22.1922}
\begin{eqnarray}
{\tt i}\partial_t \psi(t,{\boldsymbol x}) = H(t,{\boldsymbol x}) \psi(t,{\boldsymbol x}),
\label{eq:dirac_eq2}
\end{eqnarray}
where $H(t,{\boldsymbol x})$ is the Dirac Hamiltonian operator in curved space given by
\begin{eqnarray}
H(t,{\boldsymbol x})  &=&-{\tt i} [g^{00}(x)]^{-1} \gamma^{0}(x)  
\gamma^{i}(x)   \left[ \partial_{i} + \Omega_{i}(x) - {\tt i} e A_{i}(x)\right]  \nonumber \\
& &   +[g^{00}(x)]^{-1} \gamma^{0}(x) m - \mathbb{I}_{4} \left[{\tt i} \Omega_{0}(x) + eA_{0}(x) \right]  .
\end{eqnarray}
Defining generalized Dirac matrices as
\begin{eqnarray}
\beta(t,\boldsymbol{x}) &:=& [g^{00}(x)]^{-1} \gamma^{0}(x), \\
\alpha^{i}(t, \boldsymbol{x}) &:=& [g^{00}(x)]^{-1} \gamma^{0}(x)\gamma^{i}(x),
\end{eqnarray}
the Dirac Hamiltonian becomes
\begin{eqnarray}
H(t,{\boldsymbol x})  &=&  
-{\tt i}\boldsymbol{\alpha}(t,{\boldsymbol x}) \cdot  \left[\nabla + \boldsymbol{\Omega}(t,{\boldsymbol x}) - {\tt i} e \boldsymbol{A}(t,{\boldsymbol x})\right] 
 +  \beta (t,{\boldsymbol x}) m  \nonumber \\
& & -\mathbb{I}_{4} \left[ {\tt i}\Omega_{0}(t,{\boldsymbol x}) + eA_{0}(t,{\boldsymbol x}) \right]  .
\label{eq:dirac_hamil}
\end{eqnarray}
At this point, it is convenient to simplify the notation further by defining a function $F$ that allows for writing the Hamiltonian as
\begin{eqnarray}
H(t,{\boldsymbol x})  &=&  
-{\tt i}\boldsymbol{\alpha}(t,{\boldsymbol x}) \cdot \nabla  
+  \beta (t,{\boldsymbol x}) m  + F(t,{\boldsymbol x})   \, ,
\label{eq:dirac_simplified}
\end{eqnarray} 
where 
\begin{eqnarray*}
F(t,{\boldsymbol x}) := 
-{\tt i}\boldsymbol{\alpha}(t,{\boldsymbol x}) \cdot  \left[ \boldsymbol{\Omega}(t,{\boldsymbol x}) - {\tt i} e \boldsymbol{A}(t,{\boldsymbol x})\right] 
- \mathbb{I}_{4} \left[ {\tt i}\Omega_{0}(t,{\boldsymbol x}) + eA_{0}(t,{\boldsymbol x}) \right]  . 
\end{eqnarray*}
This is the general form of the Dirac Hamiltonian in curved space-time. The functions $\beta (t,{\boldsymbol x})$, $\boldsymbol{\alpha}(t,{\boldsymbol x})$ and $F(t,{\boldsymbol x})$ need to be determined \textit{a priori} from the metric and/or from the electromagnetic field potential. Some explicit examples are presented below. It is well-known that the preceding Hamiltonian is not self-adjoint with respect to the covariant inner product when the metric is time-dependent \cite{Leclerc_2007,PhysRevD.22.1922,PhysRevLett.44.1559,PhysRevD.79.024020}, casting some doubts on the conservation of probability. This can be remedied by adding a new term in the Hamiltonian as \cite{Leclerc_2007,PhysRevD.79.024020}
\begin{eqnarray}
H'(t,{\boldsymbol x}) = H(t,{\boldsymbol x}) + \frac{{\tt i}}{2} \partial_{t} \ln \left(\sqrt{|g(x)| g^{00}(x)} \right).
\end{eqnarray} 
The new Hamiltonian $H'$ is now self-adjoint. The new term can be interpreted in many ways, such as in the pseudo-hermitian operator formalism \cite{PhysRevD.79.024020,PhysRevD.82.104056,PhysRevD.83.105002} or as the time-dependence of the position eigenstates \cite{PhysRevD.79.024020}. However, the generality of these results are disputed by other authors, who proposed a different approach \cite{arminjon2010basic}. In this work, we do not dwell into this controversy as it is outside the scope of this article. For the sake of simplicity, and throughout the rest of this article, we will assume the metric is time-independent, allowing us to write the Hamiltonian as
\begin{eqnarray}
H(t,{\boldsymbol x})  &=&  
-{\tt i}\boldsymbol{\alpha}({\boldsymbol x}) \cdot \nabla  
+  \beta ({\boldsymbol x}) m  + F(t,{\boldsymbol x})   , \\
	F(t,{\boldsymbol x}) &=& 
	-{\tt i}\boldsymbol{\alpha}({\boldsymbol x}) \cdot  \left[ \boldsymbol{\Omega}({\boldsymbol x}) - {\tt i} e \boldsymbol{A}(t,{\boldsymbol x})\right] 
	- \mathbb{I}_{4} \left[ {\tt i}\Omega_{0}({\boldsymbol x}) + eA_{0}(t,{\boldsymbol x}) \right]  . 
\end{eqnarray}
This Hamiltonian is now self-adjoint and is the starting point for the development of the numerical schemes. 
The main difference between the flat and curved space versions of the Dirac equation is twofold:  1) the Dirac matrices are space dependent functions $\boldsymbol{\alpha}({\boldsymbol x})$ and $\beta({\boldsymbol x})$, assumed here to be smooth, and 2) the function $F(t,{\boldsymbol x})$ contains the contribution coming from the spin affine connection and the metric. As long as the latter are smooth enough, they do not lead to any particular computational issues. On the other hand, because the Dirac matrices are space dependent, it is obviously not possible to solve directly this equation with a Fourier-based method without avoiding convolution products. However, we can rewrite the equation in pseudodifferential form as follows in $\R^3$
\begin{eqnarray}
{\tt i}\partial_t \psi(t,{\boldsymbol x})  =   -{\tt i}\boldsymbol{\alpha}({\boldsymbol x}) \cdot \mathcal{F}^{-1}_{{\boldsymbol x}} \big\{{\tt i} \boldsymbol{\xi} \mathcal{F}_{{\boldsymbol x}}\{ \psi \} (t,\boldsymbol{\xi}) \big\} \nonumber 
+ \biggl\{ \beta({\boldsymbol x}) m  
+ F(t,\boldsymbol{x})  \biggr\}\psi(t,{\boldsymbol x}),
\label{eq:dirac_eq3}
\end{eqnarray}
which will be the ground of the proposed methodology. In this last equation,  $ \mathcal{F}_{{\boldsymbol x}}\{ \cdot \} (t,\boldsymbol{\xi})$ is the Fourier-transform operator on spatial coordinates and $\boldsymbol{\xi}$ is the transform variable. In this formulation, the Dirac partial differential equation becomes an integral equation where the derivative is expressed through its spectral representation.  
Mathematically, the Dirac equation \eqref{eq:dirac_eq2} along with the Hamiltonian \eqref{eq:dirac_simplified} is a hermitian non-strictly hyperbolic system of equations \cite{leveque2002finite}. In principle, other numerical methods such as finite volumes or Galerkin could be used. An example of a finite volume discretization in 1D and its interpretation as a lattice Boltzmann method and quantum walk can be found in \cite{Succi2015}.
\section{Absorbing Layers}\label{sec:PML}
From a practical point of view, the time-dependent Dirac equation is considered on a bounded truncated physical domain denoted by $\mathcal{D}_{\textrm{Phy}}$. 
The pseudospectral method used to solve the Dirac equation naturally induces periodic boundary conditions on a bounded domain. What follows, is a general strategy which simultaneously i) avoids$/$reduces artificial wave reflections at the domain boundary, ii) limits the transfer of the wave from one side to the opposite one by periodicity.  To reach this goal, we add a layer $\mathcal{D}_{\textrm{PML}}$  surrounding $\mathcal{D}_{\textrm{Phy}}$, and stretch the coordinates in all the directions. The overall computational domain is next defined by: $\mathcal{D} = \overline{\mathcal{D}_{\textrm{Phy}}\cup \mathcal{D}_{\textrm{PML}}}$.   We refer to \cite{MOLPHYS} for the construction of PMLs for quantum wave equations and more specifically to \cite{pinaud} for the derivation and analysis of PMLs for the Dirac equation. 
Here, we outline the main features of this technique which is detailed in \cite{jcp2019}.
The first step of the implementation of PMLs is the following change of variables \cite{ZhengPML}  involving only the space variable:
\begin{eqnarray}
\widetilde{x}^{i} = x^{i} + e^{{\tt i}\theta}\int_{L^{* i}}^{x^{i}}\Sigma(s)ds,
\end{eqnarray}
where $\theta \in (0,\pi/2)$, $i=1,2,3$  and $\Sigma$ is a function to be determined. We then define
\begin{eqnarray*}
S^{i}(x^{i}) := 1+e^{{\tt i}\theta^{i}}\widetilde{\Sigma}(x^{i}),
\end{eqnarray*}
  the function $\widetilde{\Sigma}^{i}$ being given by
\begin{eqnarray*}
\widetilde{\Sigma}^{i}(x^{i}) = 
\left\{
\begin{array}{ll}
\Sigma(|x^{i}|-L^{i}), &  L^i_* \leq |x^{i}| < L^{i},\\
0, & |x^{i}| <L^i_*,
\end{array}
\right.
\end{eqnarray*}
where $L^i_*<L^{i}$, and such that $\{x^{i}\in\R \, : \, L^i_* \leq |x^{i}| < L^{i} \}$ is the absorbing layer. The partial derivatives are then transformed into
\begin{eqnarray}\label{pml2}
\partial_{i} \rightarrow \cfrac{1}{S^{i}(x^{i})}\partial_{i} =  \cfrac{1}{1+e^{{\tt i}\theta^{i}}\widetilde{\Sigma}^{i}(x^{i})}\partial_{i},
\end{eqnarray}
with $\widetilde{\Sigma}$ vanishing while $S^{i}$ is equal to $1$ in $\mathcal{D}_{\textrm{Phy}}$, respectively. On truncated domains, we will consider the transformation \eqref{pml2}, and the associated  new Hamiltonian
\begin{eqnarray}
H_{\textrm{PML}}  = - {\tt i} \boldsymbol{\alpha}({\boldsymbol x}) \cdot  {\boldsymbol T}({\boldsymbol x}) + \beta({\boldsymbol x}) m + F(t,{\boldsymbol x}),
\label{eq:hamiltonian_pml}
\end{eqnarray}
where ${\boldsymbol T} := ([S^{1}(x^{1})]^{-1} \partial_{1},[S^{2}(x^{2})]^{-1} \partial_{2},[S^{3}(x^{3})]^{-1} \partial_{3})^T$. Several types of functions can be selected. An exhaustive study of the absorbing functions $\Sigma$ is proposed in \cite{AntoineGeuzaineTang} for Schr\"odinger equations. Here are some examples:
\begin{eqnarray*}
\left.
\begin{array}{lll} 
\textrm{Type I: } \Sigma_0(x^{i}+\delta^{i})^2, &  \textrm{Type II: } \Sigma_0(x^{i}+\delta^{i})^3, & \textrm{Type III: }  -\Sigma_0/x^{i}, \\
\\
\textrm{Type IV: }  \Sigma_0/(x^{i})^2, & \textrm{Type V: } -\Sigma_0/x^{i} -\Sigma_0/\delta^{i}, &  \textrm{Type VI: } \Sigma_0/(x^{i})^2-\Sigma_0/(\delta^{i})^2,
\end{array}
\right.
\end{eqnarray*}
 where $\delta^{i}:=L^i-L_*^i$. From the pseudospectral point of view, the space-dependence of the coefficients $(S^{i})^{-1}$ again prevents the direct application of the Fourier transform on the equation, even in the flat-space case. In the latter case, the same pseudodifferential operator representation will still allow for combining the efficiency of the pseudospectral method and the computation of the non-constant coefficient Dirac equation (see \cite{jcp2019}). 
\section{Time-discretization}\label{sec:discrete}
This section is devoted to the time discretization of the Dirac equation, in flat and curved spaces. The main tool 
which is  often used is the operator splitting technique. The complications related to the appearance of spatial differential operators are relegated to the next section where the spatial discretization is discussed. 
\subsection{Time-discretization for Dirac equation in flat space}
\noindent In order to solve the Dirac equation in flat space, a natural approach consists in splitting the equation as follows. Here, an order-2 Strang splitting \cite{strang} is considered, but higher order operator splittings can naturally be used.  Let us consider a time interval from $t_n$ to $t_{n+1}$ and assuming $\psi(t_n,\cdot)$ is given, the exact formal solution to the Dirac equation in flat space is 
\begin{eqnarray}
\psi(t_{n+1},{\boldsymbol x})= \mathcal{T} \exp \left\{ -{\tt i} \int_{t_{n}}^{t_{n+1}} H_{\mathrm{flat}}(s,\boldsymbol{x}) ds \right\} \psi(t_{n},{\boldsymbol x}),
\end{eqnarray}
where $\mathcal{T}$ stands for the time-ordered exponential. The latter can be approximated to second order by the following symmetric decomposition \cite{suzuki1993general,Suzuki1990319}:
\begin{eqnarray}
\label{eq:op_flat}
\psi(t_{n+1},{\boldsymbol x})= 
e^{-{\tt i} \frac{\Delta t}{2} \left[\beta m + F_{\mathrm{flat}}(t_{n+1/2},\boldsymbol{x}) \right]} 
e^{-\Delta t \boldsymbol{\alpha} \cdot \nabla }    
e^{-{\tt i} \frac{\Delta t}{2} \left[\beta m + F_{\mathrm{flat}}(t_{n+1/2},\boldsymbol{x}) \right]}\psi(t_{n},{\boldsymbol x})
+ \mathcal{O}(\Delta t^{3}),
\end{eqnarray} 
where three exponential operators have been introduced. The first and third exponential operators can be evaluated analytically using the fact that
\begin{eqnarray}
\label{eq:exp_dirac}
e^{{\tt i}[\beta G(t,\boldsymbol{x}) + \boldsymbol{\alpha} \cdot \boldsymbol{G}(t,\boldsymbol{x})]} = 
\mathbb{I}_{4} \cos(|G|) 
+
{\tt i} \frac{[\beta G(t,\boldsymbol{x}) + \boldsymbol{\alpha} \cdot \boldsymbol{G}(t,\boldsymbol{x})]}{|G|} \sin(|G|),
\end{eqnarray}
for arbitrary functions $G$ and $\boldsymbol{G}$,  where we defined $|G| = \sqrt{G^{2}(t,\boldsymbol{x}) +  \boldsymbol{G}^{2}(t,\boldsymbol{x})}$. 
%
%
%
A common method to deal with the second differential operator is to use the Fourier transform $\mathcal{F}_{{\boldsymbol x}}$, as follows (still using an order 2-splitting). Step by step, denoting $t_{n^*}=t_n+\Delta t$, we have
\begin{eqnarray}\label{split}
\left\{
\begin{array}{lcll}
\psi(t_{n+1/2},{\boldsymbol x}) & = & e^{-{\tt i} \frac{\Delta t}{2} \left[\beta m + F_{\mathrm{flat}}(t_{n+1/2},\boldsymbol{x}) \right]}  \psi(t_{n},{\boldsymbol x}), & t\in [t_n,t_{n+1/2}],  \\
\psi(t_{n^*},{\boldsymbol x}) & = & \mathcal{F}_{{\boldsymbol x}}^{-1}\big\{e^{{\tt i}\boldsymbol{\alpha} \cdot \boldsymbol{\xi}\Delta t}  \psi(t_{n+1/2},\boldsymbol{\xi})\big\}, & t\in [t_n,t_{n^*}], \\
\psi(t_{n+1},{\boldsymbol x}) & = & e^{-{\tt i} \frac{\Delta t}{2} \left[\beta m + F_{\mathrm{flat}}(t_{n+1/2},\boldsymbol{x}) \right]}  \psi(t_{n^*},{\boldsymbol x}), & t\in [t_{n+1/2},t_{n+1}].
\end{array}
\right.
\end{eqnarray}
Then, the second step in Fourier space can also be evaluated \textit{via} \eqref{eq:exp_dirac}. After discretizing spatially, the second equation is commonly solved using the Fast Fourier Transform (FFT), resulting in an {\it operator splitting pseudospectral scheme} \cite{fft2,Bauke20112454,keitel,grobe,keitel5}.
\subsection{Time discretization for Dirac equation in curved space}
\noindent Following the same procedure as in flat space, described in the last section, an operator splitting approach can be introduced in curved space. For a time interval from $t_n$ to $t_{n+1}$ and assuming that $\psi(t_n,\cdot)$ is given, the exact formal solution to the Dirac equation in curved space is 
\begin{eqnarray}
\psi(t_{n+1},{\boldsymbol x})= \mathcal{T} \exp \left\{ -{\tt i} \int_{t_{n}}^{t_{n+1}} H(s,\boldsymbol{x}) ds \right\} \psi(t_{n},{\boldsymbol x}),
\end{eqnarray}
where $\mathcal{T}$ stands for the time-ordered exponential. Again, the latter can be approximated to third order accuracy by a symmetric decomposition \cite{suzuki1993general,Suzuki1990319}:
\begin{eqnarray}
\label{eq:op_split_curved}
\psi(t_{n+1},{\boldsymbol x})&=& 
e^{-{\tt i} \frac{\Delta t}{2} \left[\beta(\boldsymbol{x}) m + F(t_{n+1/2},\boldsymbol{x}) \right]} 
e^{-\Delta t \boldsymbol{\alpha} (\boldsymbol{x}) \cdot \nabla }    
e^{-{\tt i} \frac{\Delta t}{2} \left[\beta(\boldsymbol{x}) m + F(t_{n+1/2},\boldsymbol{x}) \right]}\psi(t_{n},{\boldsymbol x}) \nonumber \\
&& + \mathcal{O}(\Delta t^{3}).
\end{eqnarray} 
This has the same form as \eqref{eq:op_flat}, except for the second exponential, which now has a space-dependent Dirac matrix. The latter makes the direct use of the Fourier transform, as in the flat case, challenging because the efficient FFT cannot be used and the computational complexity would be $\mathcal{O}(N^{2})$, where $N$ is the number of lattice points. Our strategy is to approximate the exponential operator to simplify the problem and to be able to exploit the pseudodifferential form of the derivative operator. Two different types of approximation are introduced, leading to two different classes of numerical schemes:
\begin{enumerate}
	\item \textbf{Crank-Nicolson approximation:}
	This is obtained by formally approximating the exponential operator in its lowest order unitary form (1/1 Pad\'e's approximant of the exponential function):
	\begin{eqnarray}
	\label{eq:pade_approx}
	e^{-\Delta t \boldsymbol{\alpha} (\boldsymbol{x}) \cdot \nabla } = 
	\frac{\mathbb{I}_{4} - \frac{\Delta t}{2} \boldsymbol{\alpha} (\boldsymbol{x}) \cdot \nabla}
	{\mathbb{I}_{4} + \frac{\Delta t}{2} \boldsymbol{\alpha} (\boldsymbol{x}) \cdot \nabla}
	+ \mathcal{O}(\Delta t^{3}).
	\end{eqnarray}
	Then, the operator splitting can be implemented by the following sequence
	\begin{eqnarray}
	\label{eq:op_split_cn}
	\left\{
	\begin{array}{ll}
	\displaystyle
	\psi(t_{n+1/2},{\boldsymbol x})  =  e^{-{\tt i} \frac{\Delta t}{2} \left[\beta (\boldsymbol{x}) m + F(t_{n+1/2},\boldsymbol{x}) \right]}  \psi(t_{n},{\boldsymbol x}), & t\in [t_n,t_{n+1/2}],  \\
	\displaystyle
	\left[ \mathbb{I}_{4} + \frac{\Delta t}{2} \boldsymbol{\alpha} (\boldsymbol{x}) \cdot \nabla \right] \psi(t_{n^*},{\boldsymbol x})  =  \left[\mathbb{I}_{4} - \frac{\Delta t}{2} \boldsymbol{\alpha} (\boldsymbol{x}) \cdot \nabla\right]   \psi(t_{n+1/2},\boldsymbol{x}), & t\in [t_n,t_{n^*}], \\
	\displaystyle
	\psi(t_{n+1},{\boldsymbol x})  =  e^{-{\tt i} \frac{\Delta t}{2} \left[\beta (\boldsymbol{x}) m + F(t_{n+1/2},\boldsymbol{x}) \right]}  \psi(t_{n^*},{\boldsymbol x}), & t\in [t_{n+1/2},t_{n+1}].
	\end{array}
	\right. \nonumber \\
	\end{eqnarray}
	As shown below, this yields a semi-implicit numerical scheme.
	
	\item \textbf{Polynomial approximation:}
	This is obtained by formally approximating the exponential operator $e^{-\Delta t \boldsymbol{\alpha} (\boldsymbol{x}) \cdot \nabla }$, by a polynomial of the form 
	\begin{eqnarray}
	\sum_{q=0}^{N_{\mathrm{p}}} a_{q} P_{q}(-\Delta t \boldsymbol{\alpha} (\boldsymbol{x}) \cdot \nabla),
	\end{eqnarray}
	where $(P_{q})_{q=0,\cdots, N_{\mathrm{p}}}$ are some set of orthogonal polynomials (such as Taylor series, second order differencing or Chebychev polynomials, for example \cite{doi:10.1063/1.448136}) and $(a_{q})_{q=0,\cdots, N_{\mathrm{p}}}$ are the polynomial coefficients, fixed to have an accurate approximation of the exponential. Then, the operator splitting can be implemented by the following sequence
	\begin{eqnarray}
	\label{eq:op_split_poly}
	\left\{
	\begin{array}{lcll}
	\displaystyle
	\psi(t_{n+1/2},{\boldsymbol x})  &=&  e^{-{\tt i} \frac{\Delta t}{2} \left[\beta (\boldsymbol{x}) m + F(t_{n+1/2},\boldsymbol{x}) \right]}  \psi(t_{n},{\boldsymbol x}), & t\in [t_n,t_{n+1/2}],  \\
	\displaystyle
	\psi(t_{n^*},{\boldsymbol x})  &=& \displaystyle \sum_{q=0}^{N_{\mathrm{p}}} a_{q} P_{q}(-\Delta t \boldsymbol{\alpha} (\boldsymbol{x}) \cdot \nabla)  \psi(t_{n+1/2},\boldsymbol{x}), & t\in [t_n,t_{n^*}], \\
	\displaystyle
	\psi(t_{n+1},{\boldsymbol x})  &=&  e^{-{\tt i} \frac{\Delta t}{2} \left[\beta (\boldsymbol{x}) m + F(t_{n+1/2},\boldsymbol{x}) \right]}  \psi(t_{n^*},{\boldsymbol x}), & t\in [t_{n+1/2},t_{n+1}].
	\end{array}
	\right. \nonumber \\
	\end{eqnarray}
	As shown below, this yields explicit numerical schemes which do not require the solution to a linear system.
\end{enumerate}    
These two strategies are well-known for the Schr\"odinger equation, but are adapted here to the Dirac equation in curved space and combined to the pseudodifferential representation of the derivative. In the next sections, some methods will be given based on the discretization of Eqs. \eqref{eq:op_split_cn} and \eqref{eq:op_split_poly} with spectral accuracy in space.

The Strang splitting was introduced mostly to be consistent with the traditional splitting strategy used in flat space and to isolate the part of the equation that requires a special treatment with pseudodifferential operators. In principle, this is not mandatory and unsplit schemes could also be used in combination with an approximation of the time-ordered exponential. Nevertheless, notice that the overall second-order of accuracy in time is still preserved by the operator splitting.


\section{Space-discretization for the Dirac equation in curved space}\label{sec:space}

In this section, the spatial discretization of the Dirac equation in curved space is described, based on the two different approaches described in the last section. Throughout, we assume that the Dirac equation is solved in curved space on a truncated domain $[-a_1,a_1]\times [-a_2,a_2]\times [-a_3,a_3] \varsubsetneq \R^3$. The cases of flat space or the one with time-dependent coefficients are standard and not recalled here.  We define two sets of grid-points in real and Fourier spaces by
\begin{eqnarray*}
	\mathcal{D}^{(x)}_{N}=\big\{{\boldsymbol x}_{k} &:=& {\boldsymbol x}_{k_1,k_2,k_3}=(x^{1}_{k_1},x^{2}_{k_2},x^{3}_{k_3})\big\}_{k \in \mathcal{O}^{(x)}_{N}}, \\
		\mathcal{D}^{(\xi)}_{N}=\big\{{\boldsymbol \xi}_{p} &:=& {\boldsymbol \xi}_{p_1,p_2,p_3}=(\xi^{1}_{p_1},\xi^{2}_{p_2},\xi^{3}_{p_3})\big\}_{p \in \mathcal{O}^{(\xi)}_{N}},
\end{eqnarray*}
where $N$, $k$ and $p$ are multi-indices ($N:= (N_1,N_2,N_3)$, $k=(k_{1},k_{2},k_{3})$ and $p=(p_{1},p_{2},p_{3})$, respectively and $N_i \in 2\N^*$), and with
\begin{eqnarray*}
	\mathcal{O}^{(x)}_{N} &=&\bigg\{ k \in \N^3/ \left(k_{i}=0,\cdots,N_i-1 \right)_{i=1,2,3} \bigg\}, \\
	%
	\mathcal{O}^{(\xi)}_{N} &=& \bigg\{ p \in \N^3/ \left(p_{i}=-\frac{N_{i}}{2},\cdots,\frac{N_{i}}{2}-1 \right)_{i=1,2,3}\bigg\}. 
\end{eqnarray*}
The set $\mathcal{D}^{(x)}_{N}$ defines a mesh with equidistant point positions in each dimension with sizes (for $i=1,2,3$)
\begin{equation*}
\displaystyle  x^{i}_{k_i+1}-x^{i}_{k_i}=h_i=2a_i/N_i \, .
\end{equation*}
One can deduce that the discrete wavenumbers in Fourier space are given by (for $i=1,2,3$)
\begin{eqnarray*}
\xi^{i}_{p_{i}}&=&p_{i}\pi/a_i \, . 
\end{eqnarray*}
The wave function $\psi(t,\boldsymbol{x})$ is discretized spatially by a projection onto the spatial mesh while $\widetilde{\psi}(t,\boldsymbol{\xi})$ is discretized on the momentum mesh. Thus, we denote by $\psi^{n}_{k}$ the approximate wavefunction at time $t_{n}$ and position ${\boldsymbol x}_k$ and by $\widetilde{\psi}^{n}_{p}$ the wave function in momentum space at time $t_{n}$ and momentum $\boldsymbol{\xi}_{p}$. The discrete wave functions $\psi^{n}_{k}$ and $\widetilde{\psi}^{n}_{p}$ are related by the discrete Fourier transform pair:
\begin{eqnarray*}
	\widetilde{\psi}^{n}_{p} =  \displaystyle \sum_{k=0}^{N-1}\psi^{n}_{k}e^{-{\tt i}\boldsymbol{\xi}_{p} \cdot (\boldsymbol{x}_{k}+\boldsymbol{a})}
	\;\;,\;\;
	\psi^{n}_{k} =  \displaystyle \frac{1}{N} \sum_{p=-N/2}^{N/2-1}\widetilde{\psi}^{n}_{p}e^{{\tt i}\boldsymbol{\xi}_{p} \cdot (\boldsymbol{x}_{k}+\boldsymbol{a})} \, ,
\end{eqnarray*}
%
where $\boldsymbol{a}=(a_1,a_2,a_3)$. Armed with this notation,  we can write the partial discrete Fourier coefficients in each dimension as, for $i=1,2,3$:
\begin{eqnarray*}
		 \widetilde{\psi}^{n}_{k|k_{i} \rightarrow p_{i}} &=&  
		 \boldsymbol{\mathcal{F}}_{i}(\psi^{n}_{k}) 
	 = \sum_{k_i=0}^{N_i-1}\psi^{n}_{k}e^{-{\tt i}\xi^{i}_{p_{i}}(x^{i}_{k_i}+a_i)} \\
		%
		\psi^{n}_{k}  &=& 
		\boldsymbol{\mathcal{F}}^{-1}_{i} (\widetilde{\psi}^{n}_{k})
		= 
		\cfrac{1}{N_i}\sum_{p_{i}=-N_i/2}^{N_i/2-1}\widetilde{\psi}^{n}_{k|k_{i} \rightarrow p_{i}}e^{{\tt i}\xi^{i}_{p_{i}}(x^{i}_{k_i}+a_i)},
\end{eqnarray*}
where the notation $k|k_{i} \rightarrow p_{i}$ means that the index $k_{i}$ in the set $k$ is replaced by the index $p_{i}$ and where the partial discrete Fourier transform operator in the $i$ coordinate is  denoted by $\boldsymbol{\mathcal{F}}_{i}(\cdot)$.

In order to approximate the partial derivatives, we use pseudospectral approximations of the pseudodifferential representation of the derivative operators. This leads to the following approximate first-order partial derivatives:
\begin{eqnarray}\label{trick}
\partial_{i}\psi(t_n,{\boldsymbol x}_{k})  \approx  
\big\{[[\partial_i]]\psi^n\big\}_{k}  :=  
\cfrac{1}{N_i}\sum_{p_{i}=-N_i/2}^{N_i/2-1}{\tt i}\xi^{i}_{p_{i}}\widetilde{\psi}^{n}_{k|k_{i} \rightarrow p_{i}}e^{{\tt i}\xi^{i}_{p_{i}}(x^{i}_{k_i}+a_i)}.
%
%
\end{eqnarray}
This is the spectral representation of the derivative which, under standard assumptions on the smoothness of the wave function,  has a spectral accuracy \cite{q1,GOTTLIEB200183}. Another representation of the spectral derivative can be found in \ref{app:spec} in terms of the differentiation matrix. In the next sections, we will exploit these relations to obtain accurate numerical schemes. This approach not only allows to select the spatial steps as large as wanted, but it  also  preserves the very high spatial accuracy, the parallel computing structure and the scalability of the split method developed in \cite{cpc2012}. In practice, we use the FFT to implement the Discrete Fourier Transform (DFT). This strategy is now combined with the two time discretization proposed in the last section to solve the Dirac equation in curved space.

\subsection{Numerical scheme I: Crank-Nicolson scheme}\label{subsec:CN}

A Crank-Nicolson scheme is obtained by discretizing \eqref{eq:op_split_cn} and by using the spectral operator for computing the derivative. The main idea then consists in approximating $\mathcal{F}^{-1}_{\boldsymbol{x}} \big\{ {\tt i}\boldsymbol{\xi} \mathcal{F}_{\boldsymbol{x}}\{ \psi \} (t,\boldsymbol{\xi}) \big\}$ by using the discrete operator $[[\nabla]]$ defined in \eqref{trick}. This yields
\begin{eqnarray}
\label{eq:cn}
\left.
\begin{array}{rcll}
\psi^{n+1/2}_{k} & = & e^{-{\tt i} \frac{\Delta t}{2} \left[\beta_{k} m + F^{n+1/2}_{k} \right]}  \psi^{n}_{k}, & t\in [t_n,t_{n+1/2}],  \\
G_{k} \psi^{n^*}_{k} & = & \widetilde{G}_{k}  \psi^{n+1/2}_{k}, & t\in [t_n,t_{n^*}], \\
\psi^{n+1}_{k} & = & e^{-{\tt i} \frac{\Delta t}{2} \left[\beta_{k} m + F^{n+1/2}_{k} \right]}  \psi^{n^*}_{k}, & t\in [t_{n+1/2},t_{n+1}],
\end{array}
\right.
\end{eqnarray} 
where we defined $\boldsymbol{\alpha}_{k}=\boldsymbol{\alpha}({\boldsymbol x}_{k})$, $ \beta_{k}=\beta({\boldsymbol x}_{k})$ and $F^{n}_{k} = F(t_{n}, {\boldsymbol x}_{k})$. We also introduced the operators
\begin{eqnarray}
G_{k} &:=& \mathbb{I}_{4} + \frac{\Delta t}{2} \boldsymbol{\alpha}_{k} \cdot [[\nabla]], \\
\widetilde{G}_{k} &:=& \mathbb{I}_{4} - \frac{\Delta t}{2} \boldsymbol{\alpha}_{k} \cdot [[\nabla]] ,
\end{eqnarray}
for convenience. An explicit example of this procedure for a given 2D metric can be found in \ref{app:cn_2d}.

Then, the numerical solution can be obtained by implementing the following algorithm.
From time $t_n$ to $t_{n+1}$, and adding the PML ${\bf T}$ (see \eqref{eq:hamiltonian_pml}), the 3-steps scheme \eqref{eq:cn} explicitly reads:
\begin{itemize}
	\item {\it Step 1.} From $t_n$ to $t_{n+1/2}$ such that $t_{n+1/2} = t_n + \Delta t/2$,  with initial data  $\psi^{n}_{k}$:
	\begin{eqnarray}\label{step1}
	\psi^{n+1/2}_{k} & = & e^{-{\tt i} \frac{\Delta t}{2} \left[\beta_{k} m + F^{n+1/2}_{k} \right]}  \psi^{n}_{k}.
	\end{eqnarray}
	To perform this step, an exact or approximate expression of the exponential operator is required. This can be achieved in different ways, depending on the metric chosen and the form of $\beta$ and $F$. The first thing to note here is that both $\beta$ and $F$ are 4-by-4 matrices. One possibility is then to use one of the numerical techniques described in \cite{doi:10.1137/S00361445024180} to compute the exponential of the matrix.

	\item {\it Step 2.} The second step, given by
	\begin{eqnarray}
	G_{k} \psi^{n^*}_{k} & = & \widetilde{G}_{k}  \psi^{n+1/2}_{k},
	\end{eqnarray}
	can be written as a linear system of equations, by using the discrete pseudospectral representation of the derivative. Naively, one would construct the matrix $G_{k}$ explicitly as done in \ref{APXA}, and solve the corresponding linear system. However, this procedure has the same computational complexity as the evaluation of convolution products. The computational efficiency of this step can be improved significantly by using a Krylov iteration solver (GMRES, conjugate gradient). This follows a technique developed before and we refer the reader to \cite{AntoineGeuzaineTang} for more details.

	\item {\it Step 3.} The third step is given by
	\begin{eqnarray}\label{step3}
	\psi^{n+1}_{k} & = & e^{-{\tt i} \frac{\Delta t}{2} \left[\beta_{k} m + F^{n+1/2}_{k} \right]}  \psi^{n^*}_{k}.
	\end{eqnarray}
	The same technique as in Step 1 can be used  to evaluate the matrix exponential. 
	
\end{itemize}

In flat space, it is well-known that the $\ell^2-$norm of the 4-spinor must be preserved, while in static curved space, the $\ell_{\gamma}^2-$norm is preserved.


\begin{prop}
Assume that $(\alpha^{i,n}_{k})_{i=1,2,3}$ are hermitian, and that $F$ is a bounded function. The numerical scheme \eqref{step1}-\eqref{step3} is unconditionally $\ell^2$-stable, and preserves the $\ell^2-$norm in flat space.
\end{prop}
{\bf Proof.} The proof is straighforward, as it mainly relies on i) the hermitivity of the Dirac matrices  $(\alpha^{i}_{k})^{\dagger} = \alpha^{i}_{k}$, with $i=1,2,3$, for steps \eqref{step1}, \eqref{step3}, ii) the definition of $G $ and $\widetilde{G}$. Using the matrix representation of the derivative given in \ref{app:spec}, we start from the hermitian transpose of the differentiation matrices:
\begin{eqnarray}
\overline{A}^{i}_{k_{i} k'_{i}}  &=& 
-\cfrac{1}{N_i}
\sum_{p_{i}=-N_i/2}^{N_i/2-1}
{\tt i}\xi^{i}_{p_{i}}e^{{-\tt i}\xi^{i}_{p_{i}}(x^{i}_{k_i} - x^{i}_{k'_i})} \, .
\end{eqnarray}
Setting $p'_{i}=-p_{i}$ and using the fact that $\xi^{i}_{-p'_{i}} = -\xi^{i}_{p'_{i}}$, we obtain
\begin{eqnarray}
\overline{A}^{i}_{k_{i} k'_{i}}  &=& 
-\cfrac{1}{N_i}
\sum_{-p_{i}=-N_i/2}^{N_i/2-1}
{\tt i}\xi^{i}_{-p_{i}}e^{{-\tt i}\xi^{i}_{-p_{i}}(x^{i}_{k_i} - x^{i}_{k'_i})} \, , \\
 &=& 
\cfrac{1}{N_i}
\sum_{-p_{i}=-N_i/2}^{N_i/2-1}
{\tt i}\xi^{i}_{p_{i}}e^{{\tt i}\xi^{i}_{p_{i}}(x^{i}_{k_i} - x^{i}_{k'_i})} \, .
\end{eqnarray}
%
This implies that, by construction: $\overline{[[\nabla ]]} = -[[\nabla]]$, and as a consequence  
\begin{eqnarray*}
\left.
\begin{array}{lcl}
\widetilde{G}^{\dagger} & = & \mathbb{I} - (\Delta t/2) \big({\boldsymbol \alpha}\cdot [[ \nabla ]]\big)^{\dagger} \\
& = & \mathbb{I} + (\Delta t/2) {\boldsymbol \alpha}\cdot [[ \nabla ]] .
\end{array}
\right.
\end{eqnarray*}
As the coefficients of the Dirac equation are bounded in space and time, we  trivially conclude about the $\ell^2$-stability. In flat space, ${\boldsymbol \alpha}$ is constant and $F$ is purely real the $\ell^2-$norm is trivially preserved (Steps 1 to 3 are unitary). $\Box$


\subsection{Numerical II: polynomial scheme}

Polynomial schemes are obtained using the same procedure as for the Crank-Nicolson method, i.e. by discretizing \eqref{eq:op_split_poly} and by using the spectral operator for computing the derivative. Again, we approximate $\mathcal{F}^{-1}_{\boldsymbol{x}} \big\{ {\tt i} \boldsymbol{\xi} \mathcal{F}_{\boldsymbol{x}}\{ \psi \} (t,\boldsymbol{\xi}) \big\}$ by using the discrete operator $[[\nabla]]$ defined in \eqref{trick}. This yields
\begin{eqnarray}
\label{eq:poly}
\left.
\begin{array}{rcll}
\psi^{n+1/2}_{k} & = & e^{-{\tt i} \frac{\Delta t}{2} \left[\beta_{k} m + F^{n+1/2}_{k} \right]}  \psi^{n}_{k}, & t\in [t_n,t_{n+1/2}],  \\
\psi^{n^*}_{k} & = &  \sum_{q=0}^{N_{\mathrm{p}}} a_{q} P_{q}(-\Delta t \boldsymbol{\alpha}_{k} \cdot [[\nabla]])   \psi^{n+1/2}_{k}, & t\in [t_n,t_{n^*}], \\
\psi^{n+1}_{k} & = & e^{-{\tt i} \frac{\Delta t}{2} \left[\beta_{k} m + F^{n+1/2}_{k} \right]}  \psi^{n^*}_{k}, & t\in [t_{n+1/2},t_{n+1}].
\end{array}
\right.
\end{eqnarray} 
The first and second steps are exactly the same as in the Crank-Nicolson scheme (see \eqref{eq:cn}) and thus are not discussed here. The second step, on the other hand, is different because it does not require a solution of a linear system. The main challenge is in computing powers of the operator $[[\nabla]]$. This is performed by using FFTs, where the number of FFT pairs  is given by the order of the polynomial.

Every scheme of this form is explicit and thus, is \textit{a priori} at best conditionally stable. As a matter of fact,  numerical experiments often show an instability of the numerical solution. However, the Crank-Nicolson scheme derived above naturally requires the solution to a large linear system at each time iteration. In order to improve the efficiency while keeping a reasonable accuracy, we propose a 2-steps {\it polynomial scheme} with directional splitting. In the following, in order to simplify the presentation, we will assume that $\boldsymbol{\alpha}(\boldsymbol{x})$ is of the form
\begin{eqnarray}
\label{eq:alpha_form}
{\boldsymbol \alpha}({\boldsymbol x}) & = & {\bf a}({\boldsymbol x})\cdot {\boldsymbol \alpha}, 
\end{eqnarray}
where ${\bf a}({\boldsymbol x}):=(a^{1}({\boldsymbol x}),a^{2}({\boldsymbol x}),a^{3}({\boldsymbol x}))^T$ and where $(a^{i})_{i=1,\cdots,3}$ are space dependent scalar functions. Using directional splitting and the form \eqref{eq:alpha_form} for the matrices, the scheme \eqref{eq:poly} is slightly modified to
\begin{eqnarray}
\label{eq:poly_mod}
\left.
\begin{array}{rcll}
\psi^{n+1/2}_{k} & = & e^{-{\tt i} \frac{\Delta t}{2} \left[\beta_{k} m + F^{n+1/2}_{k} \right]}  \psi^{n}_{k}, & t\in [t_n,t_{n+1/2}],  \\
\psi^{n^{*}_{1}}_{k} & = &  \sum_{q=0}^{N_{\mathrm{p}}} a_{q} P_{q}(-\Delta t a^{1}_{k}\alpha^{1}  [[\partial_{1}]])   \psi^{n+1/2}_{k}, & t\in [t_n,t_{n^*}], \\
\psi^{n^{*}_{2}}_{k} & = &  \sum_{q=0}^{N_{\mathrm{p}}} a_{q} P_{q}(-\Delta t a^{2}_{k}\alpha^{2} [[\partial_{2}]])   \psi^{n^{*}_{1}}_{k}, & t\in [t_n,t_{n^*}], \\
\psi^{n^{*}_{3}}_{k} & = &  \sum_{q=0}^{N_{\mathrm{p}}} a_{q} P_{q}(-\Delta t a^{3}_{k}\alpha^{3} [[\partial_{3}]])   \psi^{n^{*}_{2}}_{k}, & t\in [t_n,t_{n^*}], \\
\psi^{n+1}_{k} & = & e^{-{\tt i} \frac{\Delta t}{2} \left[\beta_{k} m + F^{n+1/2}_{k} \right]}  \psi^{n^{*}_{3}}_{k}, & t\in [t_{n+1/2},t_{n+1}].
\end{array}
\right.
\end{eqnarray} 
Let us remark  that PMLs can easily be included by simply replacing $a^{i}({\boldsymbol x})$ by $a^{i}({\boldsymbol x})/S^{i}(x^{i})$, for $i=1,2,3$ in the scheme below. Steps 1. and 5. are identical to the {\bf Numerical Scheme I} \eqref{subsec:CN}. The principle of the second step of {\bf Numerical Scheme II} consists in approximating the evolution of each direction by a Taylor expansion and by diagonalizing the Dirac matrix.  We denote by $\Lambda=\textrm{diag}(1,1,-1,-1)$ and $\Pi^{i}$ the transition matrices, such that $\alpha^{i} = \Pi^{i}\Lambda \Pi^{i,\dagger}$, for $i=1,2,3$. In addition, in the $i$-direction, we set $\phi^{n}_{k}=\Pi^{i,\dagger}\psi^{n}_{k}$.
Then, the time evolution is  approximated as follows (with $n^{*}_{0} = n+1/2$):
\begin{eqnarray}\label{rem2}
\left.
\begin{array}{lcl}
\phi^{n^{*}_{i}}_{k} & = & \phi^{n^{*}_{i-1}}_{k} -\Delta ta^{i}_{k} \Lambda ([[\partial_{i}]]\phi^{n^{*}_{i-1}})_{k} + \mathcal{O}(\Delta t),\\
& = & a^{i}_{k} \boldsymbol{\mathcal{F}}_{i}^{-1}\big[\big(1-{\tt i} \Delta t \Lambda \xi^{i}\big) \boldsymbol{\mathcal{F}}_i(\phi^{n^{*}_{i-1}}_{k})\big] + \big(1- a^{i}_{k} \big)\phi^{n^{*}_{i-1}}_{k} + \mathcal{O}(\Delta t) \, .
\end{array}
\right.
\end{eqnarray}
The first line of this equation is a polynomial Taylor scheme and is usually unstable while the second line is just a re-writing of the first line. Stability is recovered when the first term on the right-hand-side is approximated by an exponential, as $1-{\tt i} \Delta t \Lambda \xi^{i} = e^{-{\tt i} \Delta t \Lambda \xi^{i}} + \mathcal{O}(\Delta t^{2})$, which is accurate to second order in time. Then, the scheme reads
\begin{eqnarray}\label{eq:second_scheme}
\left.
\begin{array}{lcl}
\phi^{n^{*}_{i}}_{k}
& = & a^{i}_{k} \boldsymbol{\mathcal{F}}_{i}^{-1}\big[e^{-{\tt i} \Delta t \Lambda \xi^{i}} \boldsymbol{\mathcal{F}}_i(\phi^{n^{*}_{i-1}}_{k})\big] + \big(1- a^{i}_{k} \big)\phi^{n^{*}_{i-1}}_{k} \, .
\end{array}
\right.
\end{eqnarray}
Finally, to recover the wave function, we set $\psi^{n^{*}_{i}}_{k}=\Pi^{i}\phi^{n^{*}_{i}}_{k}$. We proceed similarly in the other directions. As it is not totally obvious, we next prove the consistency of {\bf Numerical Scheme II}. 
\begin{prop}
The {\bf Numerical Scheme II} is consistent with \eqref{eq:dirac_eq2} .
\end{prop}
{\bf Proof.} The analysis of the consistency only requires a focus on one of the steps 2-4, as the other steps are similar or standard. From
\begin{eqnarray*}
e^{-{\tt i} \Delta t \xi^{i}_p \Lambda} & = & \mathbb{I}_{4} - {\tt i} \Delta t \xi^{i}_p \Lambda  - \Delta t^2 (\xi^{i}_p)^2 \Lambda^2  + \mathcal{O}(\Delta t^3) \, ,
\end{eqnarray*}
we can write the first term of \eqref{eq:second_scheme} as
\begin{eqnarray*}\label{eq:second_scheme_bis}
\Xi^{n^{*}_{i}}_{k}
& := & \boldsymbol{\mathcal{F}}_{i}^{-1}\big[\Pi^{i}\big(\mathbb{I}_4   - {\tt i} \Delta t \xi^{i} \Lambda  - \Delta t^2 (\xi^{i})^2 \Lambda^2 \big)\Pi^{i,\dagger} \boldsymbol{\mathcal{F}}_i(\psi^{n^{*}_{i-1}}_{k})\big]  + \mathcal{O}(\Delta t^{3}) , \\
&=& \psi^{n^{*}_{i-1}}_{k} - \Delta t \alpha^{i} \partial_{i} \psi^{n^{*}_{i-1}}_{k}
+ \Delta t^2 \mathbb{I}_{4} \partial_{i}^{2} \psi^{n^{*}_{i-1}}_{k} + \mathcal{O}(\Delta t^{3}) ,
\end{eqnarray*}
since $\Pi^{i}\Lambda\Pi^{i,\dagger} = \alpha^{i}$. Thus, \eqref{eq:second_scheme} is written as
\begin{eqnarray}
\label{eq:scheme_with_xi}
\psi^{n^{*}_{i}}_{k}
& = & a^{i}(\boldsymbol{x}_{k}) \Xi^{n^{*}_{i-1}}_{k} + \big(1- a^{i}(\boldsymbol{x}_{k}) \big)\psi^{n^{*}_{i-1}}_{k} \\
& = & \psi^{n^{*}_{i-1}}_{k} -  a^{i}(\boldsymbol{x}_{k}) \Delta t \alpha^{i} \partial_{i} \psi^{n^{*}_{i-1}}_{k}
+  a^{i}(\boldsymbol{x}_{k}) \Delta t^2 \mathbb{I}_{4} \partial_{i}^{2} \psi^{n^{*}_{i-1}}_{k} + \mathcal{O}(\Delta t^{3}),
\end{eqnarray}
which is equivalent to
\begin{eqnarray}\label{LW}
	{\tt i} \cfrac{\psi^{n^{*}_{i}}_{k} - \psi^{n^{*}_{i-1}}_{k}}{\Delta t}
	& = &  - {\tt i} a^{i}(\boldsymbol{x}_{k})  \alpha^{i} \partial_{i} \psi^{n^{*}_{i-1}}_{k}
	+ {\tt i} a^{i}(\boldsymbol{x}_{k}) \Delta t \mathbb{I}_{4} \partial_{i}^{2} \psi^{n^{*}_{i-1}}_{k} + \mathcal{O}(\Delta t^{2}).
\end{eqnarray}
We have proven that is consistent at order $1$ in time, with
\begin{eqnarray*}
{\tt i}\partial_{t} \psi(t,{\boldsymbol x})  &=&-{\tt i} a^{i}({\boldsymbol x}) \alpha^{i}  \partial_{x}  \psi(t,{\boldsymbol x}),
\end{eqnarray*}
and this transport-like equation is equivalent to the 2 to 4 steps in the operator splitting. $\Box$
\\
\\
We now provide a stability result.
\begin{prop}
Let us assume that $ 0 \leq \sup_{{\boldsymbol  x}}a^{i}(\boldsymbol{x}_{k}) \leq C$ with $i=1,2,3$.
Then, the  {\bf Numerical Scheme II} approximating \eqref{eq:dirac_eq2}, which is assumed to be well-posed, is unconditionally $\ell^2-$stable.
\end{prop}
{\bf Proof.} Again the proof relies on Step 2. Let us focus on one of the direction and assume first that $C\leq 1$. As $\Pi^{i}$ is unitary, we trivially get : $\big|\Xi_{k}^{n} \big|^2_{2}  =  \big|\psi_{k}^{n}\big|^2_{2}$, where $|\cdot|_2$ denotes the  $\ell^2-$norm in $\C^4$ on spinor components. Consequently, this yields
\begin{eqnarray*}
\left.
\begin{array}{lclcl}
\displaystyle \| {\boldsymbol \Xi}_{k}^{n} \|_{2} & := & \displaystyle \Big(h^3\sum_{k=0}^{N} \big|\Xi_{k}^{n}\big|_2^2 \Big)^{1/2} & = & \|{\boldsymbol \psi}_{k}^{n}\|_{2} \, .
\end{array}
\right.
\end{eqnarray*}
Moreover, denoting $\mathcal{R}\{z\}$ (resp. $\mathcal{I}\{z\}$) the real (resp. imaginary) part of $z$ and using \eqref{eq:scheme_with_xi}, we get
\begin{eqnarray}
\big|\psi^{n^{*}_{i}}_{k}\big|_{2}^{2}
& = & [a^{i}(\boldsymbol{x}_{k})]^{2} \big|\Xi^{n^{*}_{i-1}}_{k}\big|_{2}^{2} + \big(1- a^{i}(\boldsymbol{x}_{k}) \big)^{2}\big|\psi^{n^{*}_{i-1}}_{k}\big|_{2}^{2} \nonumber \\
&&+  2a^{i}(\boldsymbol{x}_{k})\big(1- a^{i}(\boldsymbol{x}_{k}) \big)
\left[ \mathcal{R}\{\Xi^{n^{*}_{i-1}}_{k}\} \mathcal{R}\{\psi^{n^{*}_{i-1}}_{k}\} + \mathcal{I}\{\Xi^{n^{*}_{i-1}}_{k}\} \mathcal{I}\{\psi^{n^{*}_{i-1}}_{k}\}    \right], \nonumber \\
\label{interm}
& \leq & [a^{i}(\boldsymbol{x}_{k})]^{2} \big|\Xi^{n^{*}_{i-1}}_{k}\big|_{2}^{2} + \big(1- a^{i}(\boldsymbol{x}_{k}) \big)^{2}\big|\psi^{n^{*}_{i-1}}_{k}\big|_{2}^{2} \nonumber \\
&&+ 2\big|a^{i}(\boldsymbol{x}_{k})\big|\big|1- a^{i}(\boldsymbol{x}_{k})| \big|\Xi^{n^{*}_{i-1}}_{k}\big|_{2}\big|\psi^{n^{*}_{i-1}}_{k}\big|_{2},
\end{eqnarray}
where the Cauchy-Schwarz inequality was used to obtain \eqref{interm}.
Expanding \eqref{interm} and using $\big|\Xi^{n^{*}_{i-1}}_{k} \big|^2_{2}  =  \big|\psi^{n^{*}_{i-1}}_{k} \big|^2_{2}$, we easily  obtain the inequality:  $\big|\psi^{n^{*}_{i}}_{k} \big|^2_2 \leq  \big|\psi^{n^{*}_{i-1}}_{k} \big|^2_2$. We deduce that $\| \boldsymbol{\psi}^{n^{*}_{i}}_{k} \|_{2} \leq  \|\boldsymbol{\psi}^{n^{*}_{i-1}}_{k}\|_{2}$, for $i=1,2,3$. 
Thus, we obtain: $\| {\boldsymbol \psi}_{k}^{n^*} \|_{2} \leq  \|{\boldsymbol \psi}_{k}^{n}\|_{2}$. The  stability analysis  from Steps 1 and 5 is straightforward, hence leading to
\begin{eqnarray*}
\left.
\begin{array}{lcl}
\| {\boldsymbol \psi}_{k}^{n+1} \|_{2} & \leq &  \|{\boldsymbol \psi}_{k}^{n}\|_{2} \, .
\end{array}
\right.
\end{eqnarray*}
This concludes the proof for $C\leq 1$. For $C>1$, we proceed similarly, and conclude with the Gronwall's inequality. $\Box$
\\
\\
We can extend the above result to a second-order time scheme, by replacing the second step \eqref{eq:scheme_with_xi} by
\begin{eqnarray}
\label{step2bis}
\psi^{n^{*}_{i}}_{k}
& = & a^{i}(\boldsymbol{x}_{k}) \Xi^{n^{*}_{i-1}}_{k} + \big(1- a^{i}(\boldsymbol{x}_{k}) \big)\psi^{n^{*}_{i-1}}_{k} +  a^{i}(\boldsymbol{x}_{k}) \Delta t^2 \mathbb{I}_{4} [[\partial_{i}^{2}]]\Xi^{n^{*}_{i-1}}_{k} \, .
\end{eqnarray}
The second order in time is obtained thanks to the addition of  the rightmost term to the scheme \eqref{eq:scheme_with_xi}. The stability occurs from the implicitation of the correction$/$anti-diffusion term  $a^{i}(\boldsymbol{x}_{k}) \Delta t^2 \mathbb{I}_{4} [[\partial_{i}^{2}]]\psi$.

\subsection{Computational complexity analysis}

This paragraph is dedicated to the analysis of the computational complexity of the presented methods. In particular, we compare the complexity with i) the direct implicit method based on the direct application of the FFT on the equation, involving spatial convolution products, ii) the Crank-Nicolson scheme  and the polynomial scheme. Each {\it time iteration} requires on a $N$-point grid the following operations:
\begin{itemize}
\item  The implicit direct method requires $\mathcal{O}(N^2 + N^{\nu})$ operations, with $\nu>1$. The first term is due to the approximation of the convolution products by standard quadrature rules, and the second term comes from the numerical computation of the solution to the linear system. 
\item The Crank-Nicolson scheme  needs $\mathcal{O}(N\log N + N^{\nu})$ operations, for $\nu>1$. The first term is related
 to the approximation of the convolution products by FFT-method, while the second  one comes from the solution to  the linear system. 
\item The polynomial scheme implies  $\mathcal{O}(N\log N)$ operations since  the convolution products are computed by an
FFT-method, the rest of the scheme being linear.
\end{itemize}



%
%
%
%
%

\section{Numerical experiment for the curved-space Dirac equation}\label{sec:numerics}
Several specific examples will be considered, which basically correspond to different metrics. The objective is to demonstrate how simple and efficient is the proposed methodology. We will start with simple one-dimensional tests, then will consider more elaborated two-dimensional physical configurations. In 1-D and 2-D, the Dirac equation in curved space is slightly modified compared to 3-D. In particular, the dimension of Dirac matrices is 2-by-2, instead of 4-by-4. The numerical schemes described previously can be straightforwardly adapted to these cases.

\subsection{Static spacetime} 
\begin{figure}
\begin{center}
\includegraphics[height=4cm,keepaspectratio]{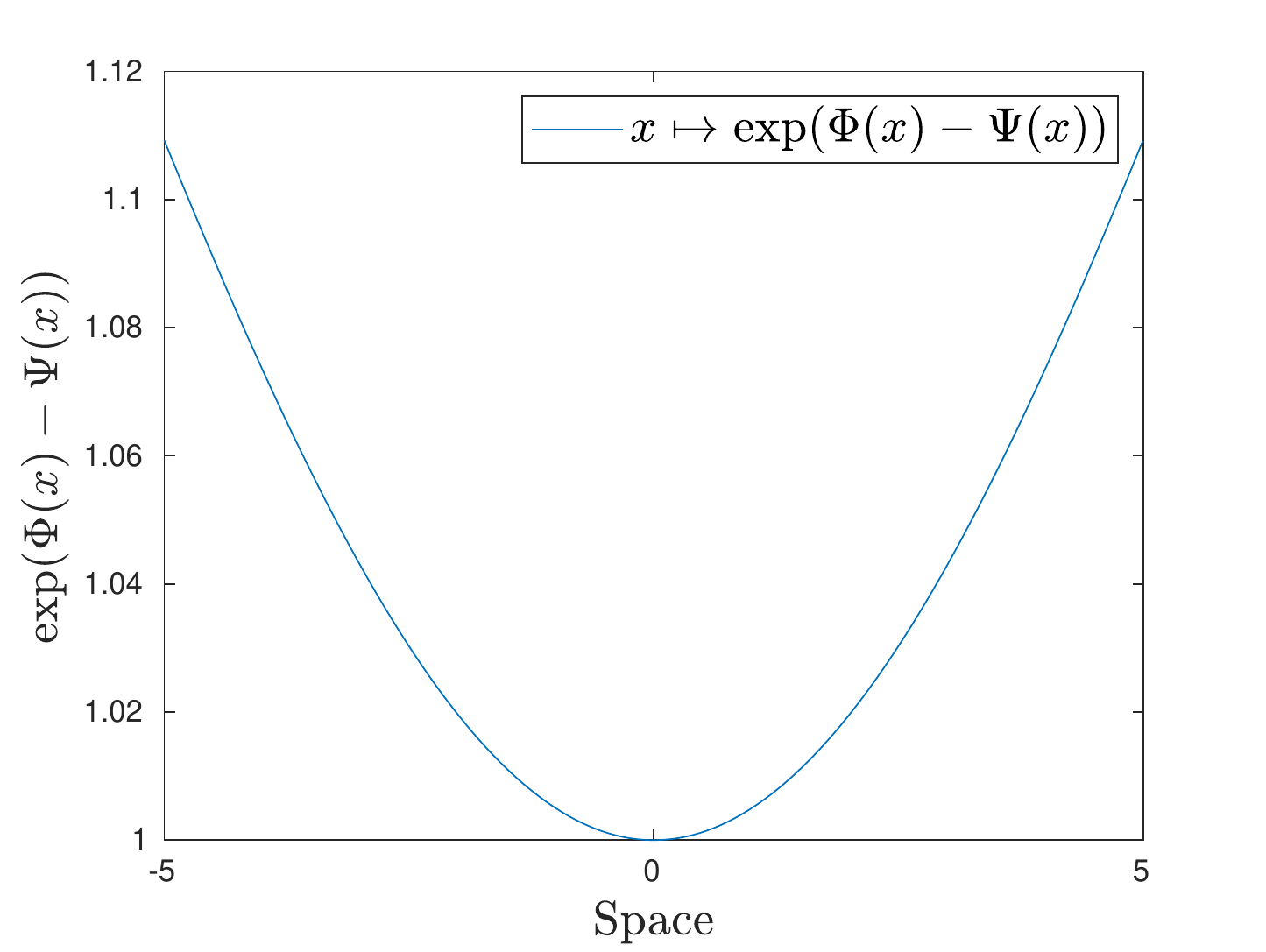}
\includegraphics[height=4cm,keepaspectratio]{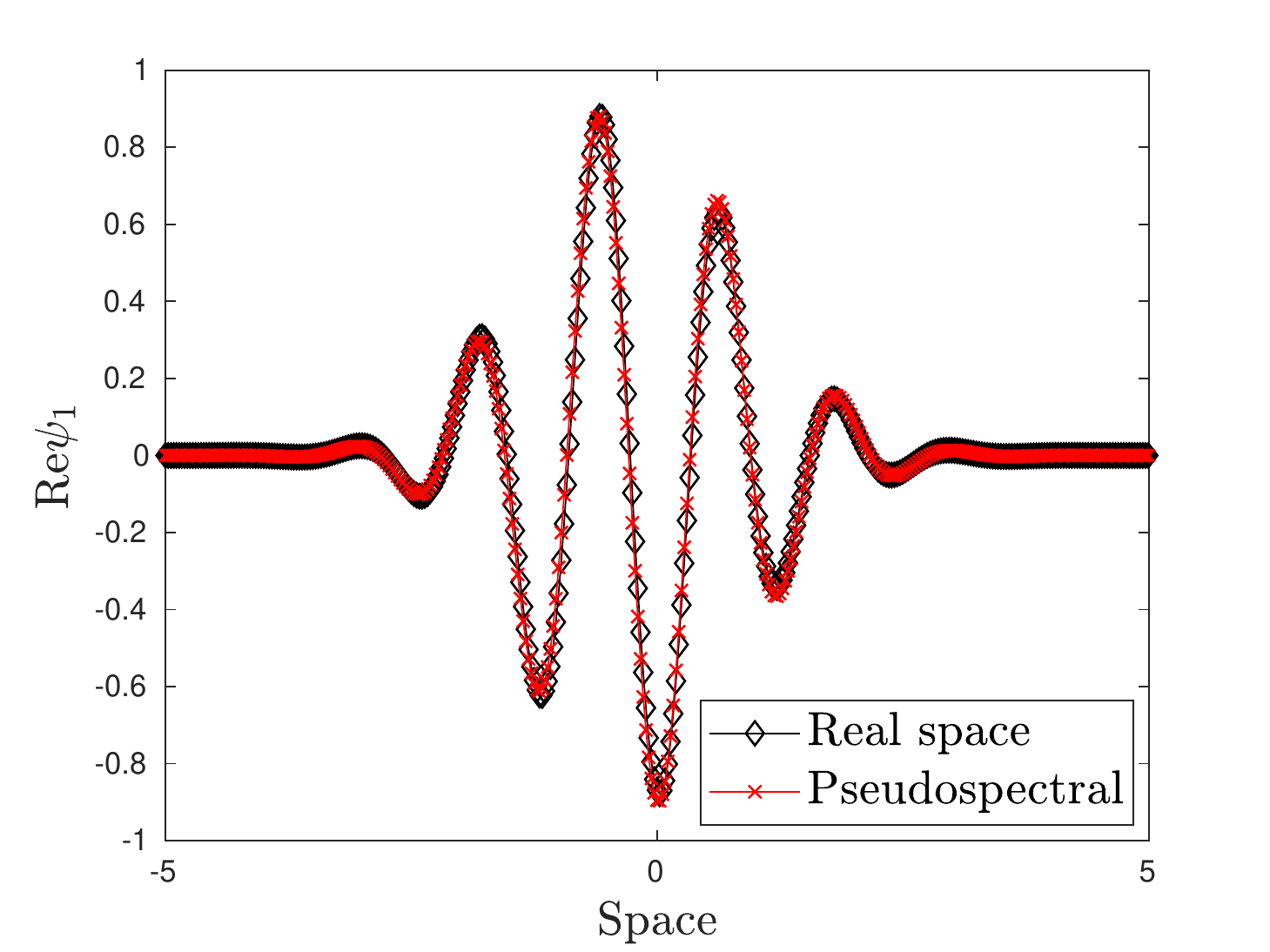}
\includegraphics[height=4cm,keepaspectratio]{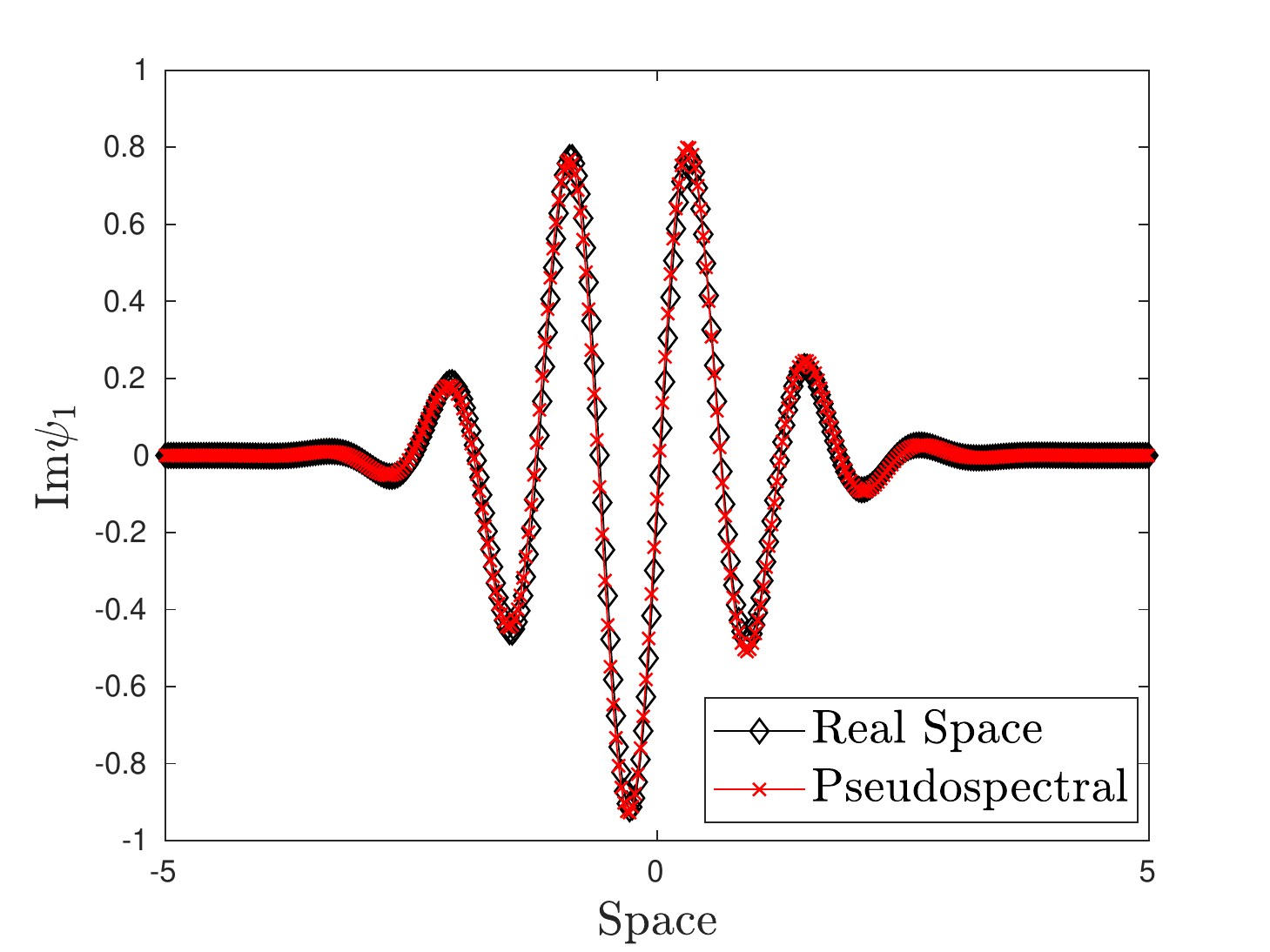}
\end{center}
\caption{{\bf Experiment 1.} (Left) Velocity field $x \mapsto \exp\big(\Phi(x)-\Psi(x)\big)$. (Middle) Real part of  $\psi_1(\cdot,T)$ with real space and pseudospectral methods. (Right) Imaginary part of $\psi_1(\cdot,T)$ with real space and pseudospectral methods. The
final time is $T=0.5$.}
\label{figvelo}
\end{figure}

\begin{figure}
\begin{center}
\includegraphics[height=4cm,keepaspectratio]{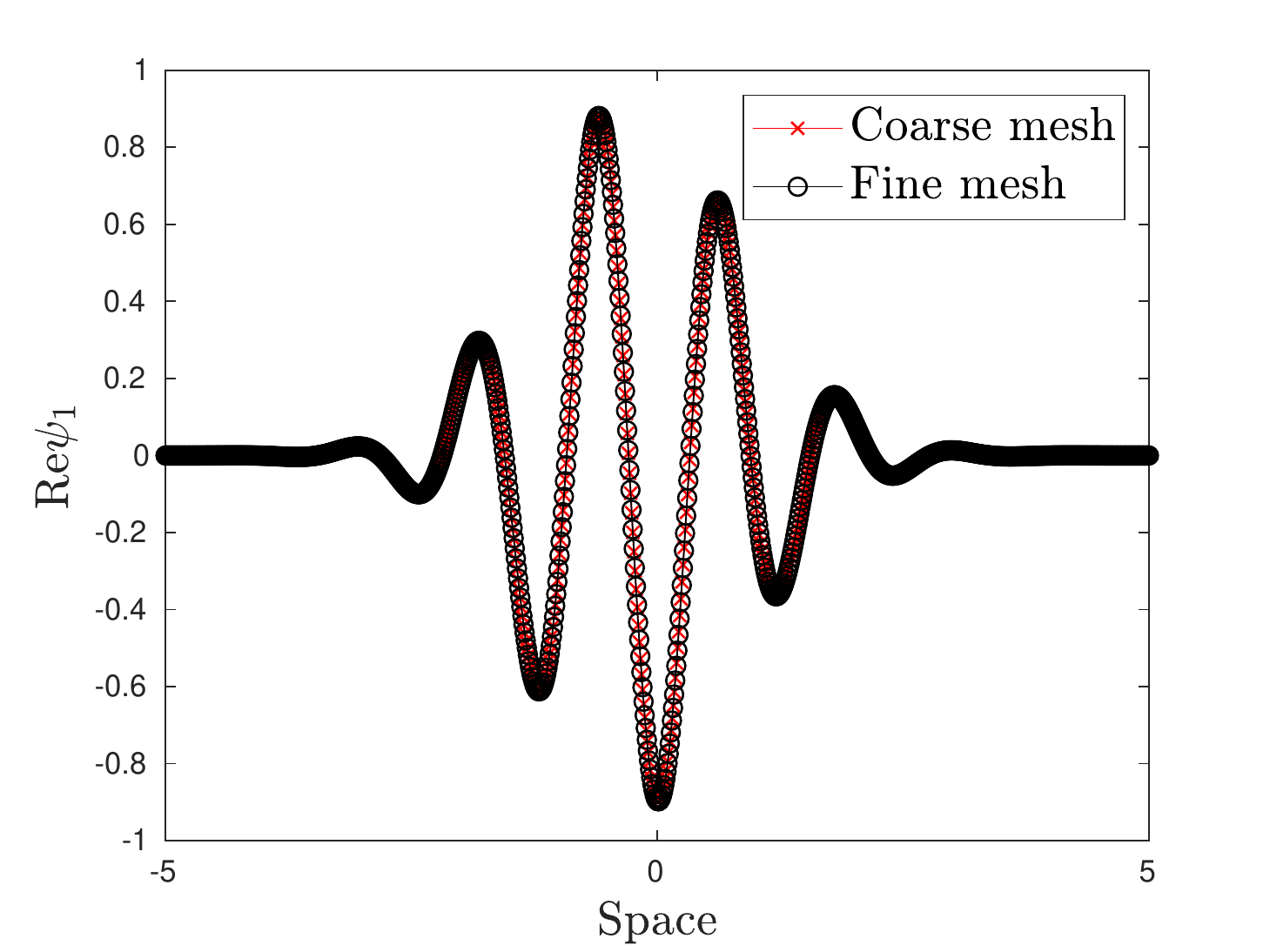}
\includegraphics[height=4cm,keepaspectratio]{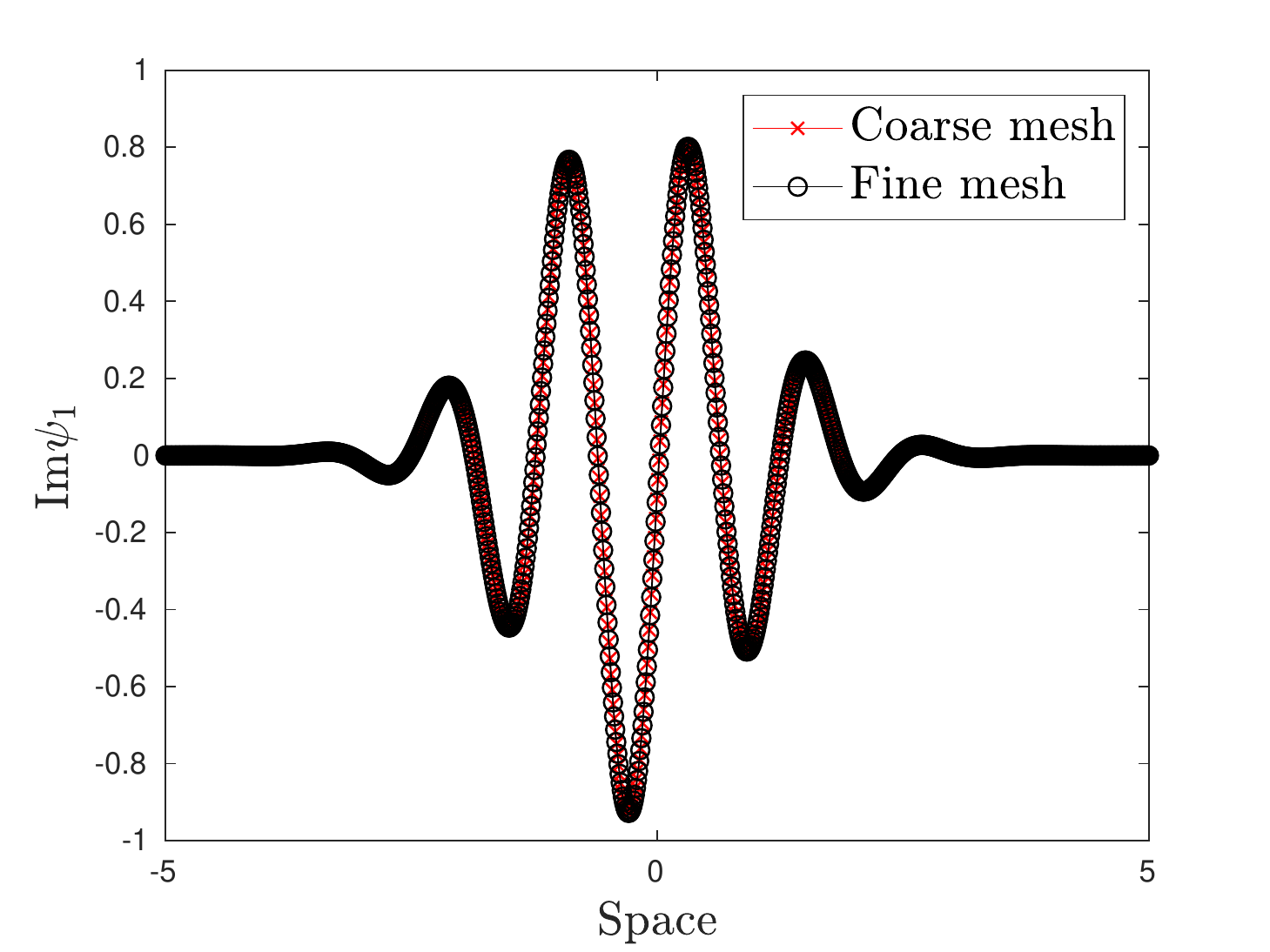}
\includegraphics[height=4cm,keepaspectratio]{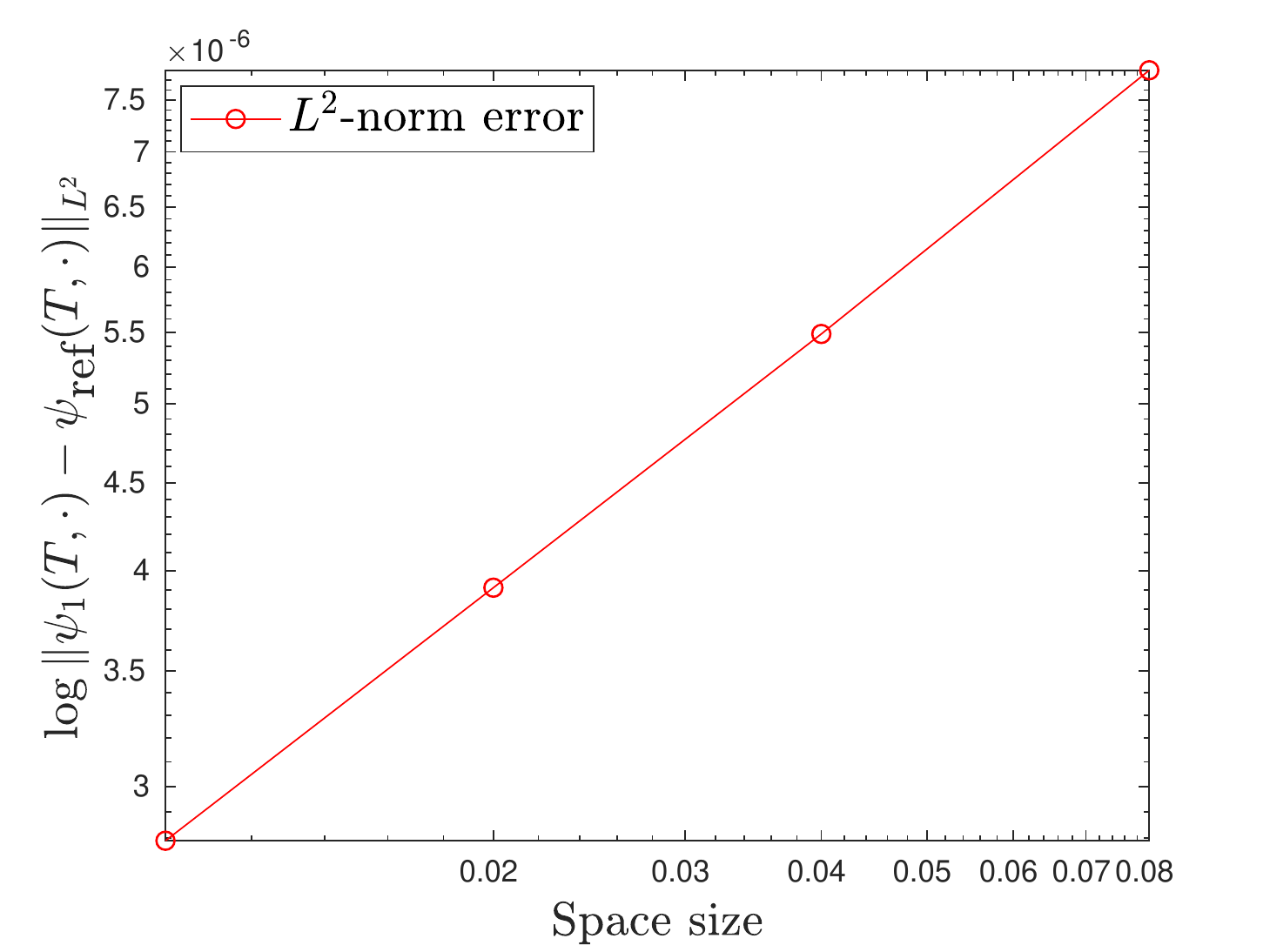}
\end{center}
\caption{{\bf Numerical Experiments 1.} (Left) Real part of  $\psi_1(\cdot,T)$ on coarse and fine grids. (Middle) Imaginary part of  $\psi_1(\cdot,T)$ on coarse and fine grids. (Right) $\ell^2$-norm error.}
\label{figveloB}
\end{figure}

We consider the metric $ds^2=e^{2\Phi(x)}dt^2-e^{2\Psi(x)}dx^2$, such that $\Phi$ and $\Psi$ are two space-dependent functions. This leads to the following one-dimensional Dirac equation \cite{Koke}
\begin{eqnarray*}
{\tt i}\partial_t\psi & = &  -{\tt i}e^{\Phi(x)-\Psi(x)}\sigma_x\Big(\partial_x + \cfrac{\Phi'(x)}{2}\Big)\psi + e^{\Phi}\sigma_zm\psi \, .
\end{eqnarray*}
\noindent{\bf Numerical Experiment 1.}
We propose a benchmark with $\Psi(x)=e^{-10^{-2}x^2}$, $\Phi(x)=e^{-5\times 10^{-3}x^2}$, and for
 $\phi_0(x)=\exp(-x^2/2+{\tt i}k_0x)$, where $k_0=5$ and  $c=1$. The computational domain is $[-5,5]$, 
 while the discretization parameters are set through: $\Delta t=5 \times 10^{-4}$ and $N_1=18027$. We compare the proposed method, with a real-space method, at $\mathrm{CFL}=0.99$, which degenerates into the Quantum-Boltzmann method for flat space \cite{cpc2012}. The second-order splitting implicit pseudospectral method ({\bf Numerical Scheme I}) is implemented by using GMRES \cite{saad}.
 We plot $x \mapsto \exp(\Phi(x)-\Psi(x))$, corresponding to a velocity field, in Fig. \ref{figvelo} (Left). We report the real and imaginary parts of the first component $\psi_1$ in Fig. \ref{figvelo}  (Middle, Right) at time $T=0.5$, corresponding to $1000$ time iterations.
Unlike  the real space method for $\mathrm{CFL}=1$, the pseudospectral method is {\it linearly stable}. As an illustration, we compare in Fig. \ref{figveloB} (Left, Middle) the real and imaginary parts of $\psi_1$ on a coarse grid ($h=1.1\times 10^{-2}$, $\Delta t=10^{-2}$) and fine grid  ($h=5.5\times 10^{-4}$, $\Delta t=5 \times 10^{-4}$). Finally in Fig. \ref{figveloB} (Right), we report in logscale the $\ell^2-$ norm of the error ($\|\psi_{h}(T,\cdot)-\psi_{\textrm{ref},h}(T,\cdot)\|_{\ell^2}$) as a function of the space size $h$.\\

\begin{figure}
\begin{center}
\includegraphics[height=4cm,keepaspectratio]{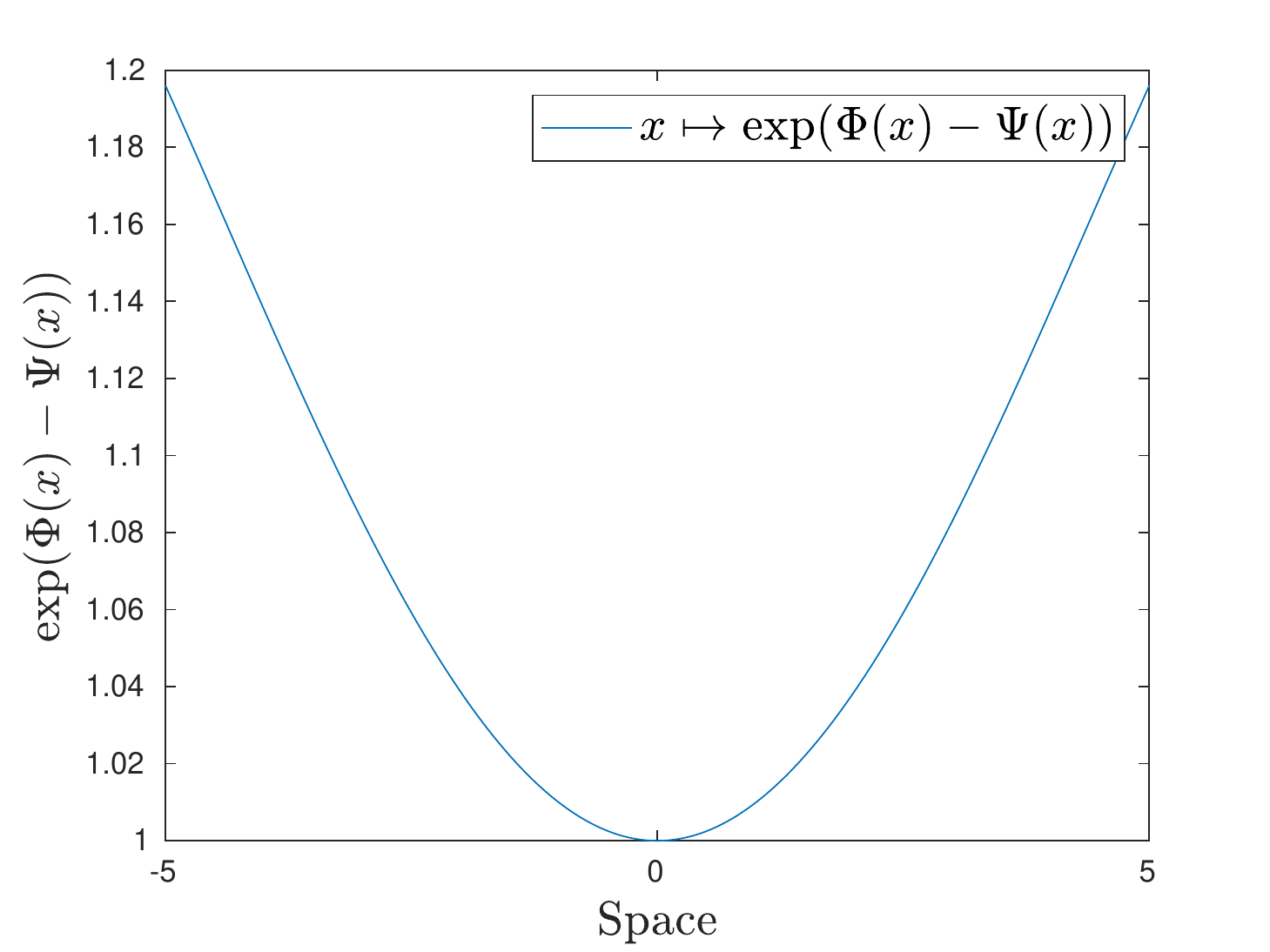}
\includegraphics[height=4cm,keepaspectratio]{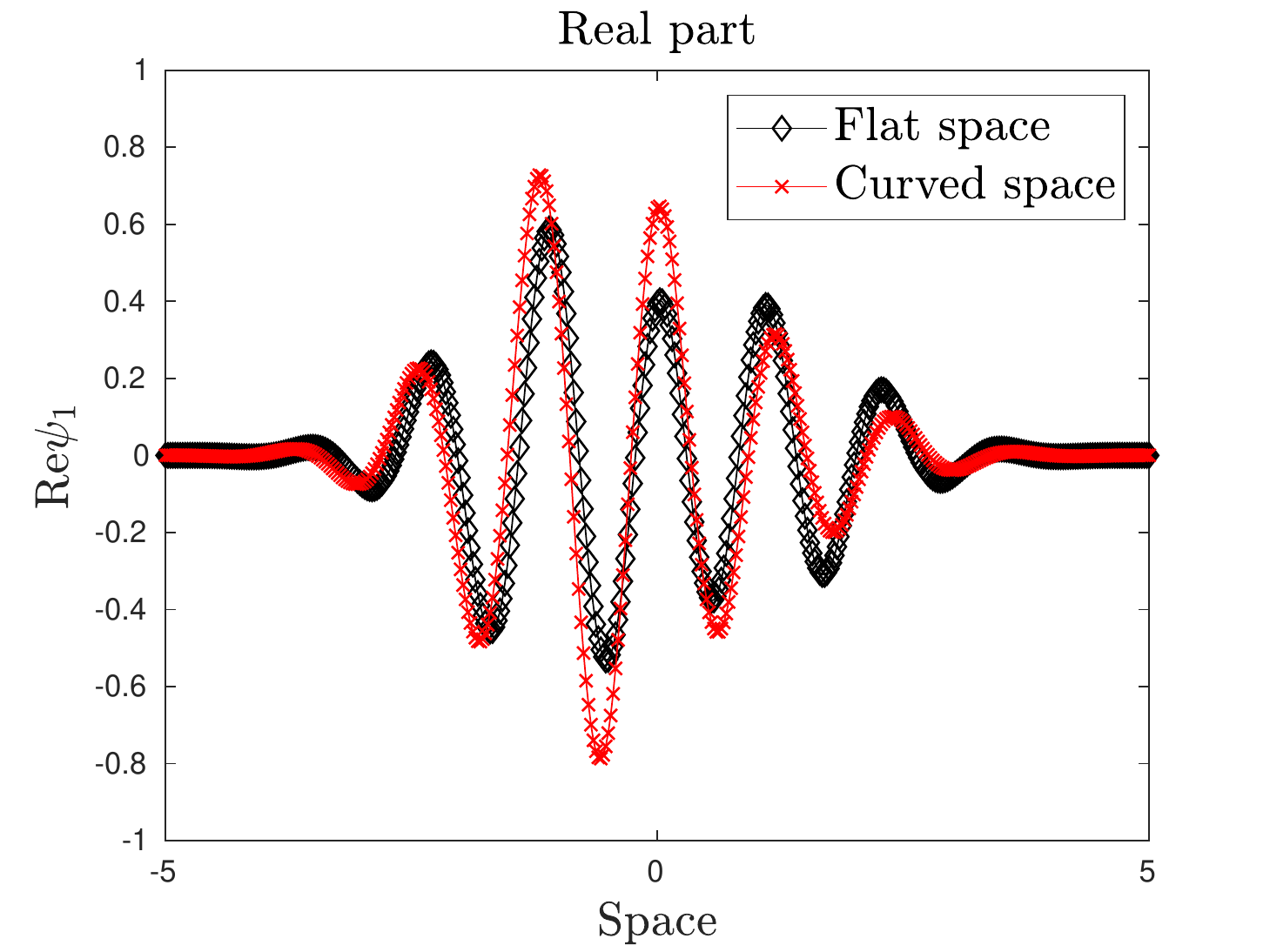}
\includegraphics[height=4cm,keepaspectratio]{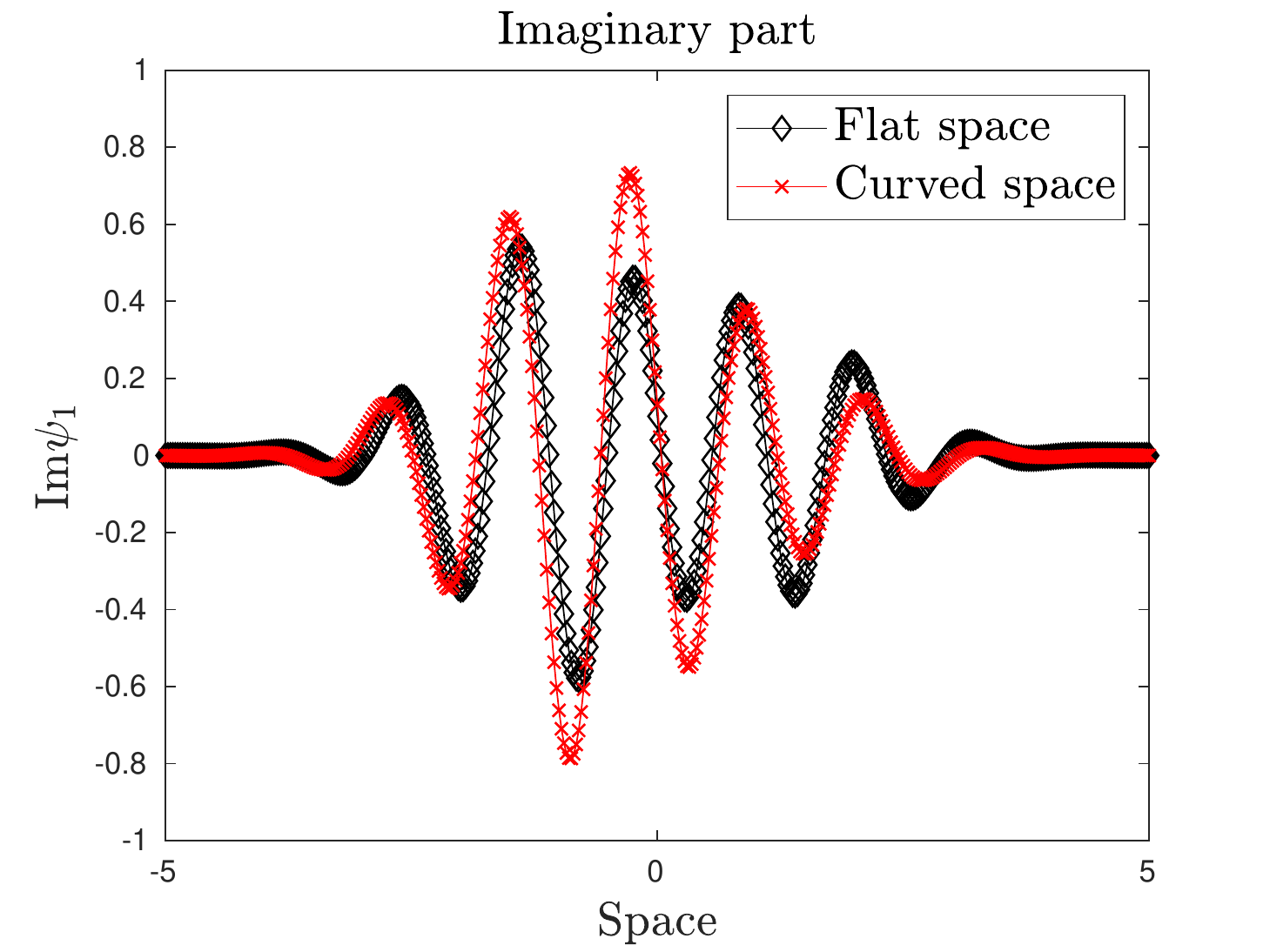}
\end{center} 
\caption{{\bf Numerical Experiments 2.} (Left) Velocity field $x \mapsto \exp\big(\Phi(x)-\Psi(x)\big)$. (Middle) Real part of $\psi_1(\cdot,T)$ for flat and curved spaces. (Right) Imaginary part of  $\psi_1(\cdot,T)$ for flat and curved spaces. The final time is $T=1$.}
\label{figvelo2}
\end{figure}

\noindent{\bf Numerical Experiments 2.} We compare the pseudospectral method in flat and curved spaces. For the curved space case, we select $\Psi(x)=\cos(x/10)e^{-10^{-2}x^2}$, $\Phi(x)=e^{-10^{-2}x^2}$, with $\phi_0(x)=\exp(-x^2/2+{\tt i}k_0x)$, where $k_0=5$ and we again take $c=1$. The computational domain is $[-5,5]$, with discretization parameters $\Delta t=5 \times 10^{-4}$ and $N_1=20001$. The second order implicit splitting pseudospectral method is again solved by using GMRES \cite{saad}.  We plot $x \mapsto \exp(\Phi(x)-\Psi(x))$, corresponding to the velocity field in Fig. \ref{figvelo2} (Left), and the real and imaginary parts of the first component $\psi_1$ in Fig. \ref{figvelo2}  (Middle, Right) at time $T=1$, corresponding to $2000$ time iterations. This illustrates the effect of the spatial curvature on the solution to the Dirac equation.\\

\begin{figure}
\begin{center}
\includegraphics[height=6cm,keepaspectratio]{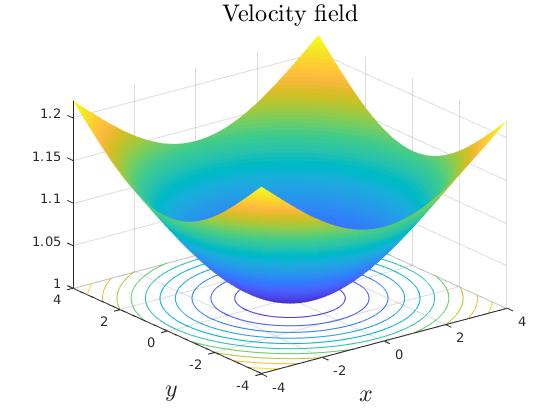}
\includegraphics[height=6cm,keepaspectratio]{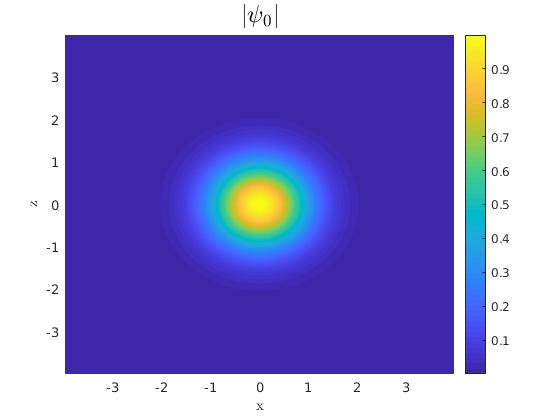}
\end{center}
\caption{{\bf Numerical Experiments 3.} (Left) Velocity field. (Right) Initial wavefunction $\psi_1(0,\cdot)$.}
\label{velo2D}
\end{figure}

\noindent{\bf Numerical Experiments 3.} 
We now consider a two-dimensional Dirac equation in curved space, defined by the metric $d{\boldsymbol s}^2=e^{2\Phi({\boldsymbol x})}dt^2-e^{2\Psi({\boldsymbol x})}d{\boldsymbol x}^2$, such that $\Phi$, $\Psi$ are two space-dependent functions, with ${\boldsymbol x}=(x,y)$,
\begin{eqnarray*}
{\tt i}\partial_t\psi & = &  -{\tt i}e^{\Phi({\boldsymbol x})-\Psi({\boldsymbol x})}\Big(\sigma_x\Big(\partial_x + \cfrac{\partial_x\Phi({\boldsymbol x})}{2}\Big) + \sigma_y\Big(\partial_y + \cfrac{\partial_y\Phi({\boldsymbol x})}{2}\Big) \Big)\psi + e^{\Phi({\boldsymbol x})}\sigma_zm\psi \, .
\end{eqnarray*}
We assume that $\Phi({\boldsymbol x})=e^{-10^{-2}\|{\boldsymbol x}\|^2}$, $\Psi({\boldsymbol x})=e^{-5\times 10^{-3}\|{\boldsymbol x}\|^2}$, with $\phi_0({\boldsymbol x})=e^{-\|{\boldsymbol x}\|^2/2+{\tt i}{\boldsymbol k}_0 \cdot {\boldsymbol x}}$, where ${\boldsymbol k}_0=(5,5)^T$ and setting $c=1$. The computational domain is $[-5,5]^2$.
We report in Fig. \ref{velo2D} (Left) the velocity field on $\mathcal{D}=[-5,5]^2$, ${\boldsymbol x} \mapsto \exp\big(\Phi({\boldsymbol x}) - \Psi({\boldsymbol x})\big)$.\\

\begin{figure}
\begin{center}
\includegraphics[height=4.8cm,keepaspectratio]{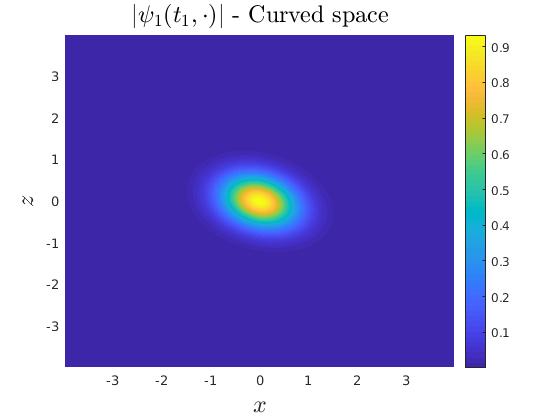}
\includegraphics[height=4.8cm,keepaspectratio]{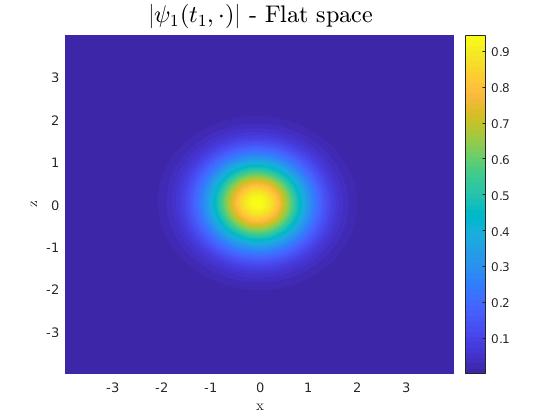}
\includegraphics[height=4.8cm,keepaspectratio]{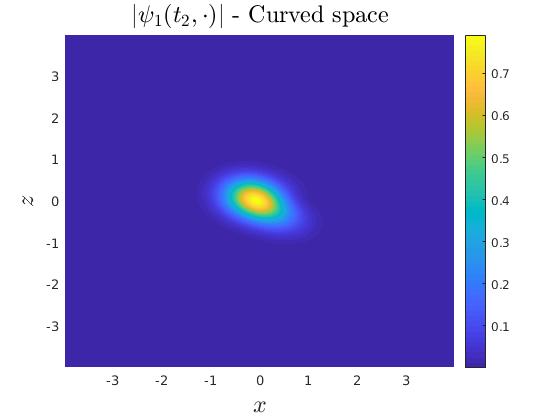}
\includegraphics[height=4.8cm,keepaspectratio]{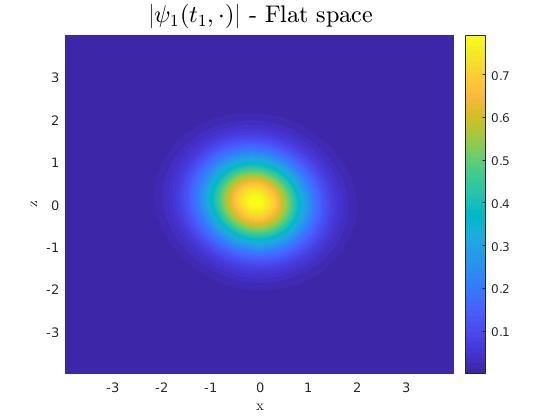}
\includegraphics[height=4.8cm,keepaspectratio]{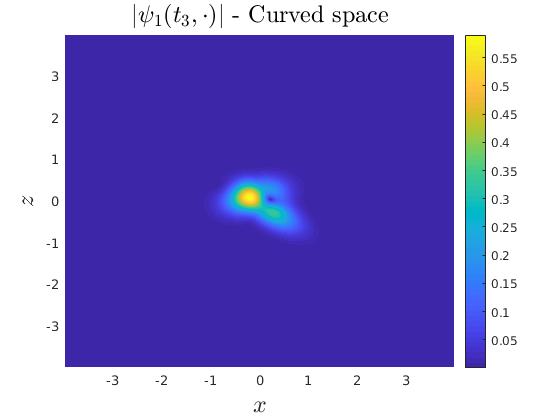}
\includegraphics[height=4.8cm,keepaspectratio]{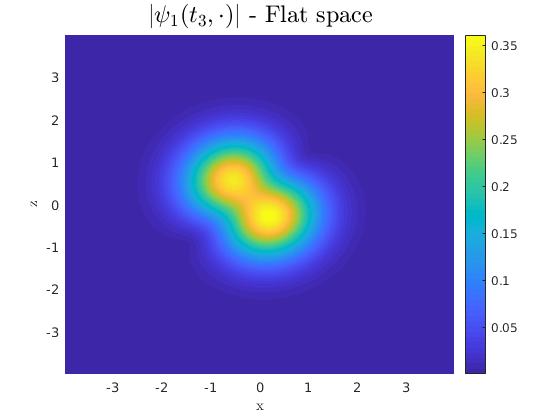}
\includegraphics[height=4.8cm,keepaspectratio]{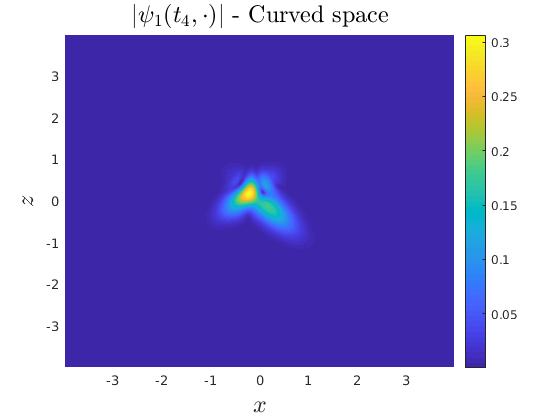}
\includegraphics[height=4.8cm,keepaspectratio]{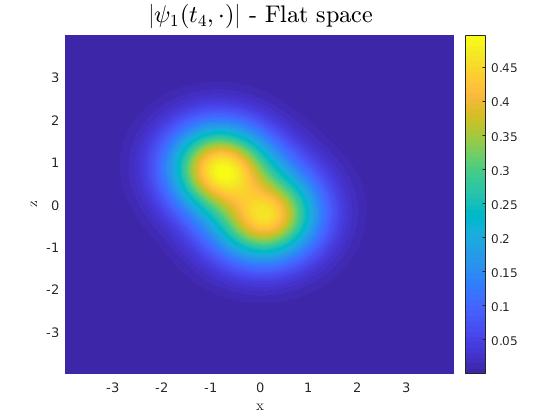}
\end{center}
\caption{{\bf Numerical Experiments 3.} Modulus of the first component of the 4-spinor (Left) Curved space: at times $t_1=0.57\times 10^{-2}$, $t_2=1.14\times 10^{-2}$, $t_3=2.28\times 10^{-2}$ and $t_4=4.56\times 10^{-2}$. (Right) Flat space: at times $t_1=1.14\times 10^{-2}$, $t_2=2.28\times 10^{-2}$ and $t_4=4.56\times 10^{-2}$. }
\label{sol2D}
\end{figure}

\begin{figure}
\begin{center}
\includegraphics[height=4.8cm,keepaspectratio]{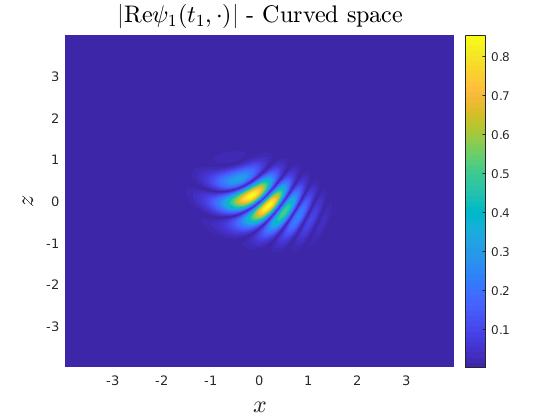}
\includegraphics[height=4.8cm,keepaspectratio]{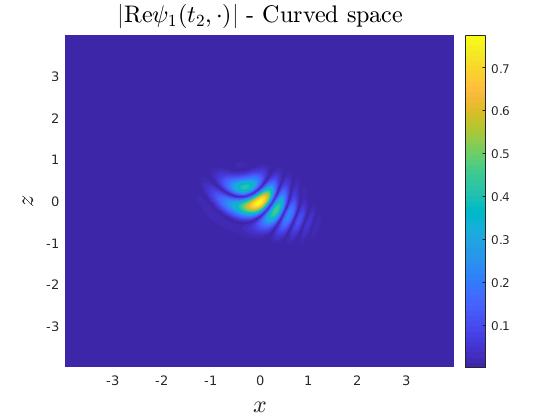}
\includegraphics[height=4.8cm,keepaspectratio]{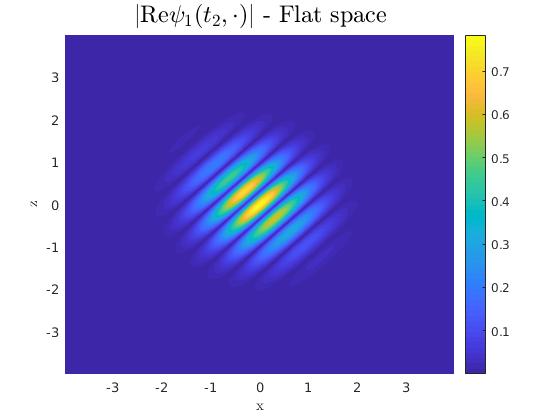}
\includegraphics[height=4.8cm,keepaspectratio]{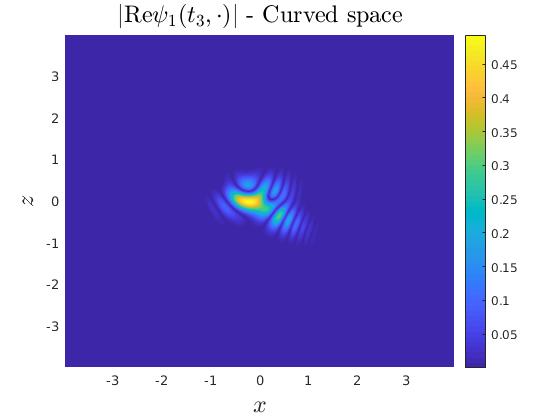}
\includegraphics[height=4.8cm,keepaspectratio]{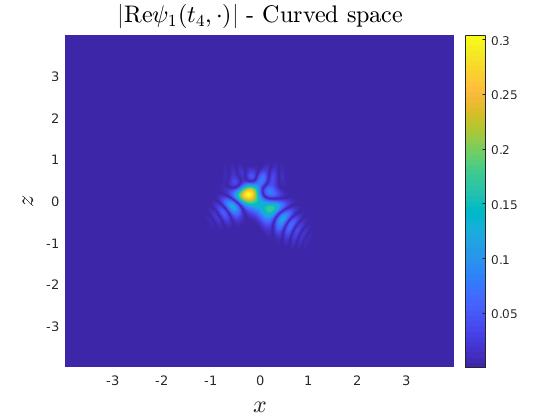}
\end{center}
\caption{{\bf Numerical Experiments 3.} Modulus of the first real part of the component of the 4-spinor (Left) Curved space: at times $t_1=0.57\times 10^{-2}$, $t_2=1.14\times 10^{-2}$, $t_3=2.28\times 10^{-2}$ and $t_4=4.56\times 10^{-2}$. (Right) Flat space: at times $t_1=1.14\times 10^{-2}$, $t_2=2.28\times 10^{-2}$ and $t_4=4.56\times 10^{-2}$. }
\label{sol2Dbis}
\end{figure}

We here implement the {\bf Numerical scheme II}. The numerical data are as follows: $\Delta t=1.14\times 10^{-4}$, and $N_1=N_z=512$. The initial data is a wavepacket
\begin{eqnarray*}
\psi_0(x,z) = \big(\phi_1(x,z),0,0,0\big)^T,
\end{eqnarray*}
where $\phi_1(x,y)=e^{-(x^2+y^2)/2 + 5{\tt i}(x-z)}$, which is plotted in Fig. \ref{velo2D} (Right). We report in Fig. \ref{sol2D} (resp. Fig. \ref{sol2Dbis}) the modulus of the first component (resp. real part of the first component) of the Dirac 4-spinor at times (in {\it atomic units}) $t_1=0.57\times 10^{-2}$, $t_2=1.14\times 10^{-2}$, $t_3=2.28\times 10^{-2}$ and $t_4=4.56\times 10^{-2}$ in the flat (Left) and curved (Right) spaces.

\subsection{Massless Dirac particles in curved graphene}

Charge carriers in graphene can be described theoretically by a 2-D Dirac equation in curved space-time, where the metric is related to the graphene sample deformation \cite{Cortijo_2007}. The dynamics of these charge carriers was recently studied numerically \cite{PhysRevB.98.155419,Debus_2018} with a lattice Boltzmann method. Using our numerical schemes, we now propose some tests with similar configurations.

\begin{figure}
\begin{center}
\includegraphics[height=6cm,keepaspectratio]{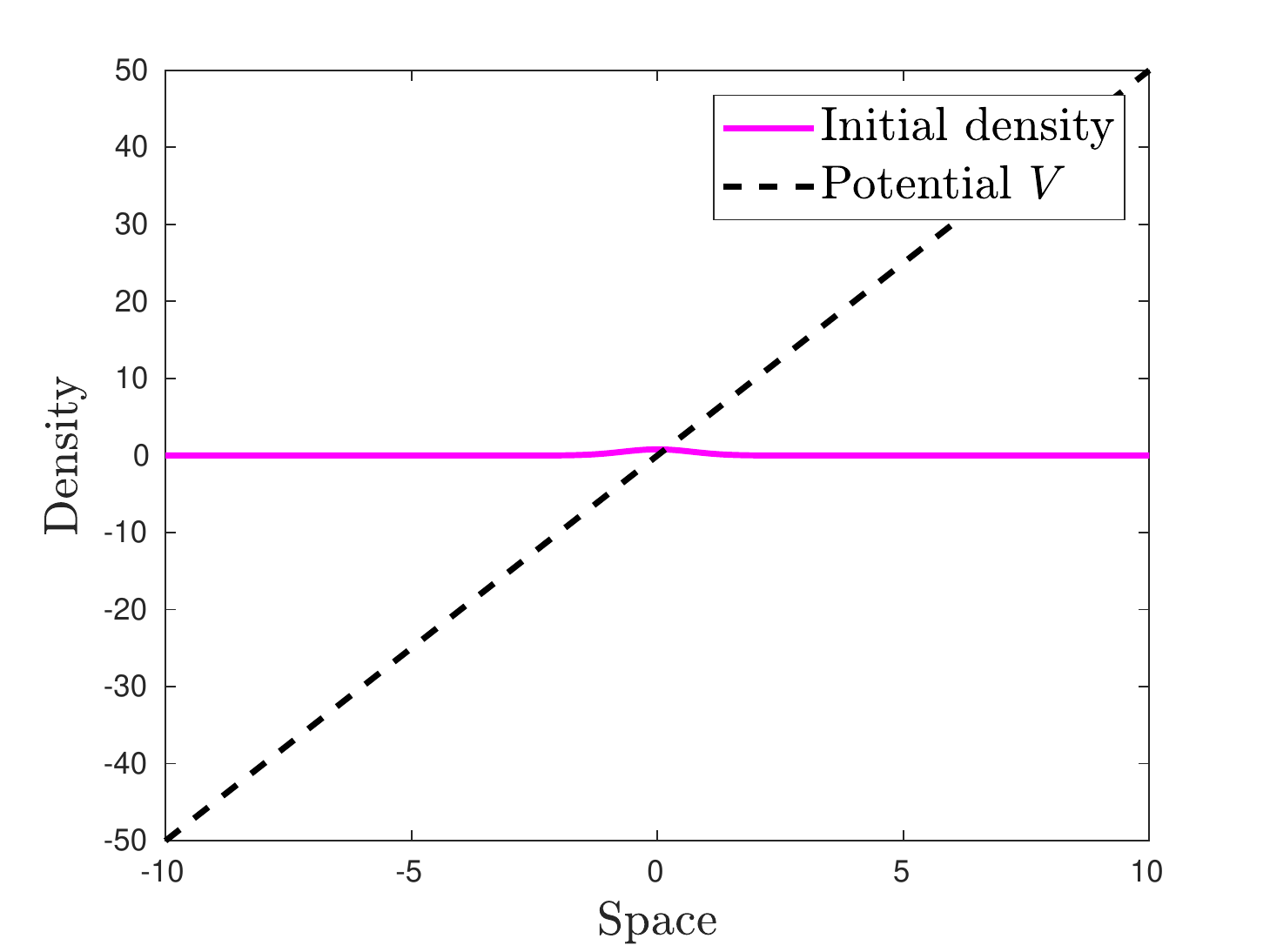}
\includegraphics[height=6cm,keepaspectratio]{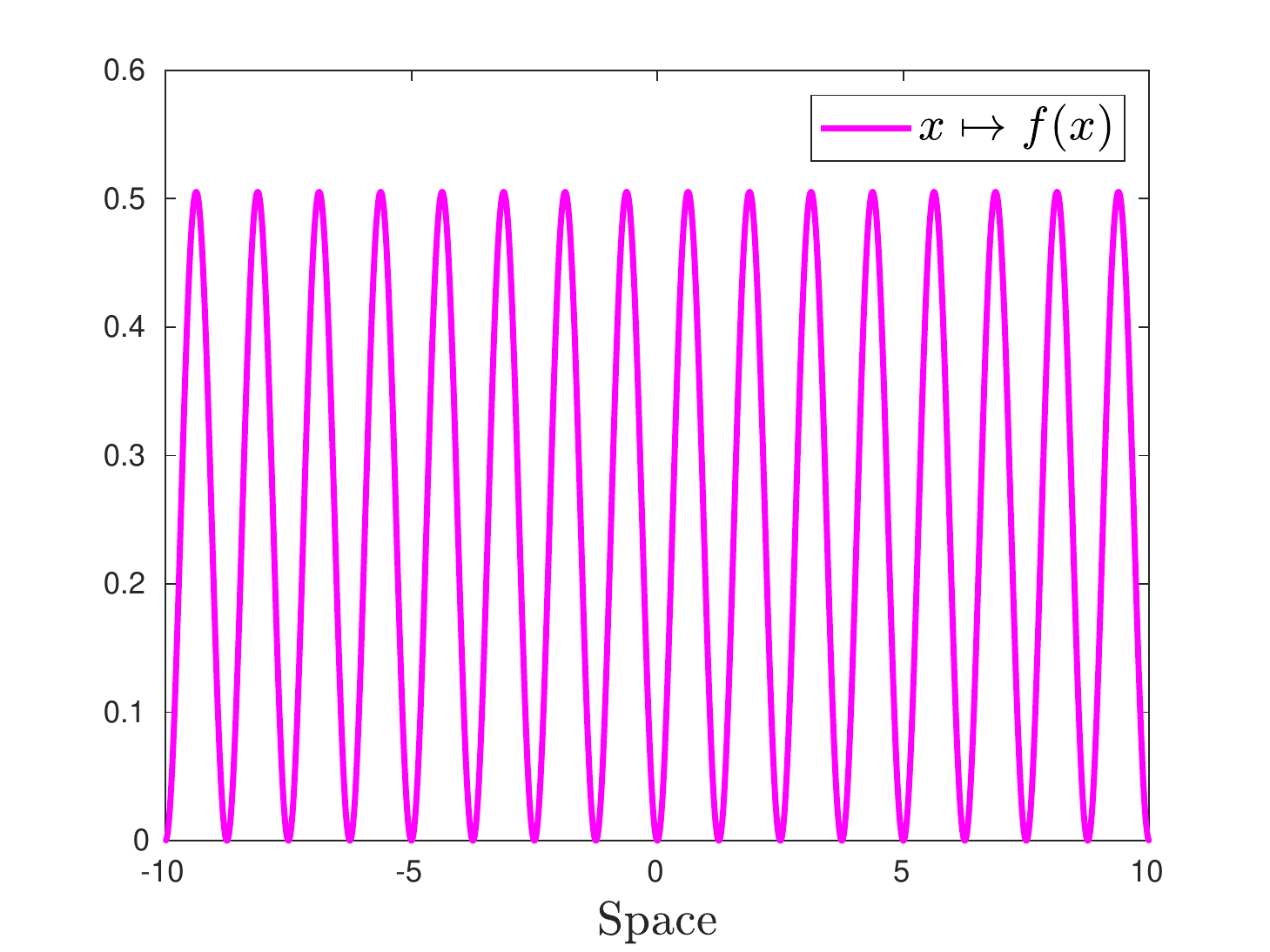}
\end{center}
\caption{{\bf Experiment 4.} (Left) Initial density (Right) Graph of $f$ and $|G|$.}
\label{figV0}
\end{figure}

\begin{figure}
\begin{center}
\includegraphics[height=6cm,keepaspectratio]{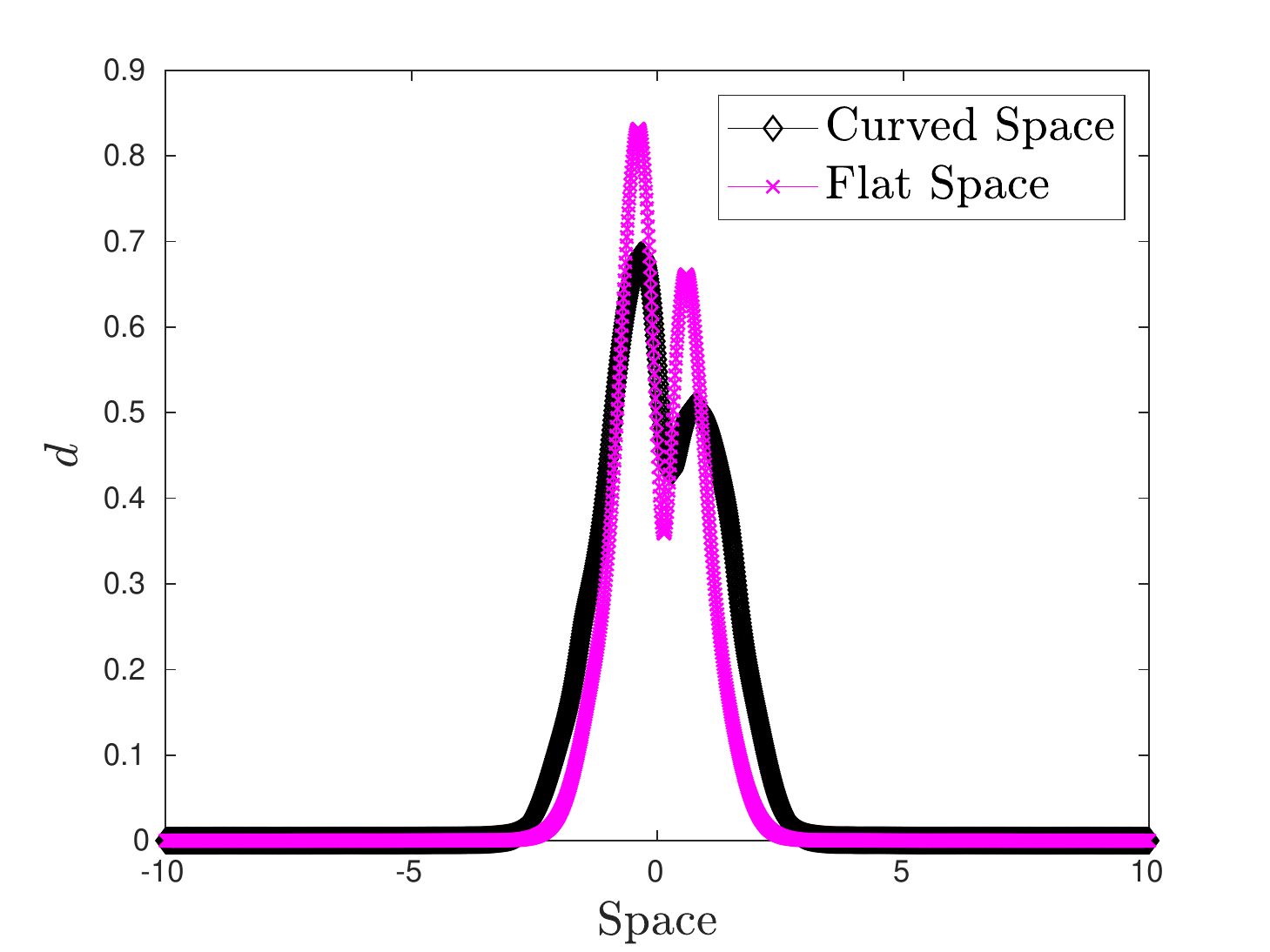}
\includegraphics[height=6cm,keepaspectratio]{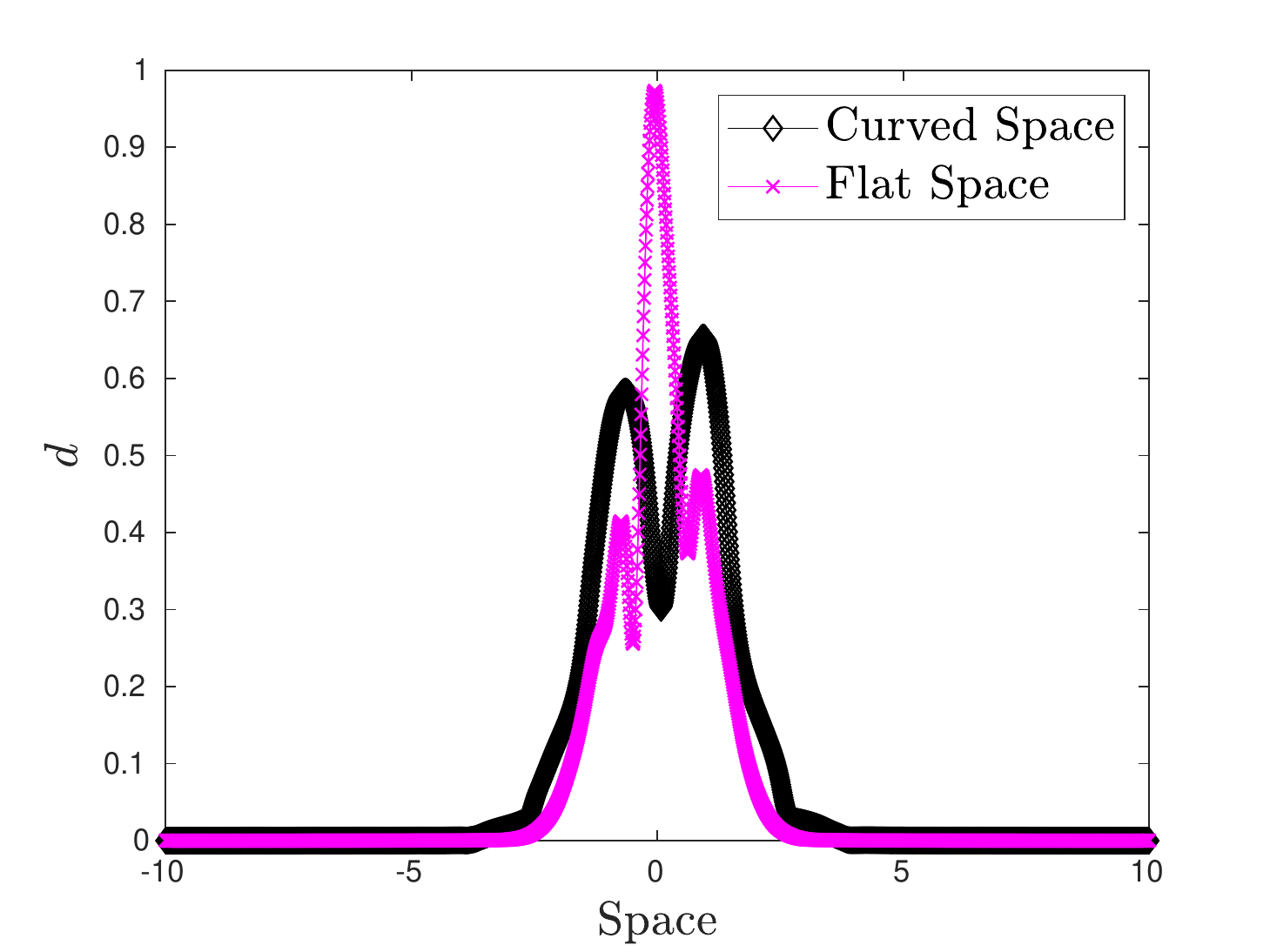}
\includegraphics[height=6cm,keepaspectratio]{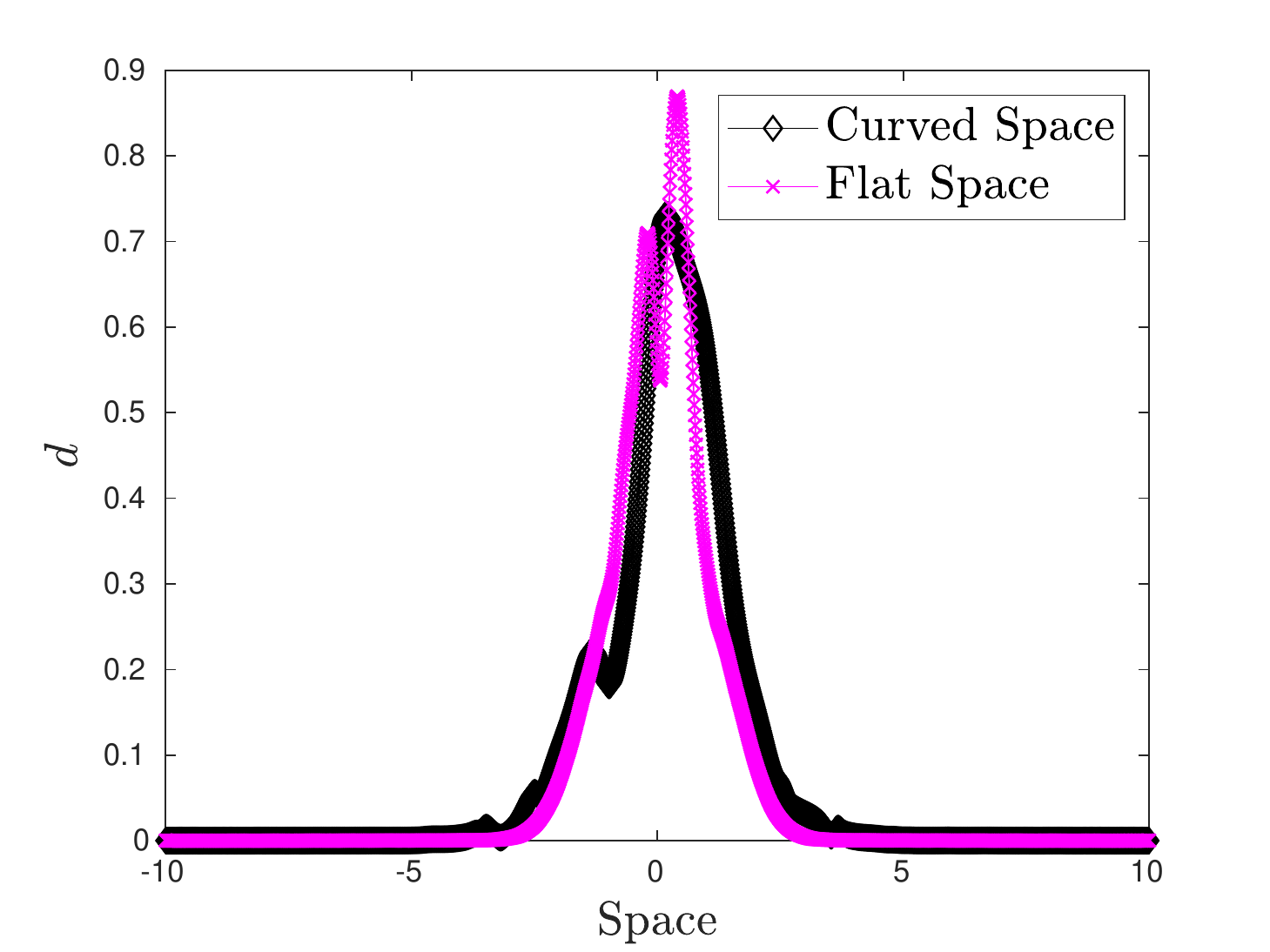}
\includegraphics[height=6cm,keepaspectratio]{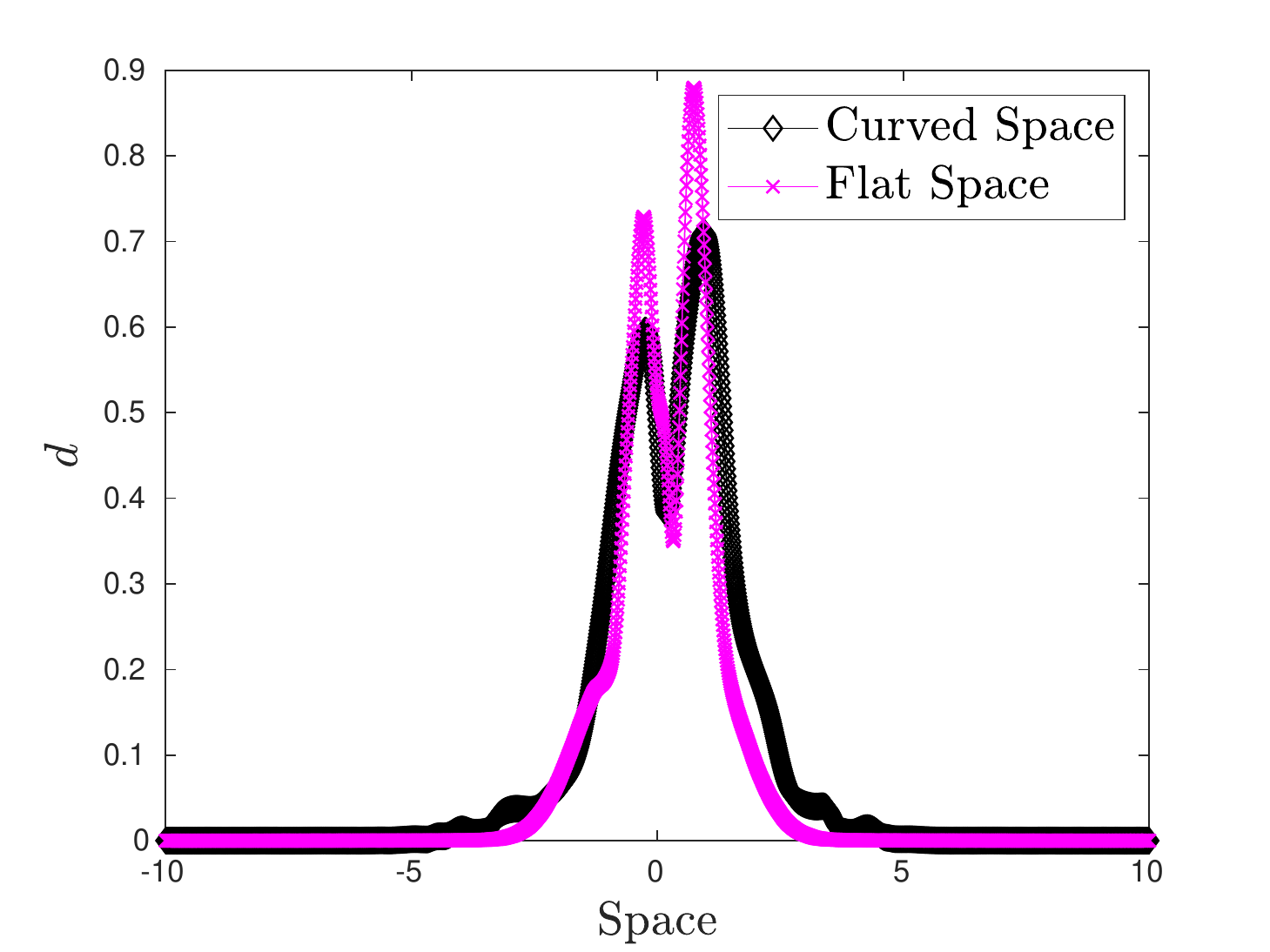}
\end{center}
\caption{{\bf Experiment 4.} From Top-Left to Bottom-Right: density in flat and curved spaces at time $t=0.4$, $t=0.8$, $t=1.2$ and $t=1.6$.}
\label{compDen0}
\end{figure}

\begin{figure}
\begin{center}
\includegraphics[height=6cm,keepaspectratio]{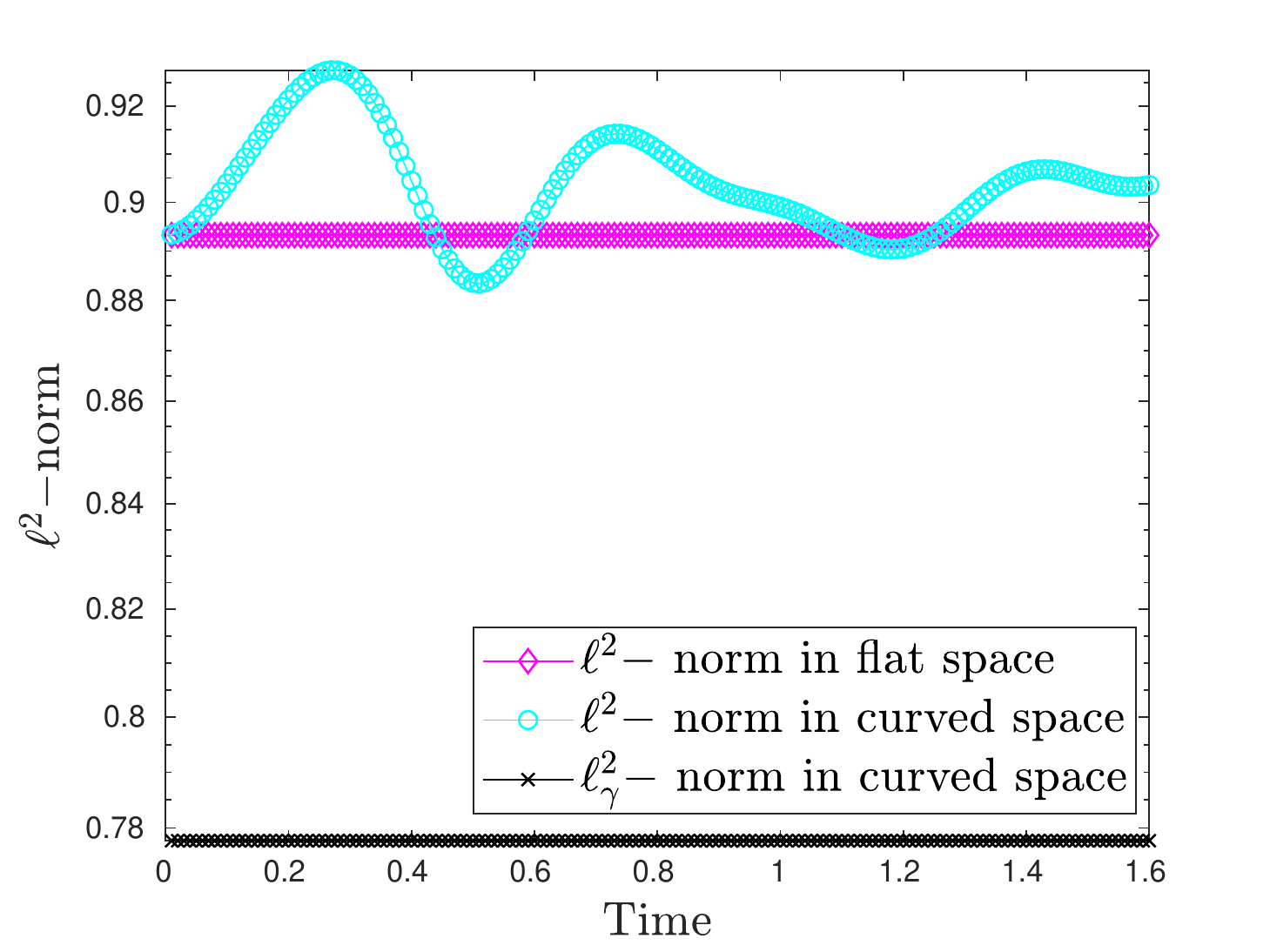}
\includegraphics[height=6cm,keepaspectratio]{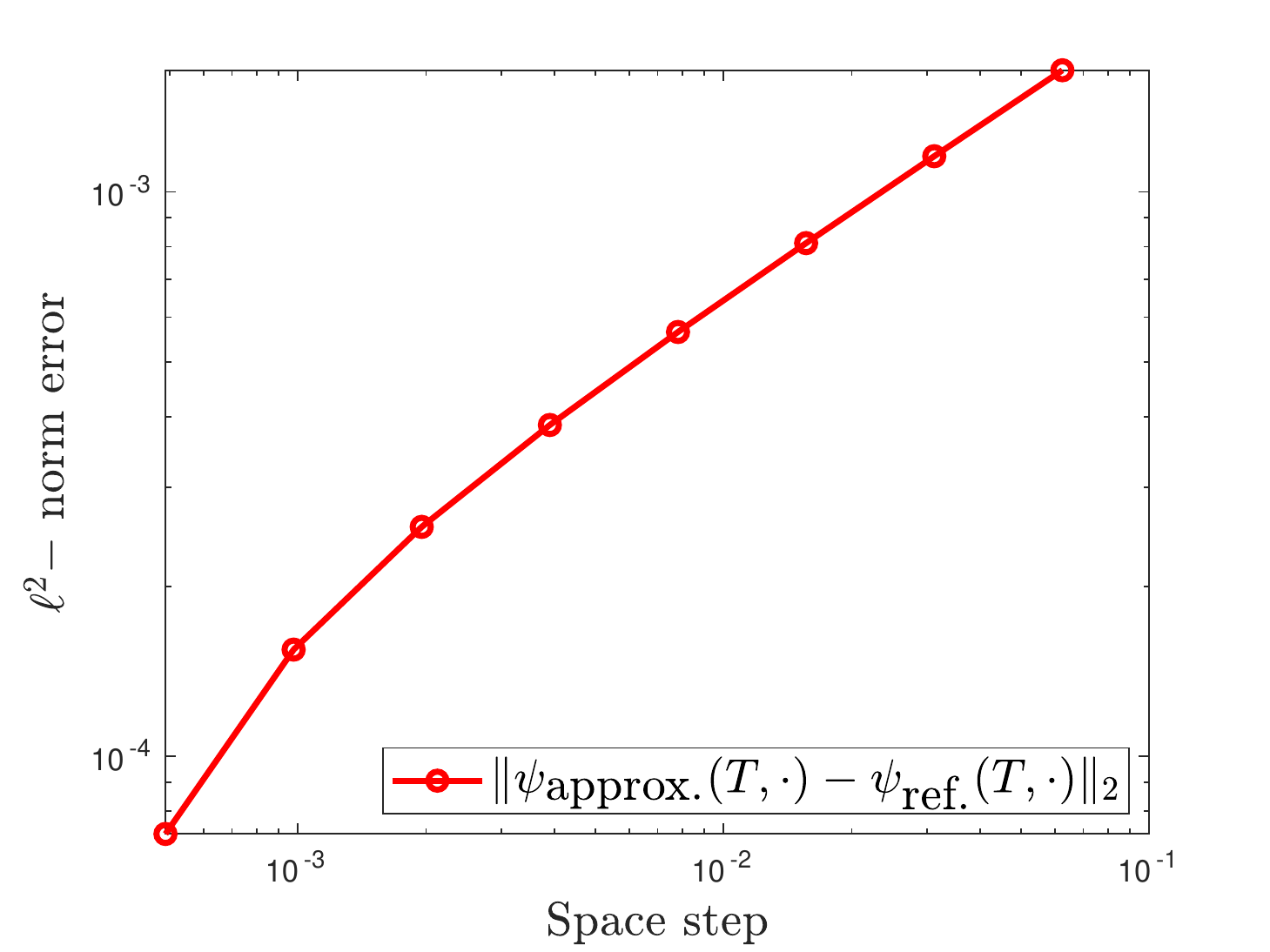}
\end{center}
\caption{{\bf Experiment 4. }(Left) $\ell^2-$norm and $\ell_{\gamma}^2-$norm in logscale of the solution as a function of time iterations in flat and curved space. (Right) Graph of convergence in logscale $(h,\|\psi_{\textrm{approx.}}(T,\cdot)-\psi_{\textrm{ref.}}(T,\cdot)\|_2$ with $h=1/2,\cdots,1/2^{10}$.}
\label{l2norm0}
\end{figure}

\noindent{\bf One-dimensional test: rippled graphene sample}. This test is dedicated to the simulation of strained graphene, and a simple comparison with non-strained graphene (corresponding to  flat space). More specifically, we consider a rippled graphene sheet parameterized by the coordinate transformation map \cite{PhysRevB.98.155419}.  We set $h(x)=a_0\cos(2\pi k_0x/\ell)$ and $f(x)=\big(h'(x)\big)^2/2 = 2\pi^2a_0^2k_0^2\sin\big(2\pi k_0 x/\ell\big)^2/\ell^2$, where $a_0$ and $k_0$ denote the amplitude and wave vector of the surface ripples, and $\ell$ is the length of the sheet. Moreover, in one-dimension, the spatial part of the metric and the tetrad are given by
\begin{eqnarray*}
g(x) & = & \big(1-f(x)\big)^2, \qquad  e_x^x(x) = \cfrac{1}{1-f(x)} \, .
\end{eqnarray*}
 The Dirac equation for modeling strained graphene with external electromagnetic potentials $(A,V)$ in 1-d then reads 
\begin{eqnarray}\label{eqC}
\partial_t \psi + \sigma_xe^x_x\big(\partial_x-{\tt i}A_x\big)\psi & = &  -{\tt i}\gamma^0(m-V)\psi \, .
\end{eqnarray}
This equation can be solved using the numerical schemes developed in this article. 
It is easy to show that the following $\ell^2_{\gamma}-$norm is conserved
\begin{eqnarray*}
\|\psi\|_{\gamma} & := & \Big(\big(1-f(x)\big)|\psi(x)|^2dx\Big)^{1/2} \, .
\end{eqnarray*}
Indeed,  we multiply \eqref{eqC} by $(1-f(x))\psi^{\dagger}$, take the real part and integrate in space, and directly get:
\begin{eqnarray*}
\cfrac{d}{dt}\int \big(1-f(x)\big)|\psi(x)|^2dx & = & 0 \, .
\end{eqnarray*}
\noindent{\bf Experiment 4.} In this first experiment,  we assume that the initial data is
\begin{eqnarray*}
\psi(0,x) = (1,{\tt i})^T\beta e^{-\beta x^2/2}/\sqrt{4\pi } \, .
\end{eqnarray*}
Numerically, we take $\beta=2$, $a_0=4\times 10^{-1}$, $k_0=2$, $c=1$, and $\ell=5$. Moreover, we fix $A_x(x)=V(x)=5x$. We plot in Fig. \ref{figV0} (Right), the graph of $x\in [-10,10] \mapsto f(x)$. Numerically, we choose the discretization parameters
to be $\Delta t=10^{-2}$ and $h=10^{-2}$. We report in Fig. \ref{compDen0} the density, defined by $d_F(t,\cdot)=|\psi_1(t,\cdot)|^2+|\psi_2(t,\cdot)|^2$ in flat space, and the density $d_C(t,\cdot)$ in curved space at different times $t=0.4$, $t=0.8$, $t=1.2$ and $t=1.6$.  In Fig. \ref{l2norm0} (Left), we report in logscale the $\ell^2-$ and $\ell_{\gamma}^2-$norms of the solution as a function of time iterations, in flat and curved space as well as illustrating the $\ell^2-$stability, and the $\ell^2-$norm conservation in flat space, and $\ell^2_{\gamma}-$norm conservation in curved space.

We also report the graph of convergence in Fig. \ref{l2norm0} (Right.), with $\Delta t=10^{-5}$, and $T=10^{-1}$ (corresponding to $1000$ time-iterations) and the computational domain is still $[-10,10]$. We represent the $\ell^2$-norm error between a solution  of reference $\psi_{\textrm{ref.}}(\cdot,T)$ (computed on a very fine mesh) and the approximate one $\psi_{\textrm{approx.}}(\cdot,T)$, computed with a mesh-size $h$ of $1/2^{i}$ with $i=4,...,11$.\\

\begin{figure}
\begin{center}
\includegraphics[height=6cm,keepaspectratio]{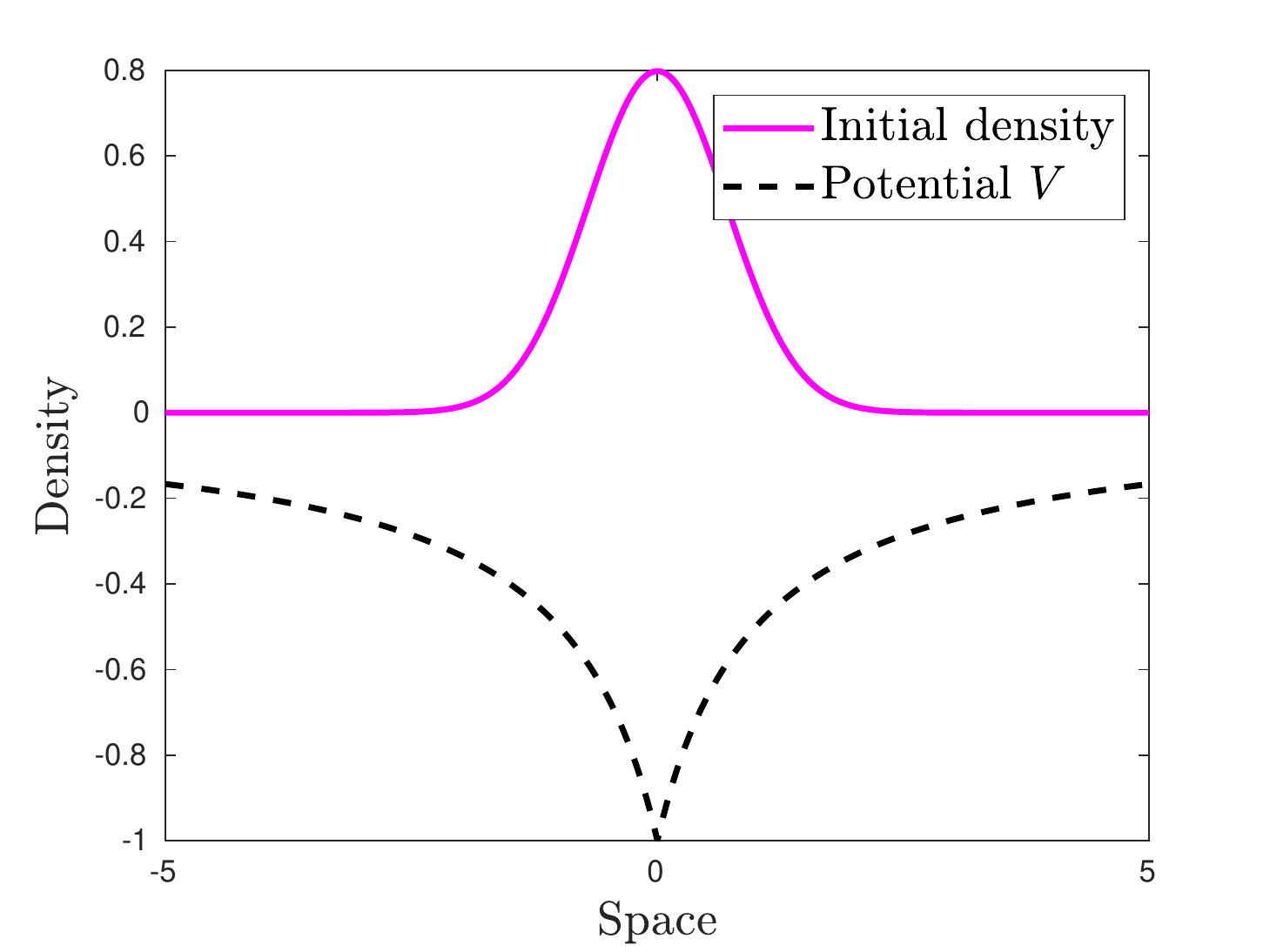}
\includegraphics[height=6cm,keepaspectratio]{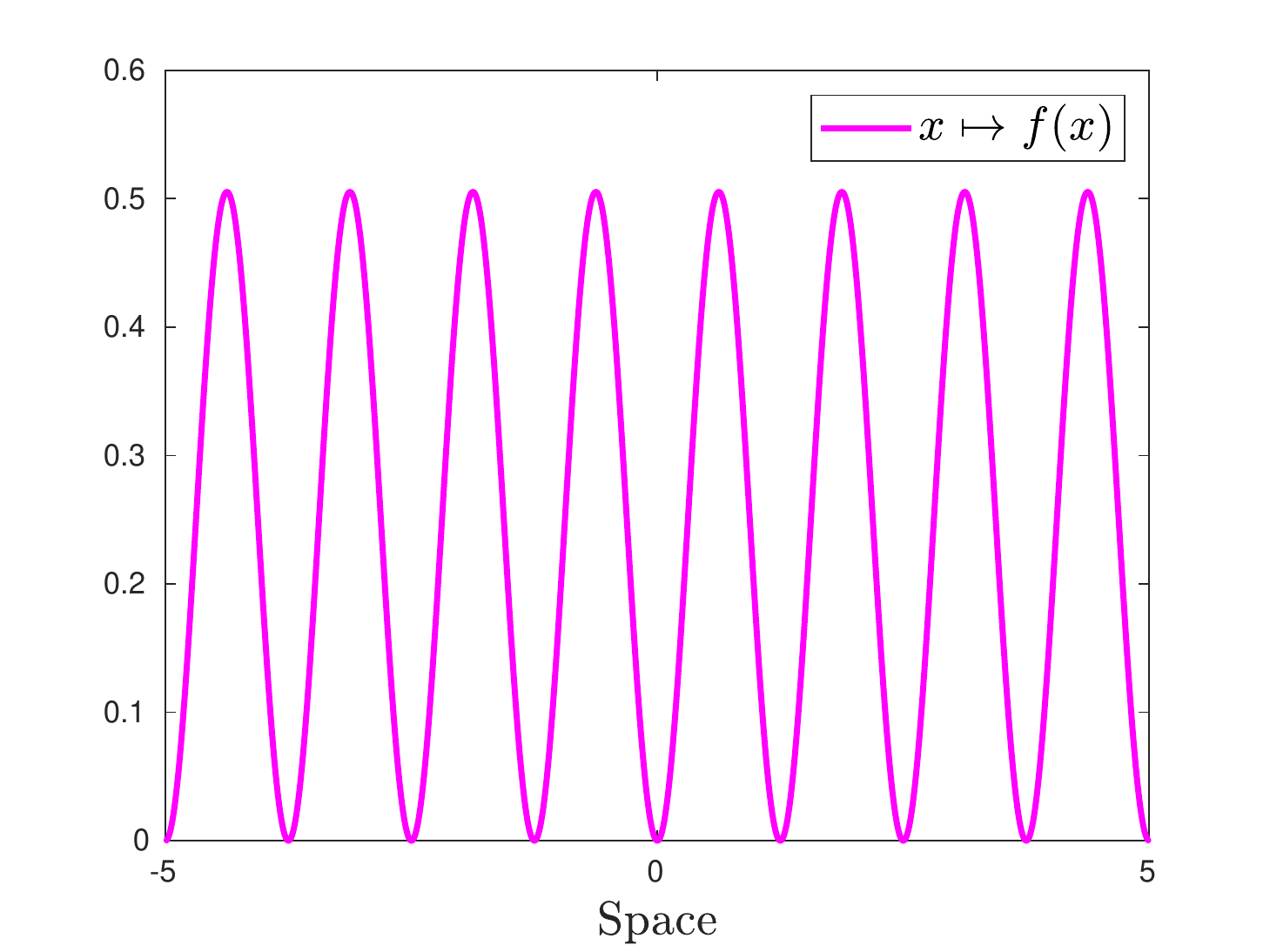}
\end{center}
\caption{{\bf Numerical Experiments 5.} (Left) Initial density and potential $V$. (Right) Graph of $f$.}
\label{figV2}
\end{figure}

\begin{figure}
\begin{center}
\includegraphics[height=6cm,keepaspectratio]{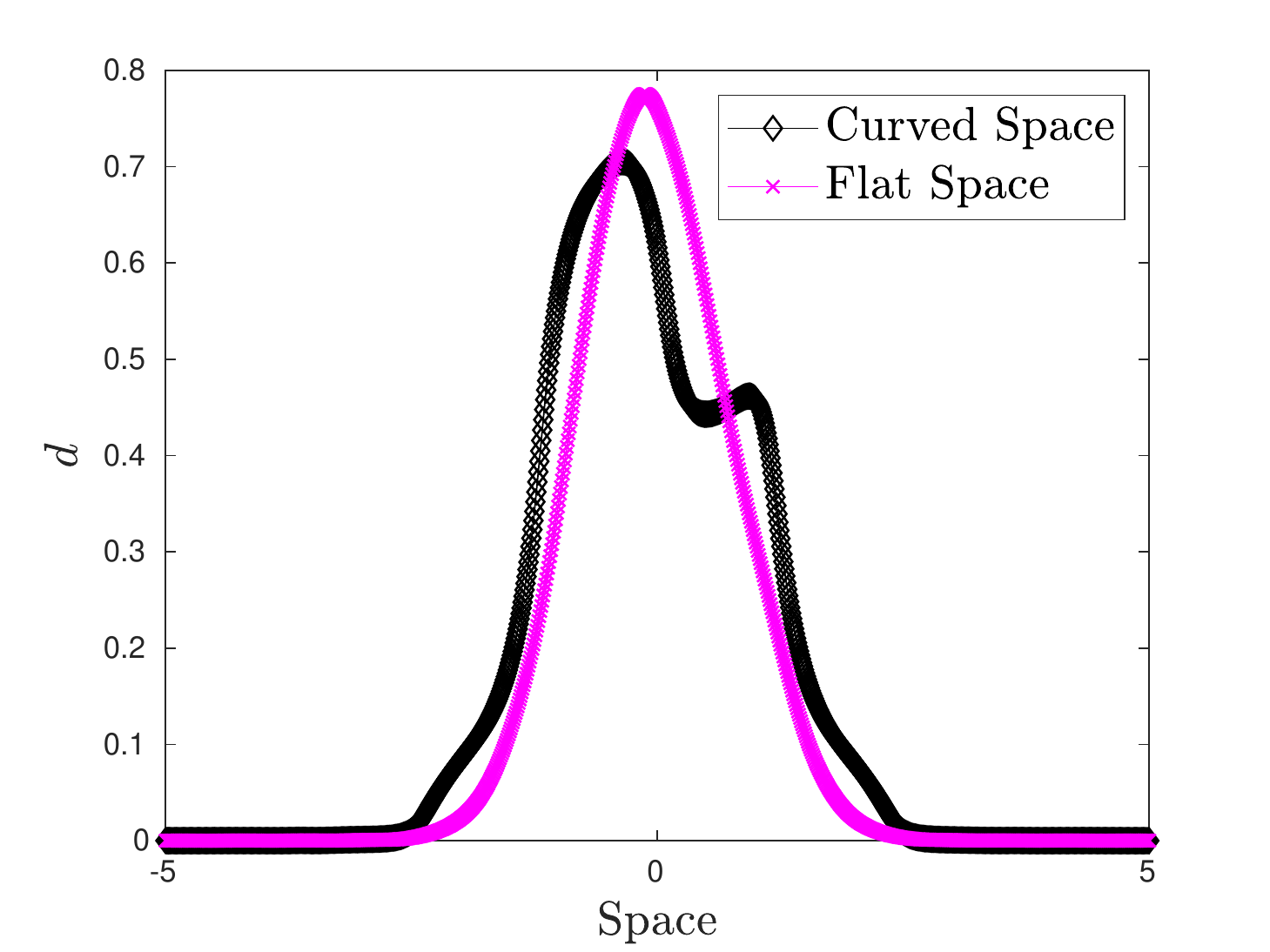}
\includegraphics[height=6cm,keepaspectratio]{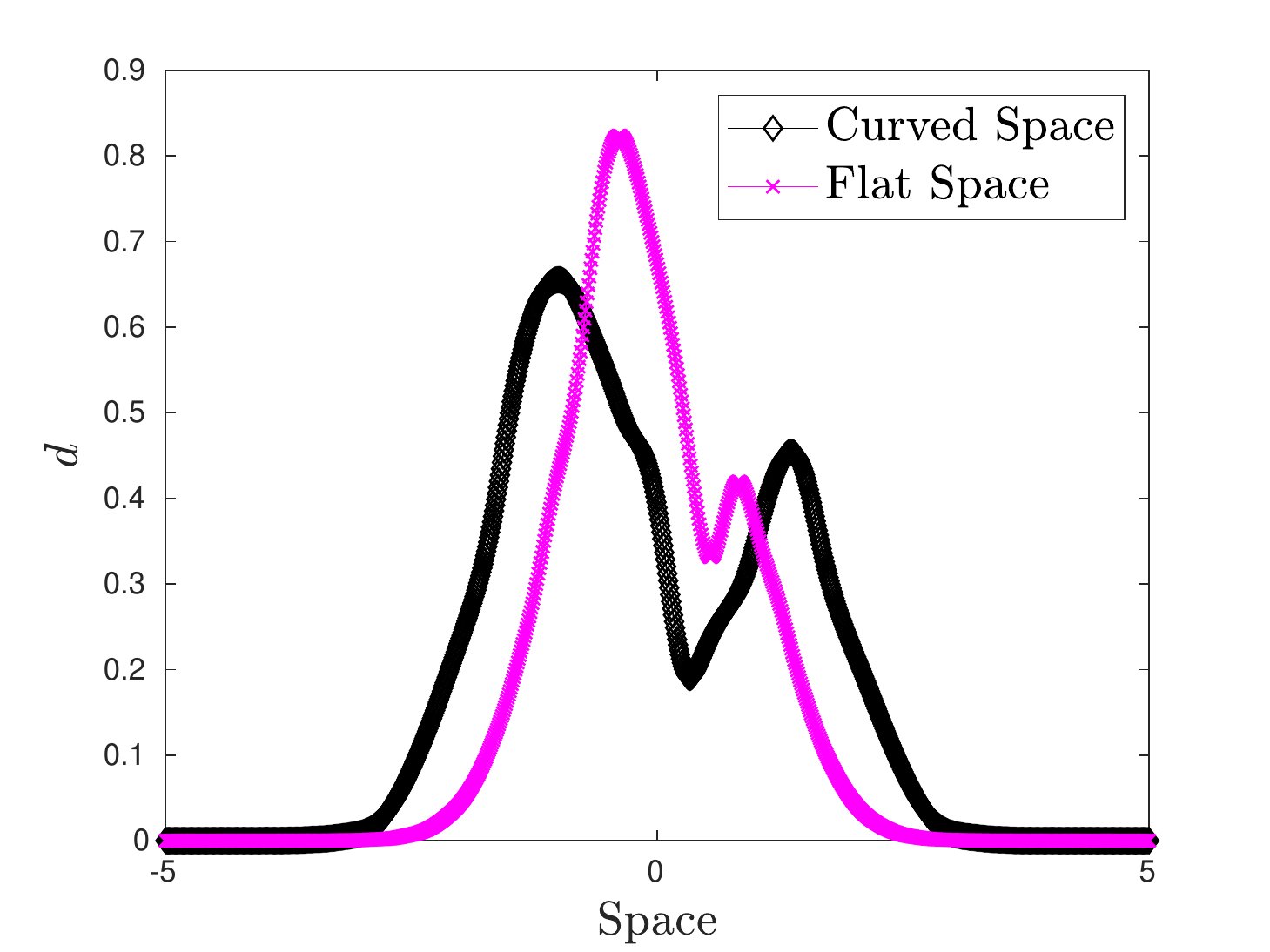}
\includegraphics[height=6cm,keepaspectratio]{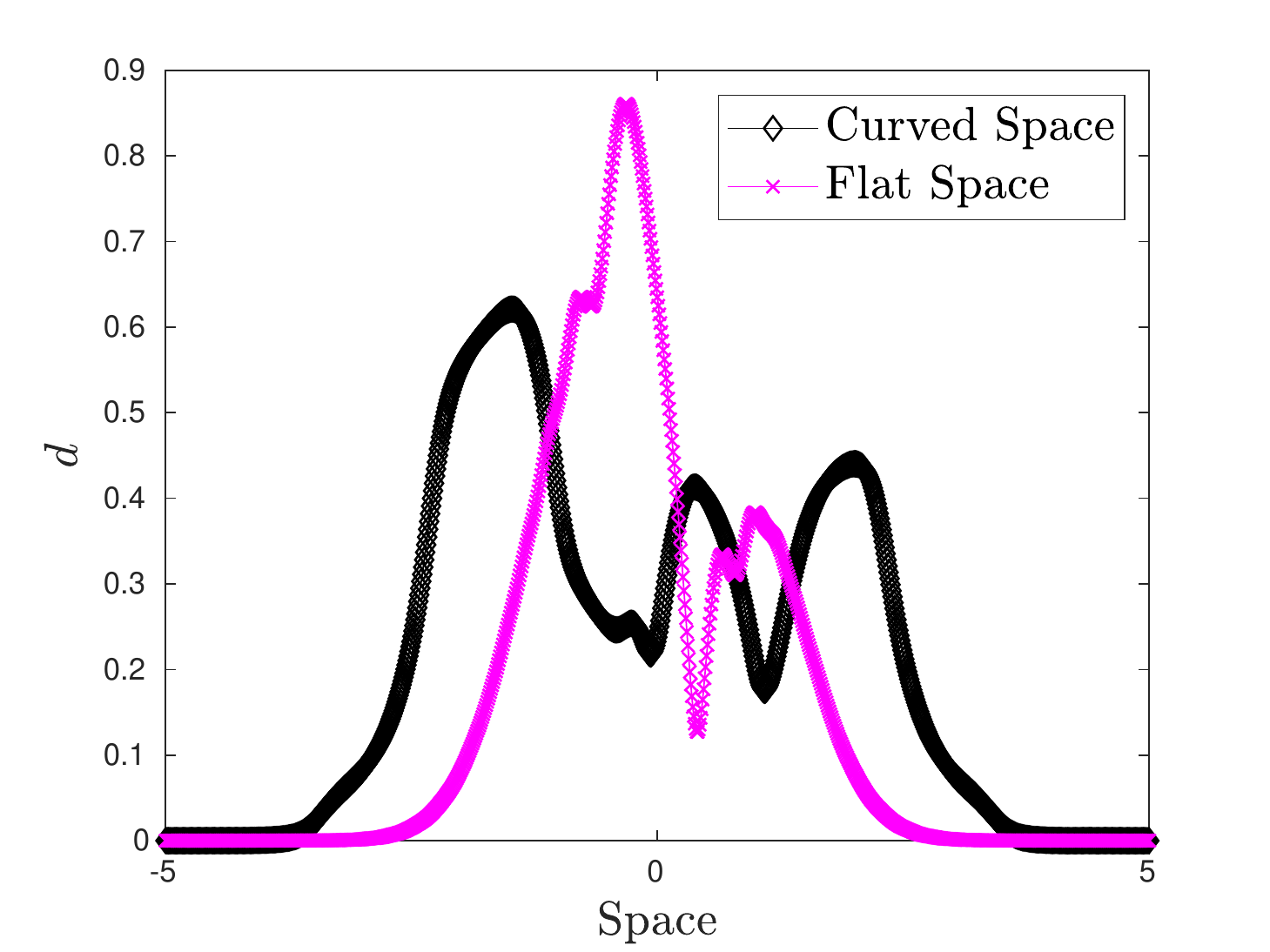}
\includegraphics[height=6cm,keepaspectratio]{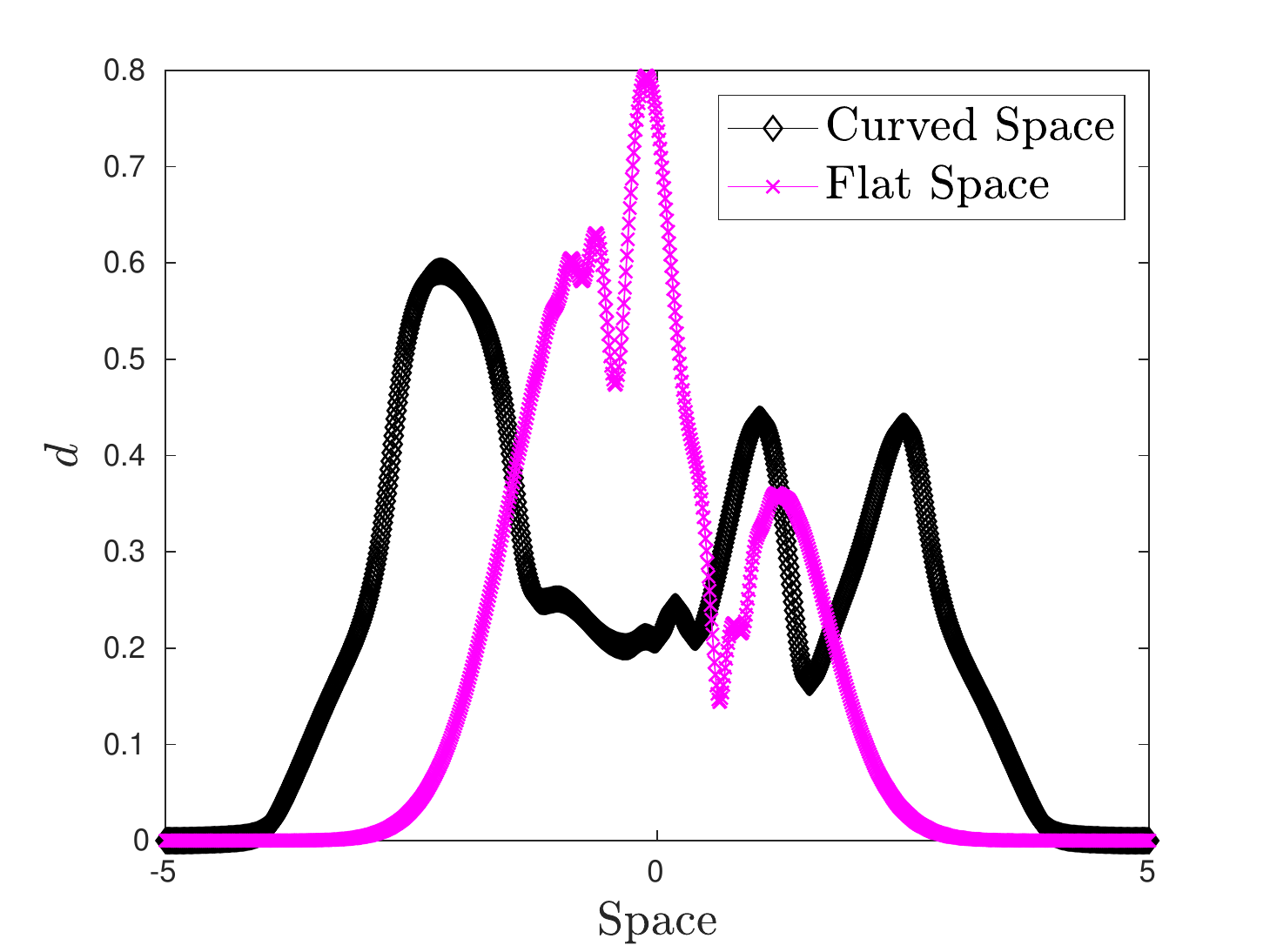}
\end{center}
\caption{{\bf Numerical Experiments 5.} From Top-Left to Bottom-Right: density in flat and curved spaces at time $t=0.2$, $t=0.4$, $t=0.6$ and $t=0.8$.}
\label{compDen2}
\end{figure}

\begin{figure}
\begin{center}
\includegraphics[height=6cm,keepaspectratio]{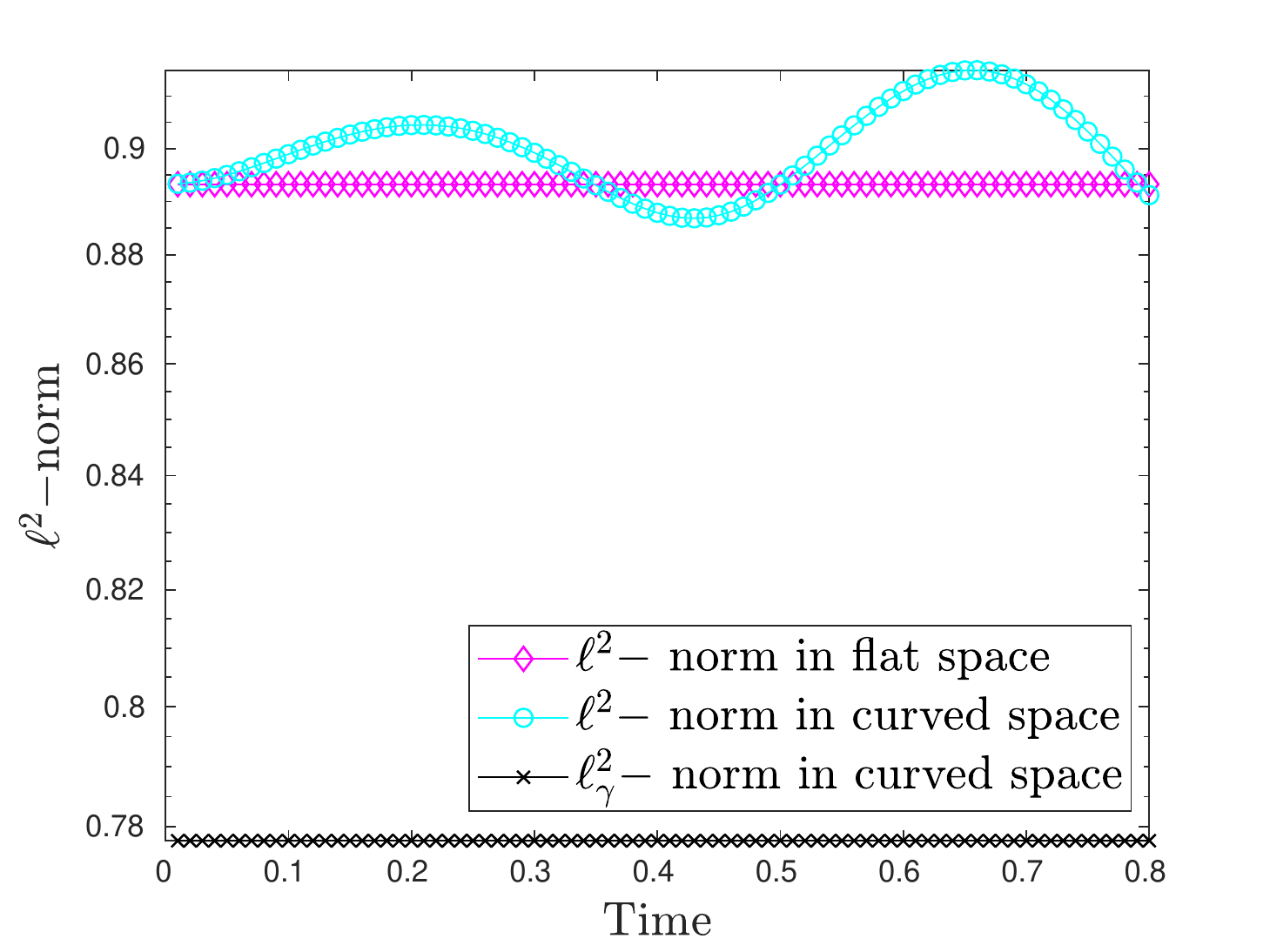}
\end{center}
\caption{{\bf Numerical Experiments 5.} $\ell^2-$norm of the solution as a function of time iterations in flat and curved space, and $\ell^2_{\gamma}$-norm of the solution as a function of time iterations in curved space.}
\label{l2norm}
\end{figure}

\noindent{\bf Numerical Experiments 5.} A more severe test is performed with different physical data. The computational domain is given by $[-5,5]$ and the initial data is
\begin{eqnarray*}
\psi(0,x) = (1,{\tt i})^T\beta e^{-\beta x^2/2}/\sqrt{4\pi } \, ,
\end{eqnarray*}
with $\beta=2$. Numerically, we take $\Delta t=10^{-2}$, and $h=10^{-2}$. Now, we take $V(x)=1/(|x|+1)$, $A_x=10x^2$ and $a_0=0.4$,  $\ell=10$, $k_0=5$.  The initial data, the potential, the functions $G$ and $f$ are reported in Fig. \ref{figV2}. We 
plot in Fig. \ref{compDen2}, the density $d_F(t,\cdot)$  in flat space, and the density $d_C(t,\cdot)$ in curved space at different times $t=0.2$, $t=0.4$, $t=0.6$ and $t=0.8$. In Fig. \ref{l2norm}, we report the $\ell^2-$norm of the solution as a function of time iterations, in flat and curved space, illustrating the $\ell^2-$stability (and $\ell^2-$norm preserving in flat space) of the proposed scheme.\\
The straining on graphene is enhanced compared to the above setting.\\

\begin{figure}
\begin{center}
\includegraphics[height=5.5cm,keepaspectratio]{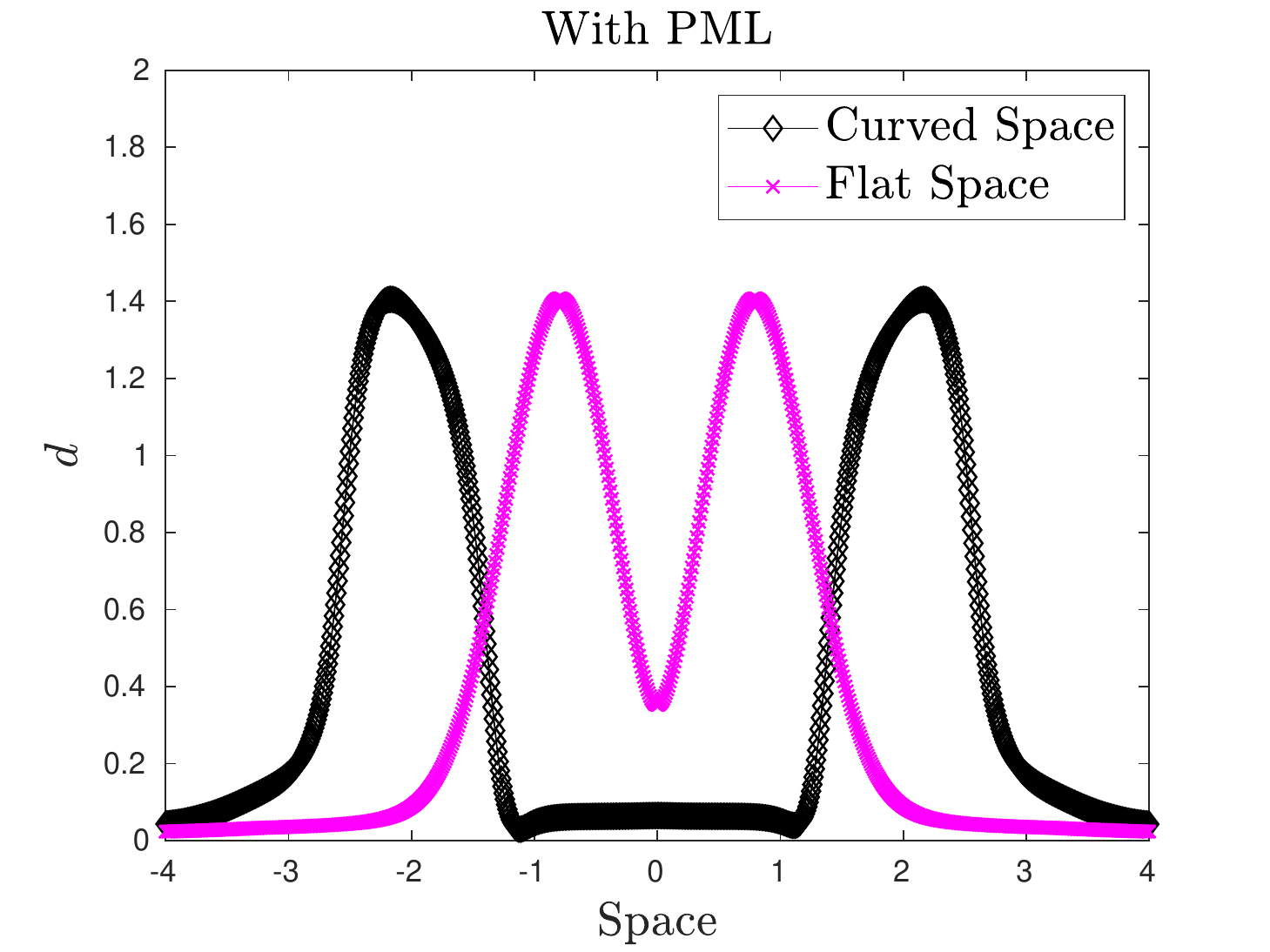}
\includegraphics[height=5.5cm,keepaspectratio]{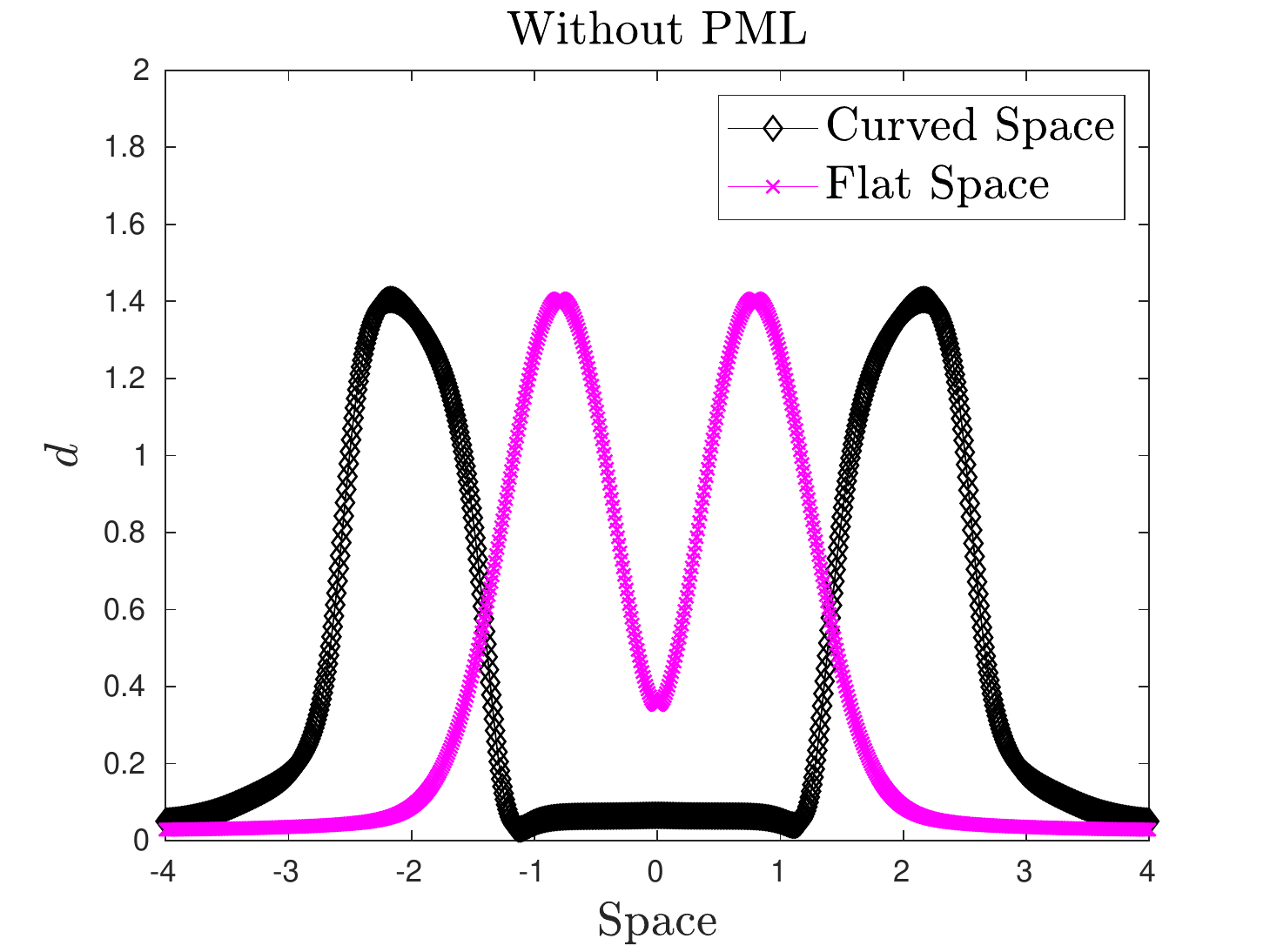}\\
\includegraphics[height=5.5cm,keepaspectratio]{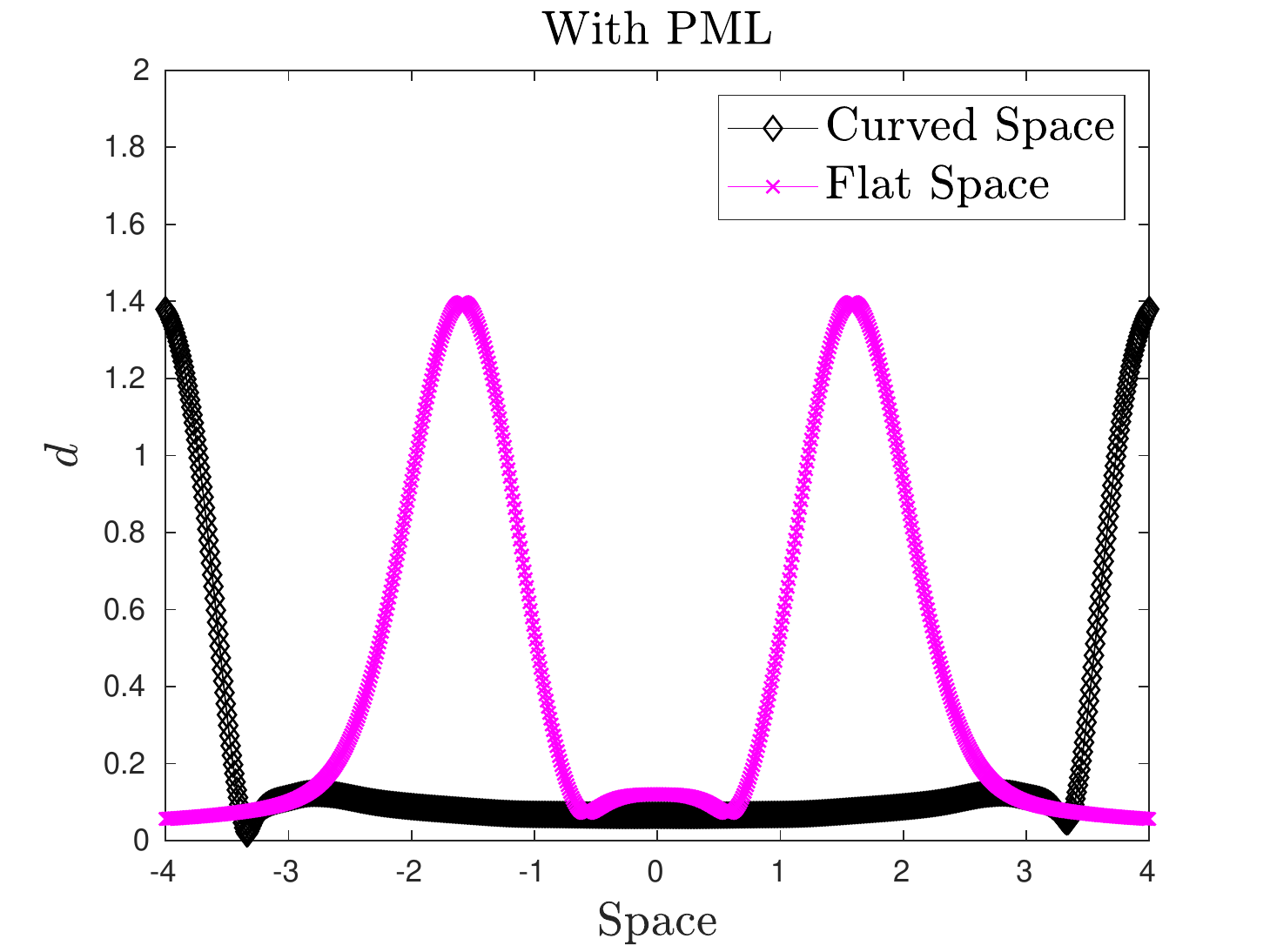}
\includegraphics[height=5.5cm,keepaspectratio]{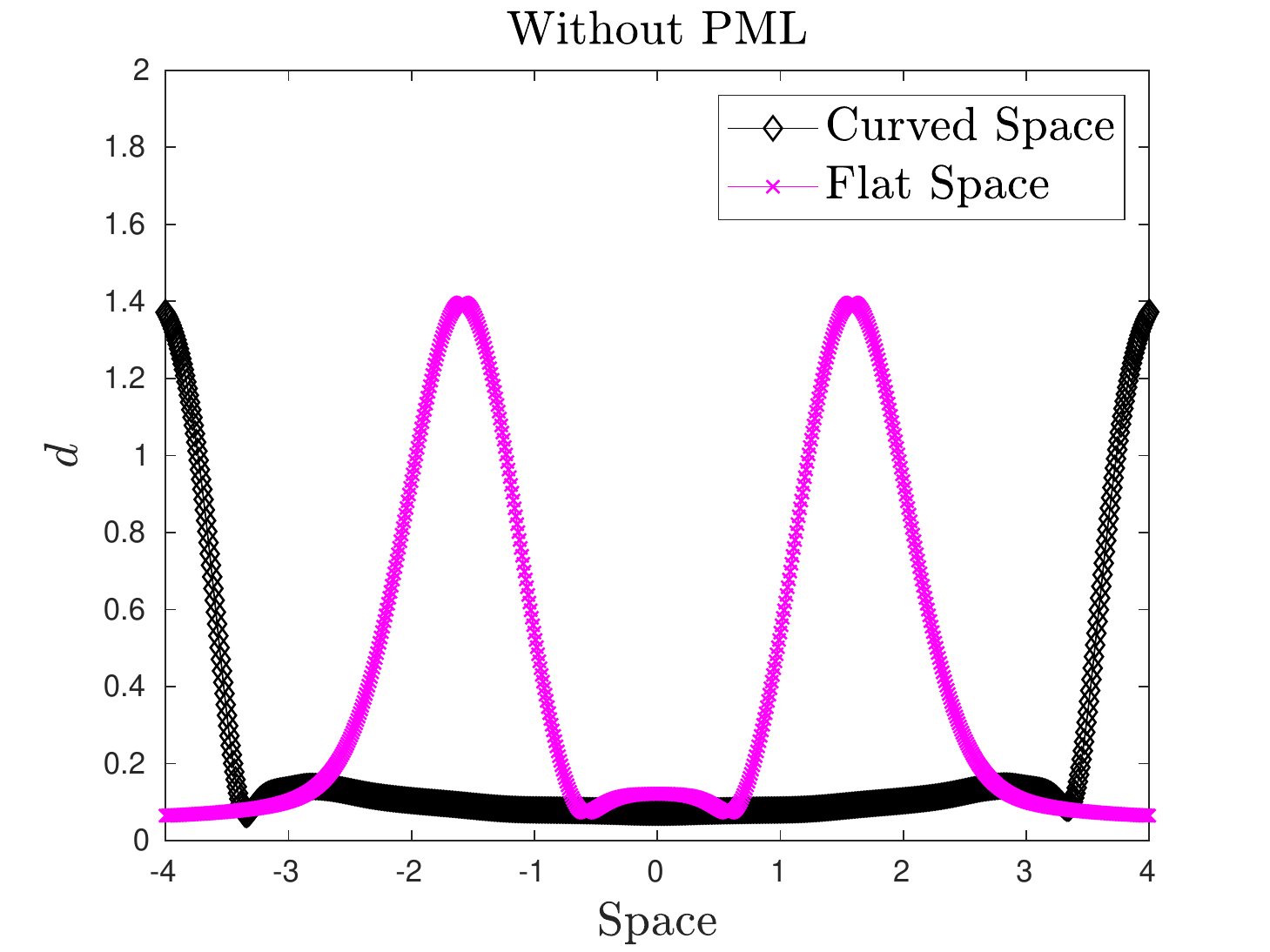}\
\includegraphics[height=5.5cm,keepaspectratio]{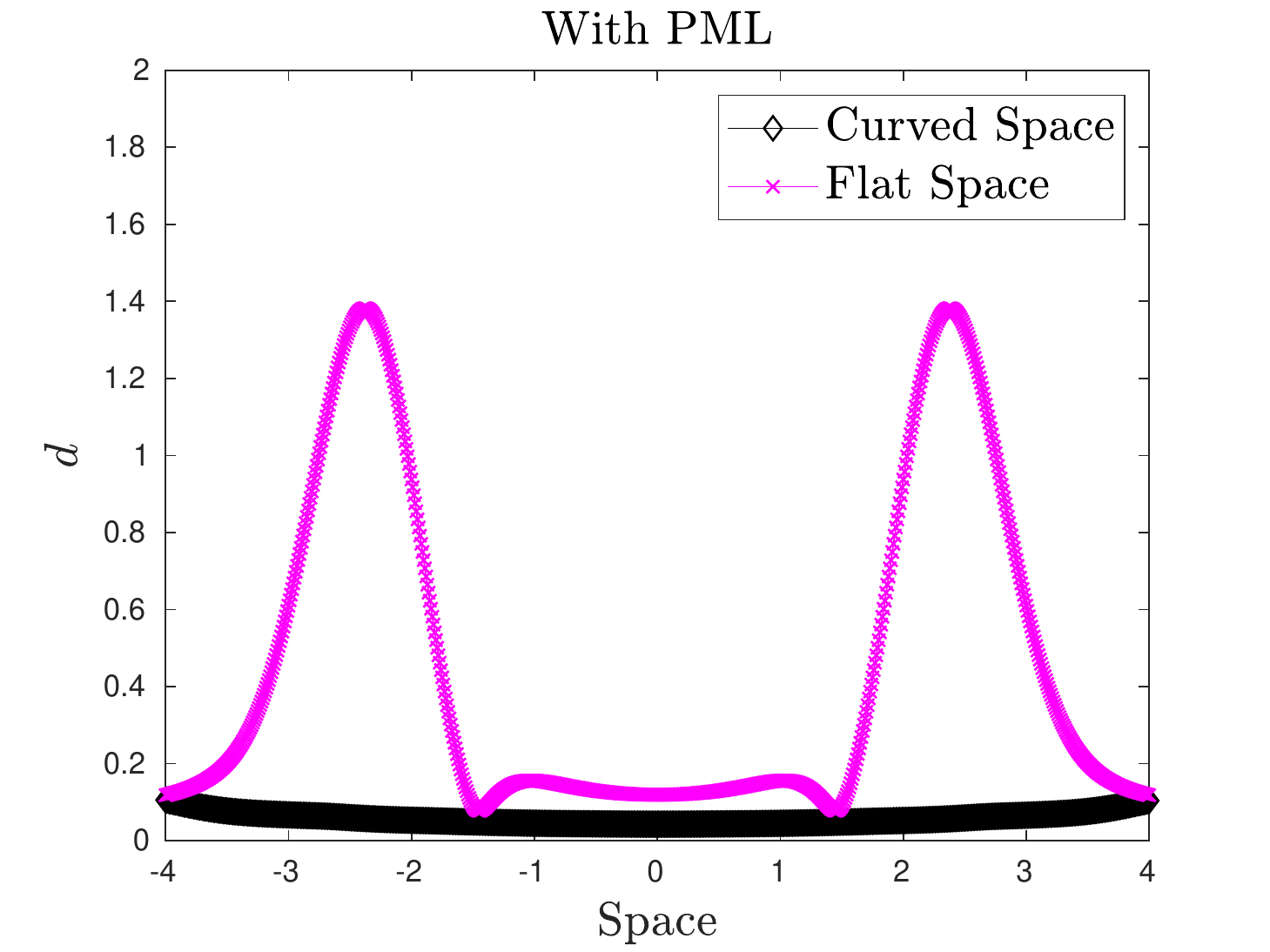}
\includegraphics[height=5.5cm,keepaspectratio]{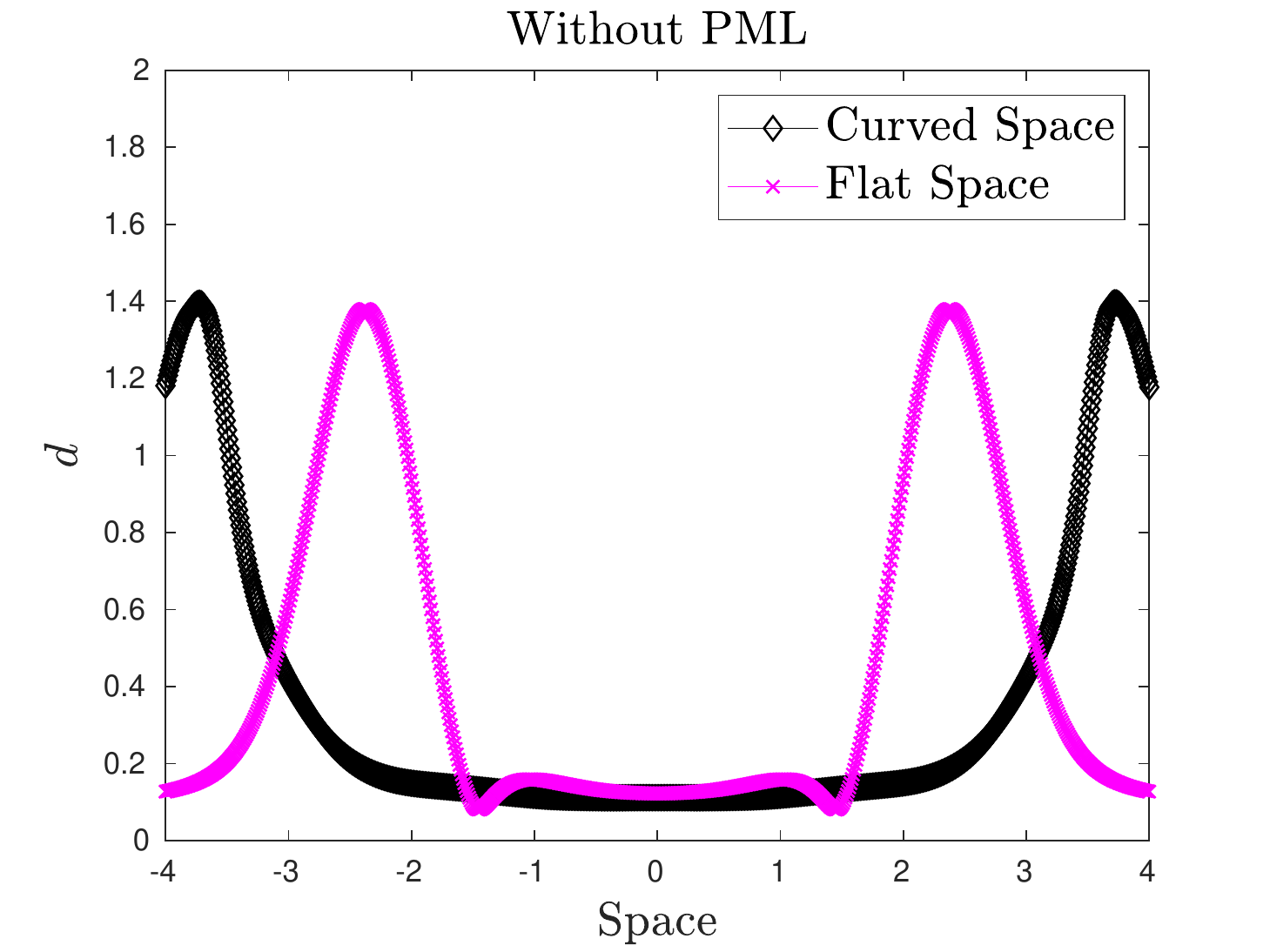}\\
\includegraphics[height=5.5cm,keepaspectratio]{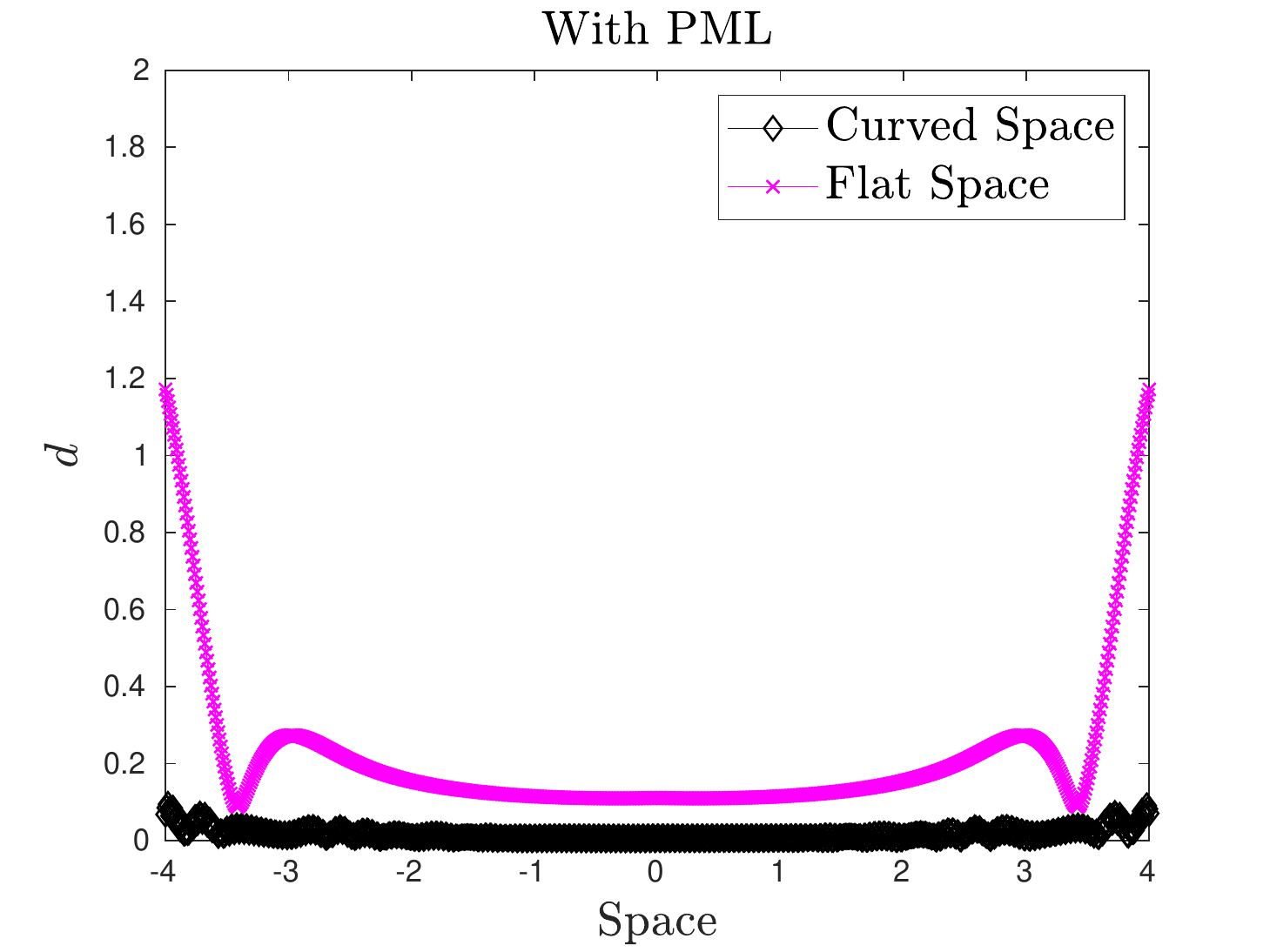}
\includegraphics[height=5.5cm,keepaspectratio]{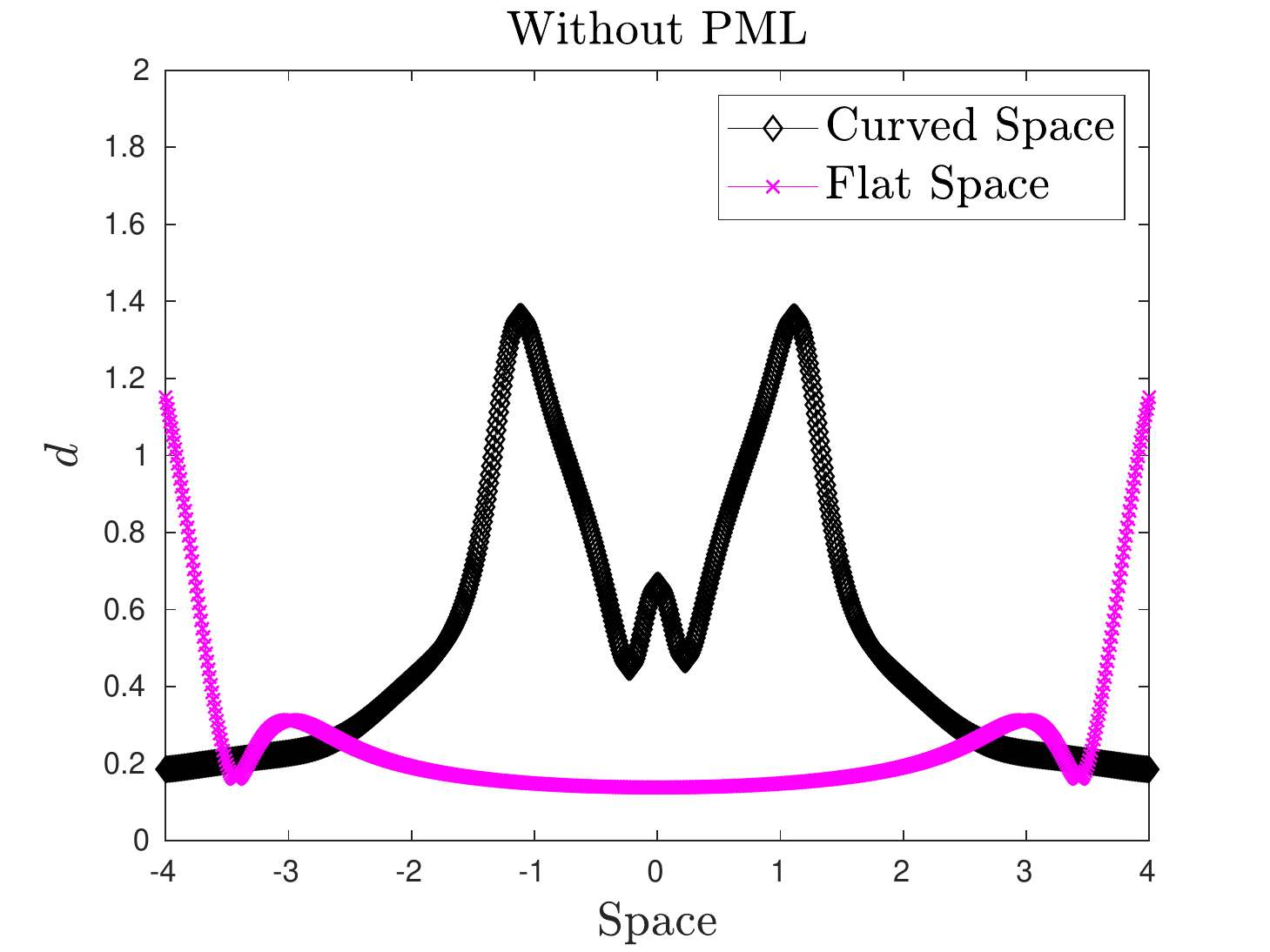}
\end{center}
\caption{{\bf Numerical Experiments 6.} Density in flat and curved spaces at time $t=0.75$, $t=1.5$, $t=2.25$ and $t=4$. (Left column) with PML. (Right column) without PML.}
\label{compDen2PML}
\end{figure}

\begin{figure}
\begin{center}
\includegraphics[height=6cm,keepaspectratio]{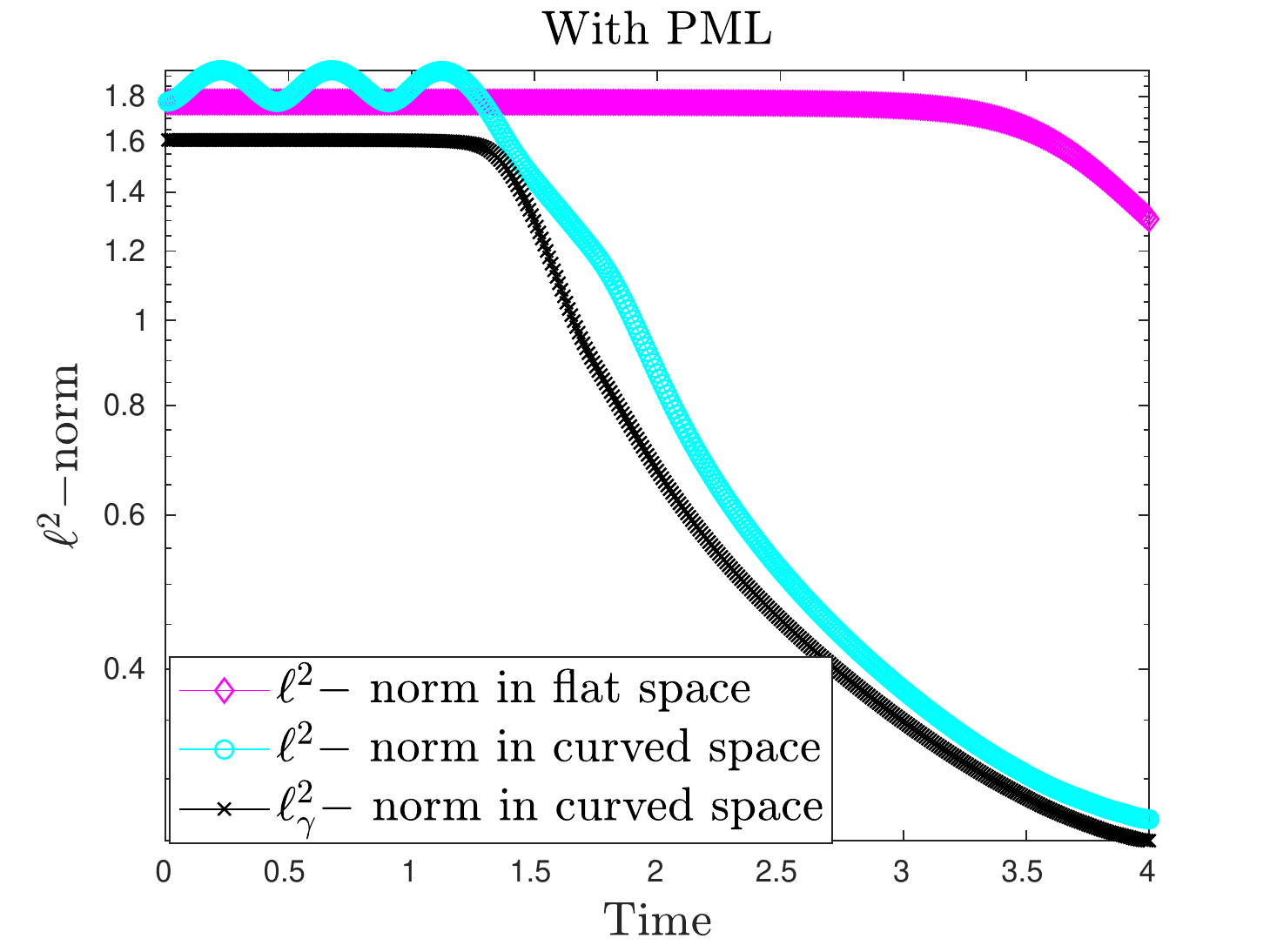}
\includegraphics[height=6cm,keepaspectratio]{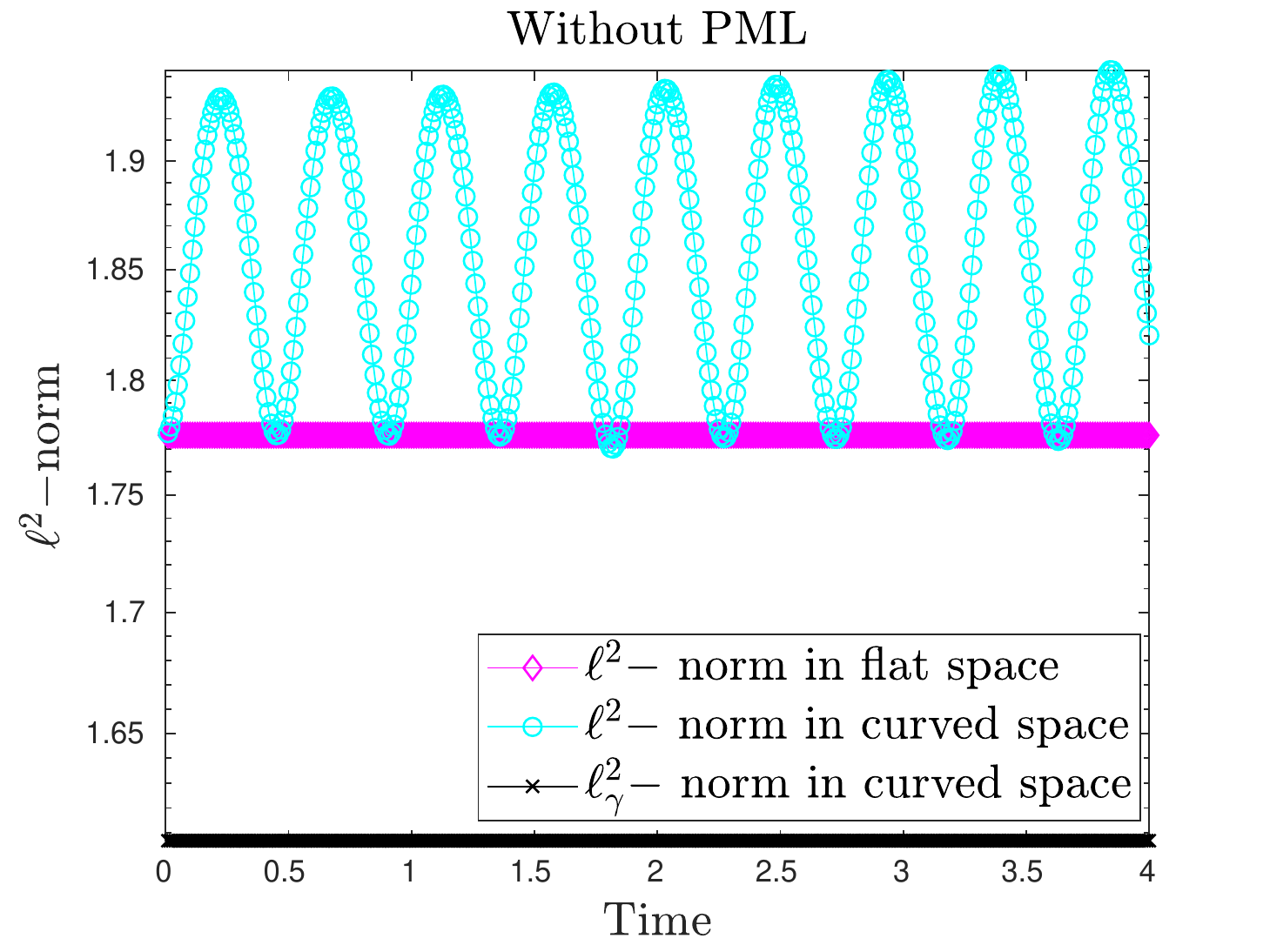}
\end{center}
\caption{{\bf Numerical Experiments 6.} $\ell^2-$norm of the solution as a function of time iterations in flat and curved space, and $\ell^2_{\gamma}$-norm of the solution as a function of time iterations in curved space. (Left) with PML. (Right) without PML.}
\label{l2normPML}
\end{figure}

\noindent{\bf Numerical Experiments 6.} The last experiment is dedicated to the use of PML in order to absorb the waves reaching the computational domain boundary. We compare the solution and its norm using and without using PML. The objective is to show that PML allow to circumvent the effect of periodic boundary conditions, by avoiding the propagation of non-physical waves entering back into the physical domain from one side to the other. We refer to \cite{jcp2019}, for a detailed study of PML for the Dirac equation using the pseudospectral method presented in this paper. The objective in this example is not to construct the best PML possible (this would require a fine study of absorbing functions and their parameters), but rather to show how efficient these PML can be. In the following, we consider absorbing functions of type I (with $\sigma_0=1$ and $\theta=0$, see Section \ref{sec:PML} for definitions), with a PML corresponding to $10\%$ of the overall computational domain. The computational domain is given by $[-4.5,4.5]$ and the initial data is
\begin{eqnarray*}
\psi(0,x) = (1,{\tt i})^T\beta e^{-\beta x^2/2}/\sqrt{4\pi } \, ,
\end{eqnarray*}
with $\beta=2$. Numerically, we take $\Delta t=10^{-2}$, and $h=10^{-2}$. Now, we take $V(x)=A_x(x)=0$, $m=0$, and $a_0=0.4$,  $\ell=5$, $k_0=2$.  We report in Fig. \ref{compDen2PML}, the density $d_F(t,\cdot)$  in flat space, and the density $d_C(t,\cdot)$ in curved space at different times $t=0.75$, $t=1.5$, $t=2.25$ and $t=4$, with (Left column) and without PML (Right column). In Fig. \ref{l2normPML}, we report the $\ell^2-$norm of the solution as a function of time iterations, in flat and curved space with and without PML. The latter illustrate the conservation of the norms, when using periodic boundary conditions without PML, and {\it vice-et-versa}.

This experiment shows that, although the proposed method naturally imposes perdiodic boundary conditions, its negative effects can be tackled thanks to Perfectly Matched Layers, which can very easily be implemented within the pseudospectral method.

\section{Conclusion}\label{sec:conclusion}
We have derived and analyzed simple pseudospectral computational methods for solving the Dirac equation in curved space  with perfectly matched layers at the computational domain boundary, and more generally for Dirac-like equations with non-constant coefficients. Interestingly, the proposed methods can easily be implemented from existing Fourier-based methods. Some
numerical one- and two-dimensional  experiments illustrating the properties of the numerical schemes were proposed. In a forthcoming paper, we will apply the developed methodology
to an extensive study of strained graphene.

\appendix

\section{Example of space-discretization in the Crank-Nicolson scheme}
\label{app:cn_2d}

As an example, let us consider a two-dimensional Dirac equation in  curved space, defined by the metric $d{\boldsymbol s}^2=e^{2\Phi({\boldsymbol x})}dt^2-e^{2\Psi({\boldsymbol x})}d{\boldsymbol x}^2$, such that $\Phi$ and $\Psi$ are two space-dependent functions. This leads to the following two-dimensional Dirac equation \cite{Koke}
\begin{eqnarray*}
	{\tt i}\partial_t\psi & = &  -{\tt i}e^{\Phi({\boldsymbol x})-\Psi({\boldsymbol x})}\Big(\sigma_x\Big(\partial_x + \cfrac{\partial_x\Phi({\boldsymbol x})}{2}\Big) + \sigma_y\Big(\partial_y + \cfrac{\partial_y\Phi({\boldsymbol x})}{2}\Big) \Big)\psi + e^{\Phi({\boldsymbol x})}\sigma_zm\psi \, .
\end{eqnarray*}
We can rewrite this equation in the form
\begin{eqnarray*}
	\left.
	\begin{array}{lcl}
		\partial_t\psi & = &  -e^{\Phi({\boldsymbol x})-\Psi({\boldsymbol x})}\sigma_x\partial_x \psi -e^{\Phi({\boldsymbol x})-\Psi({\boldsymbol x})}\sigma_y\partial_y \psi \\
		& & -\Big(e^{\Phi({\boldsymbol x})-\Psi({\boldsymbol x})}\big(\cfrac{\partial_x\Phi({\boldsymbol x})}{2}\sigma_x + \cfrac{\partial_y\Phi({\boldsymbol x})}{2}\sigma_y\big) -{\tt i}e^{\Phi({\boldsymbol x})}\sigma_zm\Big)\psi \, .
	\end{array}
	\right.
\end{eqnarray*}
We denote by ${\boldsymbol \psi}_{k}^n$ the approximate wavefunction at time $t_n$ and at ${\boldsymbol x}_{k}$. The numerical scheme reads:
\begin{eqnarray*}
	\left.
	\begin{array}{lcl}
		\cfrac{{\boldsymbol \psi}_{k}^{n+1} - {\boldsymbol \psi}_{k}^{n}}{\Delta t} & = &  -e^{\Phi({\boldsymbol x}_{k})-\Psi({\boldsymbol x}_{k})}\sigma_x[[ \partial_x]] {\boldsymbol \psi}^{n+1}_{k} -e^{\Phi({\boldsymbol x}_{k})-\Psi({\boldsymbol x}_{k})}\sigma_y[[ \partial_y ]] {\boldsymbol \psi}_{k}^{n+1} \\
		& & -\Big(e^{\Phi({\boldsymbol x}_{k})-\Psi({\boldsymbol x}_{k})}\big( \cfrac{\partial_x\Phi({\boldsymbol x}_{k})}{2}\sigma_x + \cfrac{\partial_y\Phi({\boldsymbol x}_{k})}{2}\sigma_y\big) -{\tt i}e^{\Phi({\boldsymbol x}_{k})}\sigma_zm\Big){\boldsymbol \psi}_{k}^{n+1}.
	\end{array}
	\right.
\end{eqnarray*}
Adding PML (see Section \ref{sec:PML}) to the numerical scheme formally leads to
\begin{eqnarray*}
	\left.
	\begin{array}{lcl}
		\cfrac{{\boldsymbol \psi}_{k}^{n+1} - {\boldsymbol \psi}_{k}^{n}}{\Delta t} & = &  -\cfrac{e^{\Phi({\boldsymbol x}_{k})-\Psi({\boldsymbol x}_{k})}}{S_x({\boldsymbol x}_{k})}\sigma_x[[ \partial_x]] {\boldsymbol \psi}^{n+1}_{k} -\cfrac{e^{\Phi({\boldsymbol x}_{k})-\Psi({\boldsymbol x}_{k})}}{S_y({\boldsymbol x}_{k})}\sigma_y[[ \partial_y ]] {\boldsymbol \psi}_{k}^{n+1} \\
		& & -\Big(e^{\Phi({\boldsymbol x}_{k})-\Psi({\boldsymbol x}_{k})}\big( \cfrac{\partial_x\Phi({\boldsymbol x}_{k})}{2}\sigma_x + \cfrac{\partial_y\Phi({\boldsymbol x}_{k})}{2}\sigma_y\big) -{\tt i}e^{\Phi({\boldsymbol x}_{k})}\sigma_zm\Big){\boldsymbol \psi}_{k}^{n+1} \, .
	\end{array}
	\right.
\end{eqnarray*}

\section{Representation of the spectral derivative}
\label{app:spec}

Another representation of the spectral derivative can be obtained by taking the Fourier transform on the RHS of \eqref{trick}. One obtains
\begin{eqnarray}\label{eq:spectral_rep}  
\big\{[[\partial_i]]\psi^n\big\}_{k}  =  
\cfrac{1}{N_i}
\sum_{p_{i}=-N_i/2}^{N_i/2-1}\sum_{k'_{i}=0}^{N_i-1}
{\tt i}\xi^{i}_{p_{i}}e^{{\tt i}\xi^{i}_{p_{i}}(x^{i}_{k_i} - x^{i}_{k'_i})} \psi^{n}_{k|k_{i} \rightarrow k'_{i}}.
%
%
\end{eqnarray}
This can be simplified further by noting that the sum on $(p_{i})_{i=1,2,3}$ can be performed explicitly. Therefore, the final result is that
\begin{eqnarray}\label{eq:spectral_rep_final}
\big\{[[\partial_i]]\psi^n\big\}_{k}  =  
\sum_{k'_{i}=0}^{N_i-1}
A^{i}_{k_{i} k'_{i}} \psi^{n}_{k|k_{i} \rightarrow k'_{i}} , 
%
%
\end{eqnarray}
where the differentiation matrices are given by 
\begin{eqnarray}
A^{i}_{k_{i} k'_{i}} &=& 
\cfrac{1}{N_i}
\sum_{p_{i}=-N_i/2}^{N_i/2-1}
{\tt i}\xi^{i}_{p_{i}}e^{{\tt i}\xi^{i}_{p_{i}}(x^{i}_{k_i} - x^{i}_{k'_i})} \\
&=& {\tt i}\frac{\pi}{N_{i}a_{i}} \left[ \frac{2e^{{\tt i}B^{i}_{k_{i} k'_{i}}} \sin \left( \frac{N_{i} B^{i}_{k_{i} k'_{i}}}{2}\right)}{\left(e^{{\tt i}B^{i}_{k_{i} k'_{i}}}-1\right)^{2}} + {\tt i} \frac{N_{i}\cos \left( \frac{N_{i} B^{i}_{k_{i} k'_{i}}}{2}\right)}{\left(e^{{\tt i}B^{i}_{k_{i} k'_{i}}}-1\right)}  \right],
\end{eqnarray}
where $B^{i}_{k_{i} k'_{i}}:= \frac{\pi}{a_{i}} (x_{k_{i}} - x_{k'_{i}})$.

\section{Explicit construction of the matrix in the Crank-Nicolson scheme}\label{APXA}

The second step, given explicitly by
\begin{eqnarray}
\label{eq:step2_cn}
G_{k}^{n+1/2} \psi^{n^*}_{k} & = & \widetilde{G}_{k}^{n+1/2}  \psi^{n+1/2}_{k},
\end{eqnarray}
can be written as a linear system of equations, by using the discrete pseudospectral representation of the derivative. For the LHS of \eqref{eq:step2_cn}, we use the representation given in Eqs. \eqref{trick} and \eqref{eq:spectral_rep_final} to obtain
\begin{eqnarray*}
G_{k}^{n+1/2} \psi^{n^*}_{k} & = & \psi_{k}^{n^*}  + 
\cfrac{\Delta t}{2} \sum_{i=1,2,3}  
\cfrac{\alpha_{k}^{i}}{S_{k}^{i}} \sum_{k'_{i}=0}^{N_{i}-1}
A^{i}_{k_{i},k'_{i}}
\psi^{n^{*}}_{k|k_{i} \rightarrow k'_{i}}  \, . 
%
%
%
\end{eqnarray*}
Re-arranging the sums and introducing Kronecker's symbol, the last expression can be written in the form of
\begin{eqnarray}
	\label{eq:matrix_vec}
G_{k}^{n+1/2} \psi^{n^*}_{k} & = & 
	  \sum_{k'=0}^{N-1}
  \mathcal{G}^{n}_{k k'}   \psi^{n^{*}}_{k'} ,
\end{eqnarray}
where we define the matrix representation of $G_{k}^{n+1/2}$ as
\begin{eqnarray}
\label{eq:G_mat}
\mathcal{G}^{n}_{kk'}&:=& \delta_{kk'} 
+\cfrac{\Delta t}{2} \sum_{i=1,2,3}  
\cfrac{\alpha_{k}^{i}}{S_{k}^{i}}
A^{i}_{k_{i},k'_{i}}
\delta_{kk'|k_{i}=k'_{i}} .
%
\end{eqnarray}
Equation \eqref{eq:matrix_vec} is just the matrix-vector product with a matrix defined in \eqref{eq:G_mat}.

The RHS of Eq. \eqref{eq:step2_cn} is simpler because it contains the initial data, which is a known quantity. Therefore, it is possible to use the spectral representation of the derivative in Fourier space as in Eq. \eqref{trick}. This yields
\begin{eqnarray*}
 \widetilde{G}_{k}^{n+1/2}  \psi^{n+1/2}_{k}	
&:=& \mathcal{H}^{n}_{k} ,\\
&=&
\psi_{k}^{n+1/2}  -
\cfrac{\Delta t}{2} \sum_{i=1,2,3}
\cfrac{\alpha_{k}^{i}}{S^{i}_{k}N_i}\sum_{p_{i}=-N_i/2}^{N_i/2-1}{\tt i}\xi^{i}_{p_{i}}\widetilde{\psi}^{n}_{k|k_{i} \rightarrow p_{i}}e^{{\tt i}\xi^{i}_{p_{i}}(x^{i}_{k_i}+a_i)}. 
%
%
%
\end{eqnarray*}
With these results, the linear system has the form
\begin{eqnarray}\label{step2}
\mathcal{G}^{n}\boldsymbol{\psi}^{n^*}  = \boldsymbol{\mathcal{H}}^{n} ,
\end{eqnarray}
where $\mathcal{G}^{n}$ is the matrix defined in \eqref{eq:G_mat} and the bold symbols represent vectors in real space, with components $\boldsymbol{V} = [V_{0,0,0},V_{1,0,0},\cdots,V_{N_{1}-1,N_{2}-1,N_{3}-1}]^{T}$. Solving this linear system yields  $\psi^{n^{*}}_{k}$, the value of the wave function after the second steps of the Crank-Nicolson scheme. The main problem with this approach is the evaluation of the matrix $\mathcal{G}^{n}$, which is not efficient (the computational complexity is $O(N^{2})$), and the storing of the same matrix which can be problematic since it requires $O(N^{2})$ of memory.

   \bigskip
   \medskip 

\bibliographystyle{unsrt}
\bibliography{biblio}

\end{document}